\documentclass[reqno, 11pt, a4paper]{amsart}

\usepackage{latexsym,ifthen,xspace}
\usepackage{amsmath,amssymb,amsthm}
\usepackage{enumerate, calc}
\usepackage[colorlinks=true, pdfstartview=FitV, linkcolor=blue,
  citecolor=blue, urlcolor=blue, pagebackref=false]{hyperref}
\usepackage[text={33pc,605pt},centering]{geometry}    
\usepackage{graphicx,color}

\newtheorem{theorem}{Theorem}[section]
\newtheorem{lemma}[theorem]{Lemma}
\newtheorem{prop}[theorem]{Proposition}
\newtheorem{assumption}[theorem]{Assumption}
\newtheorem{corro}[theorem]{Corollary}

\theoremstyle{definition}

\theoremstyle{remark}
\newtheorem{remark}[theorem]{Remark}

\numberwithin{equation}{section}

\DeclareMathAlphabet{\mathsl}{OT1}{cmss}{m}{sl}
\SetMathAlphabet{\mathsl}{bold}{OT1}{cmss}{bx}{sl}

\newcommand{\al}{\ensuremath{\alpha}}
\newcommand{\be}{\ensuremath{\beta}}
\newcommand{\ga}{\ensuremath{\gamma}}
\newcommand{\de}{\ensuremath{\delta}}

\newcommand{\ze}{\ensuremath{\zeta}}

\renewcommand{\th}{\ensuremath{\theta}}

\newcommand{\ka}{\ensuremath{\kappa}}
\newcommand{\la}{\ensuremath{\lambda}}

\newcommand{\si}{\ensuremath{\sigma}}

\newcommand{\om}{\ensuremath{\omega}}

\newcommand{\ve}{\ensuremath{\varepsilon}}

\newcommand{\vr}{\ensuremath{\varrho}}

\newcommand{\Si}{\ensuremath{\Sigma}}

\newcommand{\Om}{\ensuremath{\Omega}}
\newcommand{\cB}{\ensuremath{\mathcal B}}

\newcommand{\cE}{\ensuremath{\mathcal E}}
\newcommand{\cF}{\ensuremath{\mathcal F}}

\newcommand{\cH}{\ensuremath{\mathcal H}}

\newcommand{\cL}{\ensuremath{\mathcal L}}

\newcommand{\cO}{\ensuremath{\mathcal O}}

\newcommand{\bbN}{\ensuremath{\mathbb N}}

\newcommand{\bbR}{\ensuremath{\mathbb R}}

\newcommand{\bbZ}{\ensuremath{\mathbb Z}} 
\definecolor {orange} {rgb} {0.569, 0.259, 0.0}

\newcommand{\me}{\ensuremath{\mathrm{e}}}

\newcommand{\md}{\ensuremath{\mathrm{d}}}

\newcommand{\scpr}[3]{%
  \ensuremath{%
    \left\langle
      #1, #2
    \right\rangle_{\raisebox{-0ex}{$\scriptstyle \ell^{\raisebox{.1ex}{$\scriptscriptstyle 2$}} (#3)$}}
  }
}
\newcommand{\norm}[3]{%
  \ensuremath{%
    \big\lVert
      #1
    \big\rVert_{\raisebox{-.0ex}{$\scriptstyle \ell^{\raisebox{.2ex}{$\scriptscriptstyle #2$}} (#3)$}}
  }
}
\newcommand{\Norm}[2]{%
  \ensuremath{%
    \big\lVert
      #1
    \big\rVert_{\raisebox{-.0ex}{$\scriptstyle #2$}}
  }
}

\DeclareMathOperator{\mean}{\mathbb{E}}
\DeclareMathOperator{\Mean}{\mathrm{E}}
\DeclareMathOperator{\prob}{\mathbb{P}}
\DeclareMathOperator{\Prob}{\mathrm{P}}

\DeclareMathOperator{\supp}{\mathrm{supp}}

\DeclareMathOperator{\osr}{\mathrm{osr}}
\DeclareMathOperator{\osc}{\mathrm{osc}}
\DeclareMathOperator{\diam}{\mathrm{diam}}

\newcommand{\ldef}{\ensuremath{\mathrel{\mathop:}=}}

\def\indicator{{\mathchoice {1\mskip-4mu\mathrm l}%
{1\mskip-4mu\mathrm l}{1\mskip-4.5mu\mathrm l}%
{1\mskip-5mu\mathrm l}}}

\begin{document}

\title[Harnack inequalities on weighted graphs]{Harnack inequalities on weighted graphs and some applications to the random conductance model}


\author{Sebastian Andres}
\address{Rheinische Friedrich-Wilhelms Universit\"at Bonn}
\curraddr{Endenicher Allee 60, 53115 Bonn}
\email{andres@iam.uni-bonn.de}
\thanks{}

\author{Jean-Dominique Deuschel}
\address{Technische Universit\"at Berlin}
\curraddr{Strasse des 17. Juni 136, 10623 Berlin}
\email{deuschel@math.tu-berlin.de}
\thanks{}

\author{Martin Slowik}
\address{Technische Universit\"at Berlin}
\curraddr{Strasse des 17. Juni 136, 10623 Berlin}
\email{slowik@math.tu-berlin.de}
\thanks{}



\date{\today}

\dedicatory{}

\begin{abstract}
We establish elliptic and parabolic Harnack inequalities on graphs with unbounded weights. As an application we prove a local limit theorem for a continuous time random walk $X$ in an environment of ergodic random conductances taking values in $(0, \infty)$ satisfying some moment conditions.
\end{abstract}

\maketitle

\tableofcontents

\section{Introduction}
Consider the operator $L^a$ in divergence form,
\begin{align*}
 \big(L^a f\big)(x)
  \;=\;
  \sum_{i,j=1}^d
  \partial_{x_i}\Big(a_{ij}(\cdot)\, \partial_{x_j}f(\cdot)\Big)(x),
\end{align*}
acting on functions on $\bbR^d$, where the symmetric positive definite matrix $a = (a_{ij}(x))$ is bounded, measurable and uniformly elliptic, i.e.\ for some $0 < \la_1 \leq \la_2 < \infty$
\begin{align} \label{eq:ell_cont}
  \la_1\, |\xi|^2
  \;\leq\;
  \sum_{i, j = 1}^d a_{ij}(x)\, \xi_i\, \xi_j
  \;\leq\;
  \la_2\, |\xi|^2,
  \qquad \forall\, x, \xi \in \bbR^d.
\end{align}
Then, a celebrated result of Moser \cite{Mo61} states that the elliptic Harnack inequality (EHI) holds, that is there exists a constant $C_{\mathrm{EH}} \equiv C_{\mathrm{EH}}(\la_1, \la_2) < \infty$ such that for any positive function $u > 0$ which is $L^a$ harmonic on the ball $B(x_0, r)$, i.e.\ $(L^a u)(x)=0$ for all $x \in B(x_0,r)$, we have
\begin{align*}
  \max_{B(x_0,r/2)} u \;\leq\; C_{\mathrm{EH}}\, \min_{B(x_0,r/2)} u.
\end{align*}
A few years later, Moser also proved that the parabolic Harnack inequality (PHI) holds (see \cite{Mo64}), that is there exists a constant $C_{\mathrm{PH}} \equiv C_{\mathrm{PH}}(\la_1, \la_2) < \infty$ such that for any  positive  caloric function $u > 0$, i.e.\ any weak solution of the heat equation $(\partial_t-L^a)u(t,x)=0$ on $[t_0, t_0 + r^2] \times B(x_0,r)$, we have
\begin{align*}
  \max_{(t,x) \in Q_{-}} u(t,x)
  \;\leq\;
  C_{\mathrm{PH}}\, \min_{(t,x) \in Q_+} u(t,x),
\end{align*}
where $Q_- = \big[t_0 + \frac{1}{4} r^2, t_0 + \frac{1}{2} r^2 \big] \times B\big(x_0, \frac{1}{2} r \big)$ and $Q_+ = \big[t_0 + \frac{3}{4} r^2, t_0 + r^2 \big] \times B\big(x_0, \frac{1}{2} r\big)$.

Both EHI and PHI have had a big influence on PDE theory and differential geometry, in particular they play a prominent role in elliptic regularity theory. One application of Harnack inequalities -- observed by Nash in \cite{Na58} -- is to show H\"older regularity for solutions of elliptic or parabolic equations.  For instance, the PHI implies the H\"older continuity of  the fundamental solution $p^a_t(x,y)$ of the parabolic equation
\begin{align*}
  \big( \partial_t - L^a \big) p^a_t(x,y)
  \;=\;
  0,
  \qquad p^a_0(x,y) = \delta_x(y).
\end{align*}
Furthermore, the PHI is also one of the key tools in Aronson's proof \cite{Ar67} of Gaussian type bounds for $p^a_t(x,y)$.

The remarkable fact about these results is that they only rely on the uniform ellipticity and not on the regularity of the matrix $a$. We refer to the monograph \cite{SC02} for more details on this topic.

\subsection{The model}
In this paper we will be dealing with a discretised version of the operator $L^a$. We consider an infinite, connected, locally finite graph $G = (V, E)$ with vertex set $V$ and edge set $E$.  We will write $x \sim y$ if $\{x,y\} \in E$ and $e \sim e'$ for $e, e' \in E$ if the edges $e$ and $e'$ have a common vertex.  A path of length $n$ between $x$ and $y$ in $G$ is a sequence $\{x_i : i = 0, \ldots , n\}$ with the property that $x_0 = x$, $x_n = y$ and $x_i \sim x_{i+1}$.  Let $d$ be the natural graph distance on $G$, i.e. $d(x,y)$ is the minimal length of a path between $x$ and $y$.  We denote by $B(x,r)$ the closed ball with center $x$ and radius $r$, i.e. $B(x,r) \ldef \{y \in V \mid d(x,y) \leq r\}$.

The graph is endowed with the counting measure, i.e.\ the measure of $A \subset V$ is simply the number $|A|$ of elements in $A$. For any non-empty, finite $A \subset V$ and $p \in [1, \infty)$, we introduce space-averaged $\ell^{p}$-norms on functions $f\!: A \to \bbR$ by the usual formula
\begin{align*}
  \Norm{f}{p,A}
  \;\ldef\;
  \Bigg(
    \frac{1}{|A|}\; \sum_{x \in A}\, |f(x)|^p
  \Bigg)^{\!\!\frac{1}{p}}
  \qquad \text{and} \qquad
  \Norm{f}{\infty, A} \;\ldef\; \max_{x\in A} |f(x)|.
\end{align*}
For a given set $B \subset V$, we define the \emph{relative} internal boundary of $A \subset B$ by
\begin{align*}
  \partial_B A
  \;\ldef\;
  \big\{
    x \in A \;\big|\;
    \exists\, y \in B \setminus A\; \text{ s.th. }\; \{x,y\} \in E
  \big\}
\end{align*}
and we simply write $\partial A$ instead of $\partial_V A$. Throughout the paper we will make the following assumption on $G$.
\begin{assumption}\label{ass:graph}
  For some $d\geq 2$ the graph $G$ satisfies the following conditions:
  \begin{enumerate}[(i)]
    \item volume regularity of order $d$, that is there exists $C_{\mathrm{reg}} \in(0, \infty)$ such that
      \begin{align}\label{eq:ass:vd}
       C^{-1}_{\mathrm{reg}}\, r^d  \;\leq\; |B(x,r)| \;\leq\; C_{\mathrm{reg}}\, r^d,
        \qquad \forall\, x \in V,\; r \geq 1. 
      \end{align}
    \item relative isoperimetric inequality of order $d$, that is there exists $C_{\mathrm{riso}} \in (0, \infty)$ such that for all $x \in V$ and $r \geq 1$,
      \begin{align}\label{eq:ass:riso}
        \frac{|\partial_{B(x,r)} A|}{|A|}
        \;\geq\;
        \frac{C_{\mathrm{riso}}}{r},
        \qquad
        \forall\; A \subset B(x,r)\; \text{ s.th. } |A| < \tfrac{1}{2} |B(x,r)|.
      \end{align}
  \end{enumerate}
\end{assumption}
\begin{remark}
  The Euclidean lattice, $(\bbZ^d, E_d)$, satisfies the Assumption \ref{ass:graph}.
\end{remark}
Assume that the graph, $G$, is endowed with positive weights, i.e.\ we consider a family $\om = \{\om(e) \in (0, \infty) : e \in E\}$.  With an abuse of notation we also denote the \emph{conductance matrix} by $\om$, that is for $x, y \in V$ we set $\om(x,y) = \om(y,x) = \om(\{x,y\})$ if $\{x,y\} \in E$ and $\om(x,y) = 0$ otherwise.  We also refer to $\om(x,y)$ as the \emph{conductance} of the corresponding edge $\{x,y\}$ and we call $(V,E,\om)$ a \emph{weighted graph}.  Let us further define measures $\mu^{\om}$ and $\nu^{\om}$ on $V$ by
\begin{align*}
  \mu^{\om}(x) \;\ldef\; \sum_{y \sim x}\, \om(x,y)
  \qquad \text{and} \qquad
  \nu^{\om}(x) \;\ldef\; \sum_{y \sim x}\, \frac{1}{\om(x,y)}.
\end{align*}
For any fixed $\om$  we consider a \emph{reversible} continuous time Markov chain, $Y = (Y_t\!: t \geq 0 )$, on $V$ with generator $\cL^{\om}_Y$ acting on bounded functions $f\!:  V \to \bbR$ as
\begin{align} \label{def:LY}
  \big(\cL_Y^{\om}\, f\big)(x)
  \;=\;
  \mu^\om(x)^{-1} \sum_{y \sim x} \om(x,y) \, \big(f(y) - f(x)\big).
\end{align}
We denote by $\Prob_x^\om$ the law of the process starting at the vertex $x \in V$.  The corresponding expectation will be denoted by $\Mean_x^\om$.  Setting  $p^{\om}(x,y) \ldef \om(x,y) / \mu^{\om}(x)$, this random walk waits at $x$ an exponential time with mean $1$ and chooses its next position $y$ with probability $p^{\om}(x,y)$.  Since the law of the waiting times does not depend on the location, $Y$ is also called the \emph{constant speed random walk} (CSRW). Furthermore, $Y$ is a reversible Markov chain with symmetrising measure given by $\mu$. In \cite{ADS13} we mainly consider the  \emph{variable speed random walk} (VSRW) $X = (X_t\!: t \geq 0 )$, which waits at $x$ an exponential time with mean $1/\mu(x)$, with generator
\begin{align*}
  \big(\cL^{\om} f\big)(x)
  \;=\;
  \sum_{y \in \bbZ^d} \om(x,y)\, \big(f(y) - f(x)\big)
  \;=\;
  \mu^\om(x)\, \big(\cL^{\om}_Y f\big)(x).
\end{align*}

\subsection{Harnack inequalities on graphs}
In the case of uniformly bounded conductances, i.e.\ there exist $0 < \la_1 \leq \la_2 < \infty$ such that
\begin{align*}
  \la_1 \;\leq\; \om(x,y) \;\leq\; \lambda_2,
\end{align*}
both the EHI and the PHI have been derived in the setting of an infinite, connected, locally finite graph by Delmotte  in \cite{De97a, De99} with constant $C_{\mathrm{EH}} = C_{\mathrm{EH}}(\la_1, \la_2)$ and $C_{\mathrm{PH}} = C_{\mathrm{PH}}(\la_1, \la_2)$ depending on $\la_1$ and $\la_2$ only.  It is obvious that for a given ball $B(x_0,n)$ the constants $\la_1$ and $\la_2$ can be replaced by
\begin{align*}
  \max_{x,y \in B(x_0,n)} \frac{1}{\om(x,y)}
  \qquad \text{and} \qquad
  \max_{x,y \in B(x_0, n)} \om(x,y).
\end{align*}
Our main result shows that the uniform upper and lower bound can be replaced by $L^p$ bounds.

\begin{theorem}[Elliptic Harnack inequality]\label{thm:EHI}
  For any $x_0 \in V$ and $n \geq 1$ let $B(n) = B(x_0,n)$.  Suppose that $u > 0$ is harmonic on $B(n)$, i.e.\  $\cL_Y^\om\, u = 0$ on $B(n)$.  Then, for any $p,q \in (1, \infty)$ with
  \begin{align*} 
    \frac{1}{p} + \frac{1}{q} \;<\; \frac{2}{d} 
  \end{align*}
  there exists $C_{\mathrm{EH}} \equiv C_{\mathrm{EH}}\big(\Norm{\mu^\om}{p,B(n)},  \Norm{\nu^\om}{q,B(n)}\big)$ such that
  \begin{align}\label{eq:EHI}
    \max_{x \in B(n/2)} u(x)
    \;\leq\;
    C_{\mathrm{EH}}\, \min_{x \in B(n/2)} u(x).
  \end{align}
  The constant $C_{\mathrm{EH}}$ is more explicitly given by
  \begin{align*}
    C_{\mathrm{EH}}\big(\Norm{\mu^\om}{p,B(n)}, \Norm{\nu^\om}{q,B(n)}\big)
    \;=\;
    c_1\, \exp\!%
    \left(
      c_2\,
      \Big( 1 \vee \Norm{\mu^\om}{p,B(n)}\, \Norm{\nu^\om}{q,B(n)}\Big)^\ka  
    \right)
  \end{align*}
  for some positive $c_i = c_i(d,p,q)$ and $\ka = \ka(d,p,q)$.
\end{theorem}
For abbreviation we introduce  $M_{x_0,n}^\om \ldef 1 \vee \Big(\Norm{\mu^\om}{1,B(n)}/\Norm{\mu^\om}{1,B(n/2)}\Big)$.
\begin{theorem}[Parabolic Harnack inequality]\label{thm:PHI}
  For any $x_0\in V$, $t_0\geq 0$ and $n\geq 1$ let $Q(n) = [t_0,t_0+n^2]\times B(x_0,n)$.  Suppose that $u > 0$ is caloric on $Q(n)$, i.e.\  $\partial_t u - \cL_Y^{\om} u = 0$ on $Q(n)$.  Then, for any $p,q \in (1, \infty)$ with
  \begin{align*} 
    \frac{1}{p} + \frac{1}{q} \;<\; \frac{2}{d} 
  \end{align*}
  there exists $C_{\mathrm{PH}} = C_{\mathrm{PH}}\big(M_{x_0,n}^\om, \|\mu^{\om}\|_{p,B(x_0,n)}, \|\nu^{\om}\|_{q,B(x_0,n)}\big)$ such that
  \begin{align}\label{eq:PHI}
    \max_{(t,x) \in Q_{-}} u(t,x)
    \;\leq\;
    C_{\mathrm{PH}}\, \min_{(t,x) \in Q_+} u(t,x),
  \end{align}
  where $Q_- = \big[t_0 + \frac{1}{4} n^2, t_0 + \frac{1}{2}n^2\big] \times B\big(x_0, \frac{1}{2} n\big)$ and $Q_+ = \big[t_0 + \frac{3}{4} n^2, t_0 + n^2\big] \times B\big(x_0, \frac{1}{2} n\big)$. The constant $C_{\mathrm{PH}}$ is more explicitly given by
  \begin{align*}
   &C_{\mathrm{PH}}%
   \big(
     M^{\om}_{x_0,n}, \Norm{\mu^\om}{p,B(x_0,n)}, \Norm{\nu^\om}{q,B(x_0,n)}
   \big)
   \\[.5ex]
   &\mspace{36mu}=\;
    c_1\, \exp\!%
    \left(
      c_2\, \Big(M^\om_{x_0,n} \,  \big( 1 \vee \Norm{\mu^\om}{p,B(x_0,n)} \big)\,
      \big( 1 \vee \Norm{\nu^\om}{q,B(x_0,n)}\big) \Big)^{\!\ka}  
    \right)
  \end{align*}
  for some positive $c_i = c_i(d,p,q)$ and $\ka = \ka(d,p,q)$.
\end{theorem}
In the context of elliptic PDEs on $\bbR^d$ discussed above, the classical results on EHI and regularity of harmonic function have been extended to the situation, where the uniform ellipticity assumption on the coefficient (see \eqref{eq:ell_cont} ) is weakened  to a certain integrability condition very similar to ours, see \cite{EP72,Kr63,MS68,Tr71,Tr73}. In particular, the same $L^p$-condition $1/p+1/q<2/d$ is obtained. This condition also appears in  \cite{Zh13}, where an upper bound on the solution of a degenerate parabolic equation is derived by using Moser iteration and Davies' method (cf.\ \cite{CKS87}).
\begin{remark} \label{rem:PHI_pi}
  Given a speed measure $\pi^\om\!: \bbZ^d \to (0,\infty)$, one can also consider the process, $Z = (Z_t\!: t \geq 0 )$ on $V$ that is defined by a time change of the VSRW $X$, i.e. $Z_t \ldef X_{a_t}$ for $t \geq 0$, where $a_t \ldef \inf\{s \geq 0 : A_s > t\}$ denotes the right continuous inverse of the functional
  \begin{align*}
    A_t \;=\; \int_0^t \pi^\om(X_s) \, \md s, \qquad t \geq 0.
  \end{align*}
  Its generator is given by
  \begin{align} \label{def:Lpi}
    \big( \cL_{\pi}^{\om} f \big) (x)
    \;=\;
   \pi^{\om}(x)^{-1} \sum_{y \in V} \om(x,y) \big(f(y)-f(x)\big).
  \end{align}
  The maybe most natural choice for the speed measure is $\pi^{\om} = \mu^{\om}$, for which we get again the CSRW $Y$. Since we have the same harmonic functions w.r.t.\ $ \cL_{\pi}^{\om}$  for any choice of $\pi$, it is obvious that the EHI in Theorem~\ref{thm:EHI} also holds for $\cL_{\pi}^{\om}$. On the other hand, for the PHI one can show the following statement along the lines of the proof of Theorem~\ref{thm:PHI} (see \cite[Appendix~A]{CD14} and \cite[Section~3]{ADS14} for a few more details). Let $u > 0$ be caloric w.r.t.\ $ \cL_{\pi}^{\om}$  on $Q(n)$,  then for any $p,q,r \in (1, \infty)$ with
  \begin{align} \label{eq:cond_pqr} 
    \frac{1}{r} \,+\, \frac{1}{p} \frac{r-1}{r} \,+\, \frac{1}{q}
    \;<\;
    \frac{2}{d} 
  \end{align}
  there exists $C_{\mathrm{PH}} = C_{\mathrm{PH}}\big(M_{x_0,n}^\om, \Norm{\mu^{\om} / (\pi^{\om})^{1 - \frac{1}{p}}}{p,B(n)},  \Norm{\nu^\om}{q,B(n)}, \Norm{\pi^\om}{r,B(n)}\big)$ with $M_{x_0,n}^\om\ldef 1 \vee \Big(\Norm{\pi^{\om}}{1,B(n)} / \Norm{\pi^{\om}}{1,B(n/2)}\Big)$ such that
  \begin{align*}
   \max_{(t,x) \in Q_{-}} u(t,x)
    \;\leq\;
    C_{\mathrm{PH}}\, \min_{(t,x) \in Q_+} u(t,x).
  \end{align*}
  In particular, for the VSRW $X$ we have $\pi^{\om} \equiv 1$ and choosing $r=\infty$ the condition on $p$ and $q$ reads again $1/p+1/q <2/d$ in this case.
\end{remark}
Notice that the Harnack constants $C_{\mathrm{EH}}$ and $C_{\mathrm{PH}}$ in Theorem~\ref{thm:EHI} and Theorem~\ref{thm:PHI} do depend on the ball under consideration, namely on its center point $x_0$ as well as on its radius $n$. As a consequence one cannot deduce directly H\"older continuity estimates from these results. This requires ergodicity w.r.t.\ translations as an additional assumption (see Proposition~\ref{prop:cont_harmonic} and Proposition~\ref{prop:cont_caloric} below). 
\begin{assumption} \label{ass:lim_const}
  Assume that
  \begin{align*}
    \bar{\mu}
    \;=\;
    \sup_{x_0 \in V} \limsup_{n\to \infty} \Norm{\mu^{\om}}{p,B(x_0,n)}
    \;<\;
    \infty,
    \qquad
    \bar{\nu}
    \;=\;
    \sup_{x_0 \in V} \limsup_{n\to \infty} \Norm{\nu^{\om}}{q,B(x_0,n)}
    \;<\;
    \infty,
  \end{align*}
  \begin{align*}
    \overline{M}
    \;=\;
    \sup_{x_0 \in V} \limsup_{n\to \infty} M^{\om}_{x_0,n}
    \;<\;
    \infty.
  \end{align*}
\end{assumption}
In particular, Assumption \ref{ass:lim_const} implies for the Harnack constants appearing in Theorems~\ref{thm:EHI} and \ref{thm:PHI} that $C^*_{\mathrm{EH}} \ldef 2 C_{\mathrm{EH}}(\bar{\mu}, \bar{\nu}) < \infty$ and $C^*_{\mathrm{PH}} \ldef 2 C_{\mathrm{PH}}(\overline{M}, \bar{\mu}, \bar{\nu}) < \infty$ do not depend on $x_0$ and $n$. Furthermore, there exists $s^\om(x_0) \geq 1$ such that
  \begin{align*}
    \sup_{n \geq s^\om(x_0)}
    & C_{\mathrm{EH}}(\Norm{\mu^{\om}}{p,B(x_0,n)}, \Norm{\nu^{\om}}{q,B(x_0,n)} )
    \;\leq\;
    C^*_{\mathrm{EH}}
    \;<\;
    \infty,
    \\[.5ex]
    \sup_{n \geq s^\om(x_0)}
    &C_{\mathrm{PH}}(M^\om_{x_0,n},\Norm{\mu^{\om}}{p,B(x_0,n)}, \Norm{\nu^{\om}}{q,B(x_0,n)} )
    \;\leq\;
    C^*_{\mathrm{PH}}
    \;<\;
    \infty.
  \end{align*}
In other words, under Assumption \ref{ass:lim_const} the results in Theorems~\ref{thm:EHI} and \ref{thm:PHI} are Harnack inequalities becoming effective for balls with radius $n$ large enough.
\begin{corro} \label{cor:harnack}
  Suppose that Assumption \ref{ass:lim_const} holds. Then, for any $x_0 \in V$ and $n \geq s^\om(x_0)$ the elliptic and parabolic Harnack inequalities in Theorems~\ref{thm:EHI} and \ref{thm:PHI} hold with the Harnack constants in \eqref{eq:EHI} and \eqref{eq:PHI} replaced by $C^*_{\mathrm{EH}}$ and  $C^*_{\mathrm{PH}}$, respectively.
\end{corro}
In the applications to the random conductance model discussed below the conductances will be stationary ergodic random variables, so under some moment conditions Assumption~\ref{ass:lim_const} will be satisfied by the ergodic theorem.

It is well known that a PHI is equivalent to Gaussian lower and upper bounds on the heat kernel in many situations, for instance in the case of uniformly bounded conductances on a locally finite graph, see \cite{De99}.  Indeed, in a forthcoming paper \cite{ADS14} we will prove upper Gaussian estimates on the heat kernel under Assumption~\ref{ass:lim_const} following the approach in \cite{Zh13}, i.e.\ we combine Moser's iteration technique with Davies' method.  We remark that in our non-elliptic setting we cannot deduce off-diagonal Gaussian bounds directly from the PHI due to the dependence of the Harnack constant on the underlying ball. More precisely, in order to get effective Gaussian off-diagonal bounds, one needs to apply the PHI on a number of balls with radius $n$ having a distance of order $n^2$. In general, the ergodic theorem does not give the required uniform control on the convergence of space-averages of stationary random variables over such balls (see \cite{AJ75}).  However, we still obtain some near-diagonal estimates (see Proposition~\ref{prop:hke} below).

\subsection{The Method}
Since the pioneering works \cite{Mo61,Mo64} Moser's iteration technique is by far the best-established tool in order to prove Harnack inequalities.  Moser's iteration is based on two main ideas: the Sobolev-type inequality which allows us to control the $\ell^r$ norm with $r = r(d) = d/(d-2) > 1$ in terms of the Dirichlet form, and a control on the Dirichlet form of any harmonic (or caloric) function $u$.  In the uniformly elliptic case this is rather standard.  In our case where the conductances are  unbounded from above and not bounded away from zero we need to work with a dimension dependent weighted Sobolev inequality, established in \cite{ADS13} by using H\"older's inequality.  That is, the coefficient $r(d)$ is replaced by
\begin{align*}
  r(d,p,q) \;=\; \frac{d - d/p}{(d-2) + d/q}.
\end{align*}
For the Moser iteration we need $r(d,p,q) > 1$, of course, which is equivalent to $1/p + 1/q < 2/d$ appearing in statements of Theorems~\ref{thm:EHI} and \ref{thm:PHI}.  As a result we obtain a maximum inequality for $u$, that is an estimate for the $\ell^{\infty}$-norm of $u$ on a ball in terms of the  $\ell^{\al}$-norm of $u$ on a slightly bigger ball for some $\al > 0$.  Since the same holds for $u^{-1}$, to conclude the Harnack inequality we are left to link the $\ell^{\al}$-norm and the $\ell^{-\al}$-norm of $u$.  In the setting of uniformly elliptic conductances (see \cite{De97a}) this can be done by using the John-Nirenberg inequality, that is the exponential integrability of BMO functions.  In this paper we establish this link by an abstract lemma of Bombieri and Giusti \cite{BG72}  (see Lemma~\ref{lem:bombieri_giusti} below).  In order to apply this lemma, beside the maximum inequalities mentioned above, we need to establish a weighted Poincar\'e inequality, which is classically obtained by the Whitney covering technique (see e.g.\ \cite{De99, Ba04, SC02}). However, in this paper we can avoid the Whitney covering by using results from a recent work by Dyda and Kassmann \cite{DK13}.

\subsection{Local Limit Theorem for the Random Conductance Model} \label{sec:rcm_intro}
Our main motivation for deriving the above Harnack inequalities is the quenched local CLT for the random conductance model.  Consider the $d$-dimensional Euclidean lattice, $(V_d, E_d)$, for $d \geq 2$.  The vertex set, $V_d$, of this graph equals $\bbZ^d$ and the edge set, $E_d$, is given by the set of all non oriented nearest neighbour bonds, i.e. $E_d \ldef \{ \{x,y\}: x,y \in \bbZ^d, |x-y|=1\}$. 

Let $(\Om, \cF) = \big(\bbR_+^{E_d}, \cB(\bbR_+)^{\otimes\, E_d}\big)$ be a measurable space.  Let the graph $(V_d, E_d)$ be endowed with a configuration $\om = \{\om(e) : e \in E_d\} \in \Om$ of conductances.  We will henceforth denote by $\prob$ a probability measure on $(\Om, \cF)$, and we write $\mean$ to denote the expectation with respect to $\prob$.  A \emph{space shift} by $z \in \bbZ^d$ is a map $\tau_z\!: \Om \to \Om$,
\begin{align}\label{eq:def:space_shift}
  (\tau_z \om)(x, y) \;\ldef\; \om(x+z, y+z),
  \qquad \forall\, \{x,y\} \in E_d.
\end{align}
The set $\big\{\tau_x : x \in \bbZ^d\big\}$ together with the operation $\tau_{x} \circ \tau_{y} \ldef \tau_{x+y}$ defines the \emph{group of space shifts}.  We will study the nearest-neighbor \emph{random conductance model}.  For any fixed realization $\om$ it is a \emph{reversible} continuous time Markov chain, $Y = \{Y_t\!: t \geq 0 \}$, on $\bbZ^d$ with generator $\cL^{\om}_Y$ defined as in \eqref{def:LY}.  We denote by $q^{\om}(t,x,y)$ for $x, y \in \bbZ^d$ and $t \geq 0$ the transition density (or heat kernel associated with $\cL_Y^{\om}$) of the CSRW with respect to the reversible measure $\mu$, i.e.\
\begin{align*}
  q^{\om}(t,x,y) \;\ldef\; \frac{\Prob_x^\om\big[Y_t = y\big]}{\mu^{\om}(y)}.
\end{align*}
As a consequence of \eqref{eq:def:space_shift} we have $q^{\tau_z \om}(t, x, y) \;=\; q^{\om}(t, x+z, y+z)$.
\begin{assumption} \label{ass:environment}
  Assume that $\prob$ satisfies the following conditions:
  \begin{enumerate}[ (i)]
    \item $\prob\!\big[ 0 < \om(e) < \infty\big] = 1$ for all $e \in E_d$ and $\mean[\mu^{\om}(0)] < \infty$.
    \item $\prob$ is ergodic with respect to translations of $\bbZ^d$, i.e. $\prob \circ\, \tau_x^{-1} \!= \prob\,$ for all $x \in \bbZ^d$ and $\prob[A] \in \{0,1\}\,$ for any $A \in \cF$ such that $\tau_x(A) = A\,$ for all $x \in \bbZ^d$.
  \end{enumerate}
\end{assumption}
\begin{assumption}  \label{ass:moment} 
There exist $p,q \in (1, \infty]$ satisfying $1/p + 1/q < 2/d$ such that
  \begin{align}\label{eq:moment_condition}
    \mean\!\big[\om(e)^p\big] \;<\; \infty
    \quad \text{and} \quad
    \mean\!\big[\om(e)^{-q}\big] \;<\; \infty
  \end{align}
  for any $e \in E_d$.
\end{assumption}
Notice that under Assumptions \ref{ass:environment} and \ref{ass:moment} the spatial ergodic theorem gives that for $\prob$-a.e.\ $\om$,
\begin{align*}
  \lim_{n\to \infty} \Norm{\mu^{\om}}{p, B(x_0, n)}^p
  \;=\;
  \mean\!\big[ \mu^{\om}(0)^p\big]
  \;<\;
  \infty,
  \quad
  \lim_{n\to \infty} \Norm{\nu^{\om}}{q, B(x_0, n)}^q
  \;=\;
  \mean\big[\nu^{\om}(0)^q\big]
  \;<\;
  \infty,
\end{align*}
\begin{align*}
  \lim _{n \to \infty} M_{x_0,n}^{\om} \;=\; 1.
\end{align*}
In particular, Assumption \ref{ass:lim_const} is fulfilled in this case and thus for $\prob$-a.e.\ $\om$ the EHI and the PHI do hold for $\cL_Y^{\om}$  in the sense of Corollary~\ref{cor:harnack}, i.e.\ provided $n \geq s^{\om}(x_0)$.

We are interested in the $\prob$-a.s.\ or quenched long range behaviour, in particular in obtaining a quenched functional limit theorem (QFCLT) or invariance principle for the process $Y$ starting in $0$ and a quenched local limit theorem  for $Y$.  We first recall that the following QFCLT has been recently obtained as the main result in \cite{ADS13}.
\begin{theorem}[QFCLT]\label{thm:ip}
  Suppose that $d \geq 2$ and Assumptions \ref{ass:environment} and \ref{ass:moment} hold.  Then, $\prob$-a.s.~$Y_t^{(n)} \ldef \frac{1}{n} Y_{n^2 t}$ converges (under $\Prob_{\!0}^\om$) in law to a Brownian motion on $\bbR^d$  with a deterministic non-degenerate covariance matrix $\Si_Y^2$.
\end{theorem}
\begin{proof}
  See Theorem 1.3 and Remark 1.5 in \cite{ADS13}.
\end{proof}
We refer to \cite{CD14} for a quenched invariance principle for diffusions under an analogue moment condition.  In this paper our main concern is to establish a local limit theorem for $Y$. It roughly describes how the transition probabilities of the random walk $Y$ can be rescaled in order to get the Gaussian transition density of the Brownian motion, which appears as the limit process in Theorem~\ref{thm:ip}.
Write
\begin{align}\label{eq:def:GHK}
  k_t(x)
  \;\equiv\;
  k_t^{\Si_Y}(x)
  \;\ldef\;
  \frac{1}{\sqrt{(2 \pi t)^d \det \Si_Y^2} }\,
  \exp\Big(\!- x \cdot \big( \Si^2_Y \big)^{-1} x / 2t\Big)
\end{align}
for the Gaussian heat kernel with covariance matrix $\Si^2_Y$.  For $x \in \bbR^d$ write $\lfloor x \rfloor = \big( \lfloor x_1 \rfloor, \dots \lfloor x_d \rfloor\big)$.
\begin{theorem}[Quenched Local Limit Theorem]\label{thm:llt_csrw}
  Let $T_2 > T_1>0$ and $K > 0$ and suppose that $d \geq 2$ and Assumptions \ref{ass:environment} and \ref{ass:moment} hold.  Then,
  \begin{align*} 
    \lim_{n \,\to\, \infty}  
    \sup_{|x| \leq K}\, \sup_{t \,\in [T_1,T_2]}\, 
    \Big|
      n^{d}\, q^\om\big(n^2 t, 0, \lfloor nx \rfloor\big) - a\, k_t(x)
    \Big|
    \;=\;
    0,
    \qquad  \prob \text{-a.s.} 
  \end{align*}
  with $a \ldef 1 / \mean[\mu^{\om}(0)]$.
\end{theorem}
\begin{remark}
  Similarly to Remark~\ref{rem:PHI_pi} it is possible to state a local limit theorem for the process $Z$, obtained as a time-change of the VSRW $X$ in terms of a speed measure $\pi^{\om}\!: \bbZ^d \to (0,\infty)$ with generator  $\cL_{\pi}^{\om}$ defined in \eqref{def:Lpi}. Indeed, if $\pi^\om(x) = \pi^{\tau_x \om}(0)$ for every $x \in \bbZ^d$, if $0 < \mean[\pi^{\om}(0)] < \infty$ and if Assumptions \ref{ass:environment} and \ref{ass:moment} hold, we have a QFCLT also for the process $Z$ with some non-degenerate covariance matrix $\Si_Z^2 = \mean[\pi^\om(0)]^{-1} \Si^2_X$, where $\Si^2_X$ denotes the covariance matrix in the QFCLT for the VSRW $X$ (see Remark~1.5 in \cite{ADS13}). If in addition there exists $r > 1$ such that $p$, $q$ and $r$ satisfy \eqref{eq:cond_pqr} and $ \mean[\pi^{\om}(0)^r]< \infty$,  writing $q^{\om}_Z(t,x,y)$ for the heat kernel associated with $\cL_{\pi}^{\om}$ and $Z$, we have
  \begin{align*} 
    \lim_{n \,\to\, \infty}  
    \sup_{|x| \leq K}\, \sup_{t \,\in [T_1,T_2]}\, 
    \Big|
      n^{d}\, q_Z^\om\big(n^2 t, 0, \lfloor nx \rfloor\big)
      \,-\, a\, k_t^{\Si_Z}(x)
    \Big|
    \;=\;
    0,
    \qquad  \prob \text{-a.s.} 
  \end{align*}
  with $a \ldef 1/ \mean[\pi^{\om}(0)]$.
\end{remark}
Analogous results have been obtained in other settings, for instance for random walks on supercritical percolation clusters, see \cite{BH09}. In \cite{CH08}, the arguments of \cite{BH09} have been used in order to establish a general criterion for a local limit theorem to hold, which is applicable in a number of different situations like graph trees converging to a continuum random tree or local homogenisation for nested fractals.

For the random conductance model a local limit theorem has been proven in the case, where the conductances are i.i.d.\ random variables and uniformly elliptic (see Theorem~5.7 in \cite{BH09}).  Later this has been improved for the VSRW under i.i.d.\ conductances, which are only uniformly bounded away from zero (Theorem~5.14 in \cite{BD10}).  However, if the conductances are i.i.d.\ but have fat tails at zero, the QFCLT still holds (see \cite{ABDH13}), but due to a trapping phenomenon the heat kernel decay is sub-diffusive, so the transition density does not have enough regularity for a local limit theorem -- see \cite{BBHK08, BB12}.

Hence, it is clear that some moment conditions are needed. It turns out that the moment conditions in Assumption~\ref{ass:moment} are optimal in the sense that for any $p$ and $q$ satisfying $1/p + 1/q > 2/d$ there exists an environment of ergodic random conductances with $\mean[\om(e)^p] < \infty$ and $\mean[\om(e)^{-q}]<\infty$, under which the local limit theorem fails (see Theorem~\ref{thm:ex_lclt} below). In the case when the conductances are bounded from above,  i.e.\ $p = \infty$, the optimality of the condition $q > d/2$ is already indicated by the examples in \cite{Bou10b, MF06}.

However, in Section~\ref{sec:imprLB} below we will show that for some examples of conductances our moment condition can be improved if it is replaced by a moment condition on the heat kernel.  While for the example in \cite{MF06}, where the conductances are given by $\om(x,y) = \th^{\om}(x) \wedge \th^{\om}(y)$ for a family $\{\th^\om(x): x \in \bbZ^d\}$ of random variables bounded from above,  no improvement can be achieved by this method (see Proposition~\ref{prop:no_impr}), we obtain that $q > \frac{3}{4} \frac{d}{d+1}$ is sufficient if for instance $\om(x,y) = \th^{\om}(x) \vee \th^{\om}(y)$  (see Proposition~\ref{prop:impr2}).  Moreover, in the case of i.i.d.\ conductances uniformly bounded from above this procedure gives the following result (see Proposition~\ref{prop:impr1}).
\begin{theorem} \label{thm:iid_bdd}
  Assume that the conductances $\{\om(e) : e \in E_d\}$ are independent and uniformly bounded from above, i.e.\ $p = \infty$, and $\mean\big[\om(e)^{-q} \big] < \infty$ for some $q > 1/4$. Then, the QFCLT in Theorem~\ref{thm:ip} and the quenched local limit theorem in Theorem~\ref{thm:llt_csrw} hold for both the CSRW and the VSRW.
\end{theorem}
Under the assumptions of Theorem~\ref{thm:iid_bdd} a PHI and a local limit theorem have recently been proven in \cite[Theorem 1.9]{BKM15}  on one hand for the VSRW if again $q>1/4$  and on the other hand for the CSRW if $q>d/(8d-4)$, which improves the result in Theorem~\ref{thm:iid_bdd}. We refer to Remark~1.10 in \cite{BKM15} for a discussion on the optimality of these conditions for both CSRW and VSRW.

We shall prove Theorem \ref{thm:llt_csrw} by using the approach in \cite{BH09} and \cite{CH08}. The two main ingredients of the proof are the QFCLT in Theorem~\ref{thm:ip} and a H\"older-continuity estimate on the heat kernel, which at the end enables us to replace the weak convergence given by the QFCLT by the pointwise convergence in Theorem~\ref{thm:llt_csrw}. 

Finally, notice that our result applies to a random conductance model given by
\begin{align*}
  \om(x,y) \;=\; \exp\!\big(\phi(x) + \phi(y)\big),
  \qquad \{x,y\} \in E_d,
\end{align*}
where for $d \geq 3$, $\{\phi(x) : x \in \bbZ^d\}$ is the discrete massless Gaussian free field, cf.\ \cite{BS12}.  In this case the moment condition \eqref{eq:moment_condition} holds for any $p, q \in (0, \infty)$, of course.

When $d\geq 3$ we can also establish a local limit theorem for the Green kernel. Indeed, as a direct consequence from the PHI we get some on-diagonal bounds on the heat-kernel (see Proposition~\ref{prop:hke} below) from which we obtain the existence of the Green kernel $g^\om(x,y)$ defined by
\begin{align*}
  g^\om (x,y)
  \;=\;
  \int_0^\infty\mspace{-6mu} q^\om(t,x,y)\, \md t.
\end{align*}
\begin{theorem}[Local limit theorem for the Green kernel]\label{thm:llt_green}
  Suppose that $d \geq 3$ and Assumptions \ref{ass:environment} and \ref{ass:moment} hold.  Then, for any $x \not= 0$ we have  
  \begin{align*} 
    \lim_{n \,\to\, \infty}  
    \Big|
      n^{d-2}\, g^{\om}\big(0, \lfloor nx \rfloor\big) - a\, g_{\mathrm{BM}}(x)
    \Big|
    \;=\;
    0,
    \qquad \prob\text{-a.s.} 
  \end{align*}
  with $a \ldef 1 / \mean[\mu^{\om}(0)]$ and $g_{\mathrm{BM}}(x) = \int_0^\infty k_t(x)\, \md t$ denoting the Green kernel of a Brownian motion with covariance matrix $\Si^2_Y$.
\end{theorem}
On one hand, by integration Theorem~\ref{thm:llt_green} can be deduced from Theorem~\ref{thm:llt_csrw}, the on-diagonal bounds obtained in Proposition~\ref{prop:hke} or in \cite{ADS14} and suitable long-range bounds on the heat-kernel (see \cite{Da93}).
Another possibility to prove Theorem~\ref{thm:llt_green} is to combine the QFCLT in Theorem~\ref{thm:ip} with the H\"older continuity of the Green kernel obtained from the EHI in Theorem~\ref{thm:EHI} in a similar way as the local limit theorem follows from the QFCLT and the PHI. 

The paper is organised as follows:  In Section~\ref{sec:sobolev} we provide a collection of Sobolev-type inequalities needed to setup the Moser iteration and to apply the Bombieri-Giusti criterion. Then, Sections~\ref{sec:EHI} and \ref{sec:PHI} contain the proof of the EHI and PHI in Theorem~\ref{thm:EHI} and \ref{thm:PHI}, respectively. The local limit theorem for the random conductance model is proven in Section~\ref{sec:lclt}. Some examples for conductances allowing an improvement on the moment conditions are discussed in Section~\ref{sec:imprLB}. Finally, the appendix contains a collection of some elementary estimates needed in the proofs.

Throughout the paper we write $c$ to denote a positive constant which may change
on each appearance. Constants denoted $C_i$ will be the same through
each argument.

\section{Sobolev and Poincar\'{e} inequalities}
\label{sec:sobolev}
\subsection{Setup and Preliminaries}
For functions $f\!:A \to \bbR$, where either $A \subseteq V$ or $A \subseteq E$, the $\ell^p$-norm $\norm{f}{p}{A}$ will be taken with respect to the counting measure.  The corresponding scalar products in $\ell^2(V)$ and $\ell^2(E)$ are denoted by $\scpr{\cdot}{\cdot}{V}$ and $\scpr{\cdot}{\cdot}{E}$, respectively. Similarly, the $\ell^p$-norm and scalar product w.r.t.\ to a measure $\varphi: V\rightarrow \bbR$ will be denoted by $\norm{f}{p}{V,\varphi}$ and $\scpr{\cdot}{\cdot}{V,\varphi}$, i.e.\  $\scpr{f}{g}{V,\varphi}=\scpr{f}{g \cdot \varphi}{V}$.  For any non-empty, finite $A \subset V$ and $p \in [1, \infty)$, we introduce space-averaged norms on functions $f\!: A \to \bbR$ by 
\begin{align*}
  \Norm{f}{p,A,\varphi}
  \;\ldef\;
  \Bigg(
    \frac{1}{|A|}\; \sum_{x \in A}\, |f(x)|^p \; \varphi(x)
  \Bigg)^{\!\!1/p}.
\end{align*}
The operators $\nabla$ and $\nabla^*$ are defined by $\nabla f\!: E \to \bbR$ and $\nabla^*F\!: V \to \bbR$
\begin{align*}
  \nabla f(e) \;\ldef\; f(e^+) - f(e^-),
  \qquad \text{and} \qquad
  \nabla^*F (x)
  \;\ldef\;
  \sum_{e: e^+ =\,x}\! F(e) \,-\! \sum_{e:e^-=\, x}\! F(e)
\end{align*}
for $f\!: V \to \bbR$ and $F\!: E \to \bbR$, where for each non-oriented edge $e \in E$ we specify out of its two endpoints one as its initial vertex $e^-$ and the other one as its terminal vertex $e^+$.  Nothing of what will follow depend on the particular choice.  Note that $\nabla^*$ is the adjoint of $\nabla$, i.e. for all $f \in \ell^2(V)$ and $F \in \ell^2(E)$ it holds $\scpr{\nabla f}{F}{E} = \scpr{f}{\nabla^* F}{V}$.  We define the products $f \cdot F$ and $F \cdot f$ between a function, $f$, defined on the vertex set and a function, $F$, defined on the edge set in the following way
\begin{align*}
  \big(f \cdot F\big)(e) \;\ldef\; f(e^-)\, F(e),
  \qquad \text{and} \qquad
  \big(F \cdot f\big)(e) \;\ldef\; f(e^+)\, F(e).
\end{align*}
Then, the discrete analog of the product rule can be written as
\begin{align}\label{eq:rule:prod}
  \nabla(f g)
  \;=\;
  \big(g \cdot \nabla f\big) \,+\, \big(\nabla g \cdot f\big).
\end{align}
In contrast to the continuum setting, a discrete version of the chain rule cannot be established.  However, by means of the estimate \eqref{eq:A1:chain:ub1}, $|\nabla f^{\al}|$ for $f \geq 0$ can be bounded from above by
\begin{align}
  \label{eq:rule:chain1}
  \frac{1}{1 \vee |\al|}\,\big|\nabla f^{\al}\big|
  &\;\leq\;
  \big|f^{\al-1} \cdot \nabla f\big| \,+\, \big|\nabla f \cdot f^{\al-1}\big|,
  \qquad \forall\, \al \in \bbR.
  \intertext{On the other hand, the estimate \eqref{eq:A1:chain:lo} implies the following lower bound}
  \label{eq:rule:chain2}
  2\,\big|\nabla f^{\al}\big|
  &\;\geq\;
  \big|f^{\al-1} \cdot \nabla f\big| \,+\, \big|\nabla f \cdot f^{\al-1}\big|,
  \qquad \forall\, \al \geq 1.
\end{align}
The \emph{Dirichlet form} or \emph{energy} associated to $\cL^{\om}_Y$ is defined by
\begin{align} \label{eq:ibp}
  \cE^{\om}(f,g)
  \;\ldef\;
  \scpr{f}{-\cL^{\om}_Yg}{V, \mu^\om}
  \;=\;
  \scpr{\nabla f}{\om \nabla g}{E},
  \qquad
  \cE^{\om}(f) \;\equiv\; \cE^{\om}(f,f).
\end{align}
For a given function $\eta\!: B \subset V \to \bbR$, we denote by $\cE_{\eta^2}^{\om}(u)$ the Dirichlet form where $\om(e)$ is replaced by $\frac{1}{2}(\eta^2(e^+) + \eta^2(e^-))\, \om(e)$ for $e \in E$.

\subsection{Local Poincar\'{e} and Sobolev inequality}
The main objective in this subsection is to establish a version of a local Poincar\'e inequality.

Suppose that the graph $(V, E)$ satisfies the condition (ii) in Assumption~\ref{ass:graph}.  Then, by means of a discrete version of the co-area formula, the classical local $\ell^1$-Poincar\'e inequality on $V$ can be easily established, see e.g.\ \cite[Lemma~3.3]{SC97}, which also implies an $\ell^{\al}$-Poincar\'e inequality for any $\al \in [1, d)$.  Moreover, the condition (i) in Assumption~\ref{ass:graph} ensures that balls in $V$ have a regular volume growth.  Note that, due to \cite[Th\'eor\`eme~4.1]{Cou96}, the volume regularity of balls $B(x_0, n)$ and the local $\ell^{\al}$-Poincar\'e inequality on $V$ implies that for $d \geq 2$ and any $u\!: V \to \bbR$
\begin{align}\label{eq:Sobolev:Sa}
  \inf_{a \in \bbR}
  \Norm{u - a}{\frac{d \al}{d-\al}, B(x_0, n)}
  \;\leq\;
  C_1\, n\,
  \bigg(
    \frac{1}{|B(x_0, n)|}\,
    \sum_{\substack{x, y \in B(x_0, n) \\ x \sim y}} \mspace{-8mu}
    \big| u(x) - u(y) \big|^\al
  \bigg)^{\!\!1/\al}.
\end{align}
This inequality is the starting point to prove a local $\ell^2$-Poincar\'e inequality on the weighted graph $(V, E, \om)$.
\begin{prop}[Local Poincar\'{e} inequality] \label{prop:loc_poinc}
  Suppose $d \geq 2$ and for any $x_0 \in V$ and $n \geq 1$, let $B(n) \equiv B(x_0,n)$.  Then, there exists $C_{\mathrm{PI}} \equiv C_{\mathrm{PI}}(d) < \infty$ such that for any $u\!: V \to \bbR$
  \begin{align}\label{eq:local:PI}
    \Norm{u - (u)_{B(n)}}{2, B(n)}^2
    \;\leq\;
    C_{\mathrm{PI}}\, \Norm{\nu^{\om}}{\frac{d}{2}, B(n)}\;
    \frac{n^2}{|B(n)|}
    \sum_{\substack{x, y \in B(n) \\ x \sim y}} \mspace{-8mu}
    \om(x, y)\, \big( u(x) - u(y) \big)^2,
  \end{align}
  and for any $p, q \in (1, \infty]$ such that $1/p + 1/q = 2/d$,
  \begin{align}\label{eq:local:PI:weighted}
    &\Norm{u - (u)_{B(n), \mu^{\om}}}{2, B(n), \mu^{\om}}^2
    \nonumber\\[.5ex]
    &\mspace{36mu}\leq\;
    C_{\mathrm{PI}}\, \Norm{\mu^{\om}}{p, B(n)}\, \Norm{\nu^{\om}}{q, B(n)}\;
    \frac{n^2}{|B(n)|}
    \sum_{\substack{x, y \in B(n) \\ x \sim y}} \mspace{-8mu}
    \om(x, y)\, \big( u(x) - u(y) \big)^2,
  \end{align}
  where $(u)_{B(n)} \equiv (u)_{B(n),1}$ and $(u)_{B(n), \pi} \ldef \sum_{x \in B(n)} \pi(x) u(x) / \sum_{x \in B(n)} \pi(x)$ for any non-negative weight $\pi$ on $V$.
\end{prop}
\begin{remark}
We prove here a strong version of the local Poincar\'e inequality in the sense that the balls on both sides of inequality have the same radius, in contrast to weak local Poincar\'e inequalities used for instance in \cite{Ba04}. 
\end{remark}
\begin{proof}
  For any $\al \in [1, 2)$, H\"older's inequality yields
  \begin{align}\label{eq:PI:rhs}
    &\bigg(
      \frac{1}{|B(n)|}
      \sum_{\substack{x, y \in B(n)\\ x \sim y}} \mspace{-8mu} |u(x) - u(y)|^{\al}
    \bigg)^{\!\!1/\al}
    \nonumber\\
    &\mspace{36mu}\leq\;
    \Norm{\nu^{\om}}{\frac{\al}{2 - \al}}^{1 / 2}\;
    \bigg(
      \frac{1}{|B(n)|}
      \sum_{\substack{x, y \in B(n)\\ x \sim y}} \mspace{-8mu}
        \om(x, y)\, \big( u(x) - u(y) \big)^2
    \bigg)^{\!\!1/2}.
  \end{align}
  Hence, since $\inf_{a \in \bbR} \|u - a\|_{2, B(n)} = \|u - (u)_{B(n)}\|_{2, B(n)}$, the assertion \eqref{eq:local:PI} follows immediately from \eqref{eq:Sobolev:Sa} by choosing $\al = 2d / (d+2)$.

  Next, for any $p, q \in (1, \infty]$ with $1/p + 1/q = 2/d$ set $\al = 2 p d / (p d - d + 2p)$.  Then, $\al \in [1, 2)$ and $\al / (2 - \al) = q$.  Moreover, an application of H\"older's inequality yields
  \begin{align*}
    \Norm{u - (u)_{B(n), \mu^{\om}}}{2, B(n), \mu^{\om}}^2
    &\;=\;
    \inf_{a \in \bbR} \Norm{u - a}{2, B(n), \mu^{\om}}^2
    \\
    &\;\leq\;
    \Norm{\mu^{\om}}{p, B(n)}\;
    \inf_{a \in \bbR} \Norm{u - a}{2 p_*, B(n)},
  \end{align*}
  where $p_* = p / (p-1)$.  Thus, since $2 p_* = d \al / (d - \al)$, \eqref{eq:local:PI:weighted} follows by combining the last estimate with \eqref{eq:PI:rhs}.
\end{proof}
Later we will also need a weighted version of the local Poincar\'e inequality.
\begin{prop}[Local Poincar\'e inequality with radial cutoff]\label{prop:weight_poinc}
  Suppose $d \geq 2$ and for $x_0 \in V$ and $n \geq 1$ let $B(n) \equiv B(x_0, n)$.  Further, let $\eta$ be a non-negative, radially decreasing weight with $\supp \eta \subset B(n)$, i.e. $\eta(x) = \Phi(d(x_0, x))$ for some non-increasing and non-negative c\`adl\`ag function $\Phi$.  Then, for any $u\!: V \to \bbR$
  \begin{align*}
    \Norm{u - (u)_{B(n), \eta^2}}{2, B(n), \eta^2}^2
    \;\leq\;
    C_{\mathrm{PI}}\, M_1\, n^2\, \Norm{\nu^{\om}}{\frac{d}{2}, B(n)}\,
    \frac{\cE^{\om,\eta}(u)}{|B(n)|},
  \end{align*}
  where $M_1 = 8^2\, \frac{|B(n)|}{|B(n/2)|}\, \frac{\Phi(0)}{\Phi(1/2)}$, and for any $p, q \in (1, \infty]$ such that $1/p + 1/q = 2/d$,
  \begin{align*}
    \Norm{u - (u)_{B(n), \eta^2 \mu^{\om}}}{2, B(n), \eta^2 \mu^{\om}}^2
    \;\leq\;
    C_{\mathrm{PI}}\, M_2\, n^2\, \Norm{\mu^{\om}}{p, B(n)}\Norm{\nu^{\om}}{q, B(n)}\,
    \frac{\cE^{\om,\eta}(u)}{|B(n)|},
  \end{align*}
  where $M_2 = M_1 \Norm{\mu^{\om}}{1,B(n)} / \Norm{\mu^{\om}}{1, B(n/2)}$ and
  \begin{align*}
    \cE^{\om,\eta}(u)
    \;\ldef\;
    \sum_{e \in E} \min\{\eta^2(e^+), \eta^2(e^-)\}\, \om(e)\, 
    \big(\nabla u(e)\big)^2.
  \end{align*}
\end{prop}
\begin{proof}
  Given the local Poincar\'e inequality in Proposition~\ref{prop:loc_poinc}, this follows from the same arguments as in Theorem~1 and Proposition~4 in \cite[Corollary~5]{DK13}.
\end{proof}
From the relative isoperimetric inequality in Assumption~\ref{ass:graph} one can deduce that the following Sobolev inequality holds for any function $u$ with finite support
\begin{align}\label{eq:sob:S1}
  \norm{u}{\frac{d}{d-1}}{V}
  \;\leq\;
  C_{S_1}\, \norm{\nabla u}{1}{E}
\end{align}
(cf.\ Remark 3.3 in \cite{ADS13}).  Note that \eqref{eq:sob:S1} is a Sobolev inequality on an unweighted graph, while for our purposes we need a version involving the conductances as weights.  Let
\begin{align}\label{eq:def:rho}
  \rho \;=\; \rho(q,d) \;\ldef\; \frac{q d}{q(d-2) + d}.
\end{align}
Notice that $\rho(q,d)$ is monotone increasing in $q$ and $\rho(d/2,d) = 1$.  In \cite{ADS13} we derived the following weighted Sobolev inequality by using \eqref{eq:sob:S1} and the H\"older inequality.
\begin{prop}[Sobolev inequality]\label{prop:sobolev}
  Suppose that the graph $(V, E)$ has the isoperimetric dimension $d \geq 2$ and let $B \subset V$ be finite and connected.  Consider a non-negative function $\eta$ with
  \begin{align*}
    \supp \eta \;\subset\; B, \qquad
    0 \;\leq\; \eta \;\leq\; 1 \qquad \text{and} \qquad
    \eta \equiv 0 \quad \text{on} \quad \partial B.
  \end{align*}
  Then, for any $q \in [1, \infty]$, there exists $C_{\mathrm{\,S}} \equiv C_{\mathrm{\,S}}(d,q) < \infty$ such that for any $u:V \to \bbR$,
      \begin{align}\label{eq:sob:ineq:2}
        \Norm{(\eta\, u)^2}{\rho, B}
        \;\leq\;
        C_{\mathrm{\,S}}\, |B|^{\frac{2}{d}}\, \Norm{\nu^{\om}}{q,B}\;
        \Bigg(
          \frac{\cE_{\eta^2}^{\om}(u)}{|B|} \,+\,
          \norm{\nabla \eta}{\infty\!}{E}^2\,
          \Norm{u^2}{1, B, \mu^{\om}}
        \Bigg),
      \end{align}
  where $\cE_{\eta^2}^{\om}(u) \ldef \frac{1}{2} \sum_{x,y} \eta^2(x)\, \om(x,y) \, (u(y)-u(x))^2$.
\end{prop}
\begin{proof}
  See \cite[Proposition~3.5]{ADS13}.
\end{proof}

\subsection{The lemma of Bombieri and Giusti}
In order to obtain the Harnack inequalities we will basically follow the approach in \cite{De97a}. However, due to the lack of a suitable control on the BMO norm and thus of a John-Nirenberg lemma we will use the following abstract lemma  by Bombieri and Giusti (cf.\ \cite{BG72}):
\begin{lemma} \label{lem:bombieri_giusti}
  Let $\{U_\si : \si \in (0,1] \}$ be a collection of subsets of a fixed measure space endowed with a measure $m$ such that and $U_{\si'} \subset U_{\si}$ if $\si' < \si$. Fix $0 < \de < 1$, $0 < \al^* \leq \infty$ and let $f$ be a positive function on $U \ldef U_1$. Suppose that
  \begin{enumerate}[(i)]
  \item  there exists $\ga > 0$ such that for all $\de \leq \si' < \si \leq 1$ and $0 < \al \leq \min\{1, \al^*/2 \}$,
    \begin{align*}
      \bigg(
        \int_{U_{\si'}}\mspace{-3mu} f^{\al^*} \, \md m
      \bigg)^{\!\!1/\al^*} 
      \;\leq\;
      \Big(
        C_{\mathrm{BG1}}\, (\si - \si')^{-\ga}\, m[U]^{-1}
      \Big)^{\!\frac{1}{\al} - \frac{1}{\al^*}}\;
      \bigg( \int_{U_{\si}} f^{\al} \, \md m \bigg)^{\!\!1/\al},
    \end{align*}
  \item for all $\la > 0$,
    \begin{align*}
      m\big[ \ln f > \la \big]
      \;\leq\;
      C_{\mathrm{BG2}}\, m[U] \, \la^{-1},
    \end{align*}
  \end{enumerate}
  with some positive constants $C_{\mathrm{BG1}}$ and $C_{\mathrm{BG2}}$.  Then, there exists $A = A(c_1, c_2, \al^*, \de, \ga)$ such that
  \begin{align*}
    \int_{U_{\de}} f^{\al^*}\, \md m
    &\;\leq\;
    A\, m[U],
    \quad \text{if } \al^* < \infty,
    \qquad \text{and} \qquad
    \sup_{x \in U_{\de}} f(x)
    &\;\leq\;
    A,
    \quad \text{if } \al^* = \infty.
  \end{align*}
  The constant $A$ is more explicitly given by
  \begin{align*}
    A
    \;=\;
    \exp\!\Big(
      c_1 \Big(c_2 + 2\, \big(C_{\mathrm{BG1}} \vee C_{\mathrm{BG2}}\big)^3 \Big)
    \Big)
  \end{align*}
  for some positive $c_1 \equiv c_1(\de, \ga)$ and $c_2 \equiv c_2(\al^*)$.
\end{lemma}
\begin{proof}
  See Lemma 2.2.6 in \cite{SC02}.
\end{proof}

\section{Elliptic Harnack inequality} \label{sec:EHI}
In this section we prove Theorem~\ref{thm:EHI}. 

\subsection{Mean value inequalities}
\begin{lemma} \label{lem:EHI:DF_ualpha}
 Consider a connected, finite subset $B \subset V$ and a function $\eta$ on $V$ with
  \begin{align*}
    \supp \eta \;\subset\; B, \qquad
    0 \;\leq\; \eta \;\leq\; 1 \qquad \text{and} \qquad
    \eta \equiv 0 \quad \text{on} \quad \partial B.
  \end{align*}
  Further, let $u > 0$ be such that $\cL_Y^\om\, u \geq 0$ on $B$.  Then, there exists $C_2 < \infty$ such that for all $\al \geq 1$
  \begin{align}
    \frac{\cE_{\eta^2}^{\om}(u^{\al})}{|B|}
    \;\leq\;
    C_2\, \tfrac{\al^4}{(2\al-1)^2}\;
    \norm{\nabla \eta}{\infty\!}{E}^2\;
    \Norm{u^{2 \al}}{1,B, \mu^{\om}}.
  \end{align}
\end{lemma}
\begin{proof}
  Since $\cL_Y^\om\, u \geq 0$, a summation by parts yields,
  \begin{align}\label{eq:DF:est0}
    0
    \;\geq\;
    \scpr{\eta^2 u^{2\al-1}}{-\cL_Y^\om\, u}{V,\mu^{\om}}
    \;=\;
    \scpr{\nabla (\eta^2 u^{2\al-1})}{\om\, \nabla u}{E}.
  \end{align}
  As an immediate consequence of the product rule \eqref{eq:rule:prod} and \eqref{eq:A1:pol:ub}, we obtain
  \begin{align}\label{eq:DF:est1}
    &\scpr{\nabla (\eta^2 u^{2\al-1})}{\om\, \nabla u}{E}
    \nonumber\\
    &\mspace{36mu}\geq\;
    \frac{2 \al - 1}{\al^2}\; \cE_{\eta^2}^{\om}(u^\al)
    \,-\,
    \frac{1}{2}\,
    \scpr{|\nabla \eta^2|}
         {\om\, (u^{2 \al -1} \cdot |\nabla u| + |\nabla u| \cdot u^{2\al-1})}{E}.
  \end{align}
  Using \eqref{eq:A1:chain:ub2} we find that
  \begin{align*}
    &\scpr{|\nabla \eta^2|}
         {\om\, (u^{2 \al -1} \cdot |\nabla u| + |\nabla u| \cdot u^{2\al-1})}{E}
    \\[.5ex]
    &\mspace{36mu}\leq\;
    4\,
    \scpr{(\eta \cdot \sqrt{\om} + \sqrt{\om} \cdot \eta)\, |\nabla u^{\al}|}
         {\sqrt{\om}\,(u^{\al} \cdot |\nabla \eta| + |\nabla \eta| \cdot u^{\al})}
         {E}
  \end{align*}
  Thus, by applying the Young inequality $|a b| \leq \frac{1}{2}(\ve a^2 + b^2 / \ve)$ and combining the result with \eqref{eq:DF:est1}, we obtain
  \begin{align}\label{eq:DF:lb}
    &\scpr{\nabla (\eta^2 u^{2\al-1})}{\om\, \nabla u}{E}
    \nonumber\\[.5ex]
    &\mspace{36mu}\geq\;
    \bigg(\frac{2 \al - 1}{\al^2} - 4 \ve\; \bigg)\;
    \cE_{\eta^2}^{\om}(u^{\al}) \,-\, \frac{4}{\ve}\; |B|\;
    \norm{\nabla \eta}{\infty\!}{E}^2\; \Norm{u^{2 \al}}{1,B, \mu^{\om}}.
  \end{align}
  By choosing $\ve = \frac{2 \al - 1}{8 \al^2}$ and combining the result with \eqref{eq:DF:est0}, the claim follows.
\end{proof}
\begin{prop} \label{prop:EHI:mos_it}
  For any $x_0 \in V$ and $n \geq 1$, let $B(n) \equiv B(x_0,n)$.  Let $u > 0$ be such that $\cL_Y^\om\, u \geq 0$ on $B(n)$.  Then, for any $p,q \in (1, \infty]$ with
  \begin{align*} 
    \frac{1}{p} + \frac{1}{q} \;<\; \frac{2}{d} 
  \end{align*}
  there exists $\ka =\ka(d,p,q) \in (1, \infty)$ and $C_3 = C_3(d,q)$ such that for all $\be \in [2 p_*,\infty]$ and for all $1/2 \leq \si' < \si \leq 1$ we have
  \begin{align}\label{eq:EHI:mos_it}
    \Norm{u}{\be,B(\si' n)}
    \;\leq\;
    C_3\,
    \Bigg( 
      \frac{1 \vee \Norm{\mu^\om}{p,B(n)}\, \Norm{\nu^\om}{q,B(n)}}
           { (\si - \si')^{2}}
    \Bigg)^{\!\!\ka}\; \Norm{u}{2 p_*,B(\sigma n)},
  \end{align}
  where $p_* = p / (p-1)$.
\end{prop}
\begin{proof}
  For fixed $1/2 \leq \si' < \si \leq 1$, let $\{B(\si_k n)\}_k$ be a sequence of balls with radius $\si_k n$ centered at $x_0$, where
  \begin{align*}
    \si_k \;=\; \si' + 2^{-k} (\si - \si')
    \qquad \text{and} \qquad
    \tau_k \;=\; 2^{-k-1} (\si - \si'),
    \quad k = 0, 1, \ldots
  \end{align*}
  Note that $\si_k = \si_{k+1} + \tau_k$ and $\si_0 = \si$.  Set $\al_k = (\rho / p_*)^k$.  Since $1/p + 1/q < 2/d$ we have $\rho > p_*$, and, in particular, $\al_k \geq 1$ for every $k$.  Due to the discreteness of the underlying space $\bbZ^d$, we distinguish two different cases.

  Consider first the case $\tau_k n < 1$, that is $B(\si_{k+1} n ) = B(\si_k n)$.  Note that
  \begin{align*}
    \Norm{u^{2 \al_k}}{\al p_*, B(\si_{k+1} n)}
    &\;\leq\;
    \Norm{u^{2 \al_k}}{p_*, B(\si_k n)}^{1/\al}
    \Big( \max_{x \in B(\si_k n)} u(x)^{2 \al_k} \Big)^{\!\!1/\al_*}
    \\[.5ex]
    &\;\leq\;
    |B(\si_k n)|^{1/(\al_* p_*)} \Norm{u^{2 \al_k}}{p_*, B(\si_k n)}.
  \end{align*}
  Since $d/(2 \al_* p_*) \leq 1$ and $n < 1 / \tau_k$, we find
  \begin{align}\label{eq:EHI:iter:0}
    \Norm{u}{2 \al_{k+1} p_*,B(\si_{k+1} n)}
    \;\leq\;
    \bigg(
      c\, \frac{2^{2k}}{(\si - \si')^{2}}\,
    \bigg)^{\!\!1/(2\al_k)}\,
    \Norm{u}{2 \al_k p_*,B(\si_k n)}.
  \end{align} 

  Consider now the case $\tau_k n \geq 1$.  Let $\eta_k$ be a cut-off function with $\supp \eta_k \subset B(\si_{k} n)$ having the property that $\eta_k \equiv 1$ on $B(\si_{k+1} n)$, $\eta_k \equiv 0$ on $\partial B(\si_{k} n)$ and linear decaying on $B(\si_{k} n) \setminus B(\si_{k+1} n)$.  This choice of $\eta_k$ implies that $|\nabla \eta_k(e)| \leq 1/\tau_{k} n$ for all $e \in E$.  Then, by applying the Sobolev inequality \eqref{eq:sob:ineq:2} to $u^{\alpha_k}$, we get
  \begin{align*}
    &\Norm{(\eta_k\, u^{\al_k})^2}{\rho, B(\si_{k} n)}
    \nonumber\\[.5ex]
    &\mspace{36mu}\leq\;
    C_{\mathrm{S}}\, |B(\si_k n)|^{\frac{2}{d}}\, \Norm{\nu^{\om}}{q,B(\si_k n)}\,
    \Bigg(
      \frac{\cE_{\eta_k^2}^{\om}(u^{\al_k})}{|B(\si_{k} n)|}
      \,+\,
      \frac{\Norm{\mu^{\om}}{p,B(\si_k n)}}{(\tau_k n)^2}\,
      \Norm{u^{2 \al_k}}{p_*,B(\si_{k} n)}
    \Bigg).
  \end{align*}
  On the other hand, by using Lemma~\ref{lem:EHI:DF_ualpha} and H\"{o}lder's inequality with $1/p + 1/p_* = 1$, we find
  \begin{align*}
    \frac{\cE_{\eta_k^2}^{\om}(u^{\al_k})}{|B(\si_k n)|}
    \;\leq\;
    C_2\, \bigg(\frac{\al_k}{\tau_k n}\bigg)^{\!\!2}\;
    \Norm{\mu^{\om}}{p,B(\si_k n)}\;
    \Norm{u^{2 \al_k}}{p_*,B(\si_{k} n)}.
  \end{align*}
  Since $\alpha_{k+1} p_* = \alpha_k \rho$, we obtain by combining these two estimates and using the volume regularity which implies that $|B(n)| / |B(\si'n)| \leq C_{\mathrm{reg}}^2 2^d$ that
  \begin{align}\label{eq:EHI:iter:1}
    \Norm{u}{2 \al_{k+1} p_*, B(\si_{k+1} n)}
    &\;\leq\;
    \bigg(
      c\;
      \frac{2^{2k}\, \al_k^2}{(\si -\si')^{2}}\;
      \Norm{\mu^\om}{p,B(n)}\, \Norm{\nu^\om}{q,B(n)}\,
    \bigg)^{\!\!1/(2\al_k)}\;
    \Norm{u}{2 \al_k p_*,B(\si_{k} n)}
  \end{align}
  for some $c < \infty$.  By iterating the inequalities \eqref{eq:EHI:iter:0} and \eqref{eq:EHI:iter:1}, respectively, and using the fact that $\sum_{k=0}^{\infty}k/\al_k < \infty$ and $\sup_k 2^{(k+1)/2 \alpha_k} < \infty$, there exists $C_3 < \infty$ independent of $n$ such that, for any $K \in \bbN$,
  \begin{align*}
    \Norm{u}{2 \al_K p_*, B(\si' n)}
    &\;\leq\;
    C_3\, \prod_{k=0}^{K-1}
    \Bigg(
      \frac{1 \vee \Norm{\mu^\om}{p,B(n)}\, \Norm{\nu^\om}{q,B(n)}}
      {(\si -\si')^{2}}
    \Bigg)^{\!\!1/(2\al_k)}\;
    \Norm{u}{2 p_*,B(\si n)}.
  \end{align*}
  Setting $\ka \ldef \frac{1}{2} \sum_{k=0}^\infty (1/\al_k) < \infty$ and using that $B(\si_K n) \downarrow B(\si' n)$, we get
  \begin{align*}
    \max_{x \in B(\si' n)} u(x)
    &\;=\;
    \lim_{K \to \infty} \Norm{u}{2 \al_K p_*, B(\si_K n)}
    \\[.5ex]
    &\;\leq\;
    C_3\,
    \Bigg(
      \frac{1 \vee \Norm{\mu^\om}{p,B(n)}\, \Norm{\nu^\om}{q,B(n)}}
      {(\si -\si')^{2}}
    \Bigg)^{\!\!\ka}\;
    \Norm{u}{2 p_*, B(\si n)}.
  \end{align*}
  For any $\be \in [2p_*, \infty)$ the claim is immediate since $\|u\|_{\be, B(\si' n)} \leq \max_{x \in B(\si' n)} u(x)$. 
\end{proof}
\begin{corro} \label{cor:EHI:max_est1}
  Let $u > 0$ be such that $\cL_Y^{\om}\, u = 0$ on $B(n)$.  Then, for all $\al \in (0, \infty)$ and  $1/2 \leq \si' < \si \leq 1$ there exists $C_4 < \infty$ such that
  \begin{align} \label{eq:EHI:max_est1}
    \max_{x \in B(\si' n)} u(x)^{-1}
    \;\leq\;
    \Bigg(
      C_4\,
      \frac{1 \vee \Norm{\mu^{\om}}{p,B(n)}\, \Norm{\nu^{\om}}{q,B(n)}}
           {(\si - \si')^{2}}
    \Bigg)^{\!\!2 p_* \ka/\al}\; 
    \Norm{u^{-1}}{\al, B(\si n)}.
  \end{align}
\end{corro}
\begin{proof}
  Since $u > 0$ is harmonic on $B(n)$, the function $f = u^{-\al / 2p_*}$ is sub-harmonic on $B(n)$, i.e.\ $\cL_Y^{\om}\, f \geq 0$.  Thus, \eqref{eq:EHI:max_est1} follows by applying \eqref{eq:EHI:mos_it} with $\be = \infty$ to $f$.
\end{proof}
\begin{corro} \label{cor:EHI:max_est2}
  Consider a non-negative $u$ such that $\cL_Y^{\om}\, u \geq 0$ on $B(n)$.  Then, for all $\al \in (0, \infty)$ and  $1/2 \leq \si' < \si \leq 1$, there exists $C_5(\al) < \infty$ such that
  \begin{align} \label{eq:EHI:max_est2}
    \max_{x \in B(\si' n)} u(x)
    \;\leq\;   
    \Bigg(
      C_5(\al)\,
      \frac{1 \vee \Norm{\mu^{\om}}{p,B(n)}\, \Norm{\nu^{\om}}{q,B(n)}}
           {(\si - \si')^{2}}
    \Bigg)^{\!\!\ka'}\; 
    \Norm{u}{\al, B(\si n)}
  \end{align}
  with  $\ka'=\big(1\vee \frac{2p_*}{\alpha}\big) \kappa>1$.
\end{corro}
\begin{proof}
  The proof of this Corollary, based on arguments given originally in \cite[Theorem 2.2.3]{SC02}, can be found in \cite[Corollary~3.9]{ADS13}.  Nevertheless, we will repeat it here for the reader's convenience.

  In view of \eqref{eq:EHI:mos_it}, for any $\al \geq 2p_*$ the statement \eqref{eq:EHI:max_est2} is an immediate consequence of Jensen's inequality.  It remains to consider the case $\al \in (0, 2p_*)$.  For all $1/2 \leq \si' < \si \leq 1$ we have from \eqref{eq:EHI:mos_it} with $\be =\infty$,
  \begin{align}\label{eq:EHI:max_est3}
    \max_{x \in B(\si' n)} u(x)
    \;\leq\;
    C_3\,
    \Bigg( 
      \frac{1 \vee \Norm{\mu^{\om}}{p,B(n)}\, \Norm{\nu^{\om}}{q, B(n)}}
           {(\si - \si')^{2}}
    \Bigg)^{\!\!\ka}\;
    \Norm{u}{2p_*,B(\si n)}.
  \end{align}
  In the sequel, let $1/2 \leq \si' < \si \leq 1$ be arbitrary but fixed and set $\si_k = \si - 2^{-k}(\si - \si')$ for any $k \in \bbN_0$.  Now, by H\"older's inequality, we have for any $\al \in (0, 2p_*)$
  \begin{align*}
    \Norm{u}{2p_*, B(\si_k n)}
    \;\leq\;
    \Norm{u}{\al, B(\si_k n)}^{\theta}\; \Norm{u}{\infty, B(\si_k n)}^{1-\theta}
  \end{align*}
  where $\theta = \al / 2 p_*$.  Hence, in view of \eqref{eq:EHI:max_est3} and the volume regularity which implies that $|B(\si n)| / |B(\si' n)| \leq C_{\mathrm{reg}}^2 2^d$, we obtain
  \begin{align}
    \Norm{u}{\infty, B(\si_{k-1} n)}
    \;\leq\;
    2^{2\ka k}\; J\;
    \Norm{u}{\al, B(\si n)}^{\theta}\;
    \Norm{u}{\infty, B(\si_k n)}^{1 - \theta},
  \end{align}
  where we introduced $J = c\, \big(\Norm{\mu^{\om}}{p,B(n)}\, \Norm{\nu^{\om}}{q, B(n)} / (\si - \si')^{2} \big)^{\ka}$ to simplify notation.  Hence, by iteration, we get
  \begin{align}
    \Norm{u}{\infty, B(\si' n)}
    \;\leq\;
    2^{2\ka \sum_{k=0}^{i-1} (k+1)(1 - \theta)^k}\;
    \Big(
      J\; \Norm{u}{\al, B(\si n)}^{\theta}
    \Big)^{\!\sum_{k=0}^{i-1} (1 - \theta)^k}\; 
    \Norm{u}{\infty, B(\si_i n)}^{(1 - \theta)^i}.
  \end{align}
  As $i$ tends to infinity, this yields $\|u\|_{\infty, B(\si' n)} \leq 2^{2\ka / \th^2}\, J^{1 / \th}$ and \eqref{eq:EHI:max_est2} is immediate.
\end{proof}
The proof of the EHI or more precisely the verification of condition (i) in the lemma of Bombieri-Giusti requires a certain bound on the averaged $\ell^\al$-norm of a harmonic function $u$ for an arbitrarily small $\al > 0$.  However, the constant $C_5(\al)$ diverges as $\al$ tends to zero so that Corollary~\ref{cor:EHI:max_est2} cannot directly be used in order to verify condition (i) in Lemma~\ref{lem:bombieri_giusti}.  Therefore, we proceed by establishing a bound for super-harmonic functions, for which another Moser iteration is needed. 
\begin{lemma} \label{lem:EHI:DF_ualpha_super}
  Consider a connected, finite subset $B \subset V$ and a function $\eta$ on $V$ with
  \begin{align*}
    \supp \eta \;\subset\; B, \qquad
    0 \;\leq\; \eta \;\leq\; 1 \qquad \text{and} \qquad
    \eta \equiv 0 \quad \text{on} \quad \partial B.
  \end{align*}
  Further, let $u > 0$ be such that $\cL_Y^\om\, u \leq 0$ on $B$.  Then, for any $\al < \frac 1 2$,
  \begin{align} \label{eq:EHI:DFsuper}
    \frac{\cE_{\eta^2}^{\om}(u^{\al})}{|B|}
    \;\leq\;
    \frac{32 \al^2}{(1-2\al)^2} \, \osr\big( \eta^2 \big) \,
    \norm{\nabla \eta}{\infty\!}{E}^2\;
    \Norm{u^{2 \al}}{1,B, \mu^{\om}},
  \end{align}
  where $\osr(\eta) \ldef \max\big\{\eta(y) / \eta(x) \wedge 1 \mid \{x,y\} \in E,\; \eta(x) \ne 0\big\}$.
\end{lemma}
\begin{proof}
  Since $\cL_Y^\om\, u \leq 0$, a summation by parts yields,
  \begin{align}\label{eq:DFsuper:est_0}
    0
    \;\leq\; 
    \scpr{\eta^2 u^{2\al-1}}{-\cL_Y^\om\, u}{V,\mu^{\om}}
    \;=\;
    \scpr{\nabla (\eta^2 u^{2\al-1})}{\om\, \nabla u}{E}.
  \end{align}
  Further, using \eqref{eq:A2:beta} (with $\be = 1-2\al$) and \eqref{eq:chain_super} we obtain for any $\{x,y\}\in E$ that
  \begin{align*}
    &\nabla\big(\eta^2 u^{2\al-1}\big)(x,y)\, (\nabla u)(x,y)
    \\[.5ex]
    &\mspace{36mu}\leq\;
    \frac{1}{2}\, \min\big\{\eta^2(x), \eta^2(y)\big\}\,
    \big(\nabla u^{2\al-1}\big)(x,y)\, (\nabla u)(x,y)
    \\
    &\mspace{56mu} +\,
    \frac{8}{1-2\al}\, \osr\big(\eta^2\big)\, (\nabla \eta)^2(x,y) \, \big( u^{2\al}(x)+u^{2\al}(y) \big)
    \\[.5ex]
    &\mspace{36mu}\leq\;
    \frac{2\al -1}{2\al^2}\, \big(\nabla u^{\al}\big)(x,y)
    \,+\,
    \frac{8}{1-2\al}\, \osr\big(\eta^2\big)\,
    (\nabla \eta)^2(x,y)\, \big( u^{2\al}(x)+u^{2\al}(y) \big),
  \end{align*}
  and therefore
  \begin{align}\label{eq:DFsuper:est1}
    &\scpr{\nabla (\eta^2 u^{2\al-1})}{\om\, \nabla u}{E}
    \\[.5ex]
    &\mspace{36mu} \leq\;
    \frac{2\al -1}{2\al^2} \; \cE_{\eta^2}^{\om}(u^\al)
    \,+\,  \frac{16 }{1-2\al} \,|B| \,  \osr\big( \eta^2 \big) \,
    \norm{\nabla \eta}{\infty\!}{E}^2\;
    \Norm{u^{2 \al}}{1,B, \mu^{\om}},
    \nonumber
  \end{align}
  which together with \eqref{eq:DFsuper:est_0} gives the claim.
\end{proof} 
\begin{prop} \label{prop:EHI:mos_it_super}
  For any $x_0 \in V$ and $n \geq 1$, let $B(n) \equiv B(x_0,n)$.  Let $u > 0$ be such that $\cL_Y^\om\, u \leq 0$ on $B(n)$.  Then, for any $p, q \in (1, \infty]$ with
  \begin{align*} 
    \frac{1}{p} + \frac{1}{q} \;<\; \frac{2}{d}
  \end{align*}
  and any $\be \in (0, p_*/2)$ there exists $C_6 = C_6(d, p, q)$ and $\ka = \ka(d, p, q)$ such that for all $\al \in (0, \be / 2)$ and for all $1/2 \leq \si' < \si \leq 1$ we have
  \begin{align}\label{eq:EHI:mos_it_super}
    \Norm{u}{\be,B(\si' n)}
    \;\leq\;
    \Bigg( C_6 \,
      \frac{1 \vee \Norm{\mu^\om}{p,B(n)}\, \Norm{\nu^\om}{q, B(n)}}
           { (\si - \si')^{2}}
    \Bigg)^{\!\!\ka\, (1/\al - 1/\be)}\; \Norm{u}{\al, B(\si n)}.
  \end{align}
\end{prop}
\begin{proof}
  We will proceed as in \cite[Theorem~2.2.5]{SC97}.  Since $1/p + 1/q < 2/d$, let us recall that $\rho / p_* > 1$.  W.l.o.g.\ we also assume that $\rho/p_* \leq 2$ (otherwise we replace $\rho$ by $\rho' \ldef 2p_* > 1$ in the following argument).  Then, there exists $j \geq 2$ such that
  \begin{align*}
    \be\, \bigg( \frac{\rho}{p_*} \bigg)^{\!\!-j}
    \;\leq\;
    \al
    \;<\;
    \be\, \bigg( \frac{\rho}{p_*} \bigg)^{\!\!-(j-1)}.
  \end{align*}
  Further, set
  \begin{align*}
    \al_k
    \;\ldef\;
    \frac{\be}{2 p_*} \bigg( \frac \rho {p_*} \bigg)^{\!\!k-j}
    \leq\;
    \frac{\be}{2 p_*} 
    \;<\;
    \frac{1}{4},
    \qquad \forall\,k = 0, \ldots, j.
  \end{align*}
  In the sequel, let $\si_k$ and $\tau_k$ be defined as in the proof of Proposition~\ref{prop:EHI:mos_it} above.  Again we distinguish two cases. First, for $\tau_k n < 1$, we obtain as before
  \begin{align}\label{eq:EHI:iter_super:0}
    \Norm{u}{2 \al_{k+1} p_*, B(\si_{k+1} n)}
    &\;\leq\;
    \bigg(
      c\,
      \frac{2^{2k}}
      {(\si -\si')^{2}}\;
    \bigg)^{\!\!1/(2\al_k)}\,
    \Norm{u}{2 \al_k p_*,B(\si_{k} n)}.
  \end{align}

  Next consider the case $\tau_k n \geq 1$.  Let $\eta_k$ be defined as in the proof of Proposition~\ref{prop:EHI:mos_it}.  Note that $\osr(\eta^2_k) \leq 4$ for all $k$ such that $\tau_k n \geq 1$.  Then, as in the proof of Proposition~\ref{prop:EHI:mos_it}, we apply again the Sobolev inequality \eqref{eq:sob:ineq:2} to $u^{\al_k}$ and afterwards we use Lemma~\ref{lem:EHI:DF_ualpha_super}, the volume regularity and the relation $\al_{k+1} p_* = \al_k \rho$ to obtain
  \begin{align}\label{eq:EHI:iter_super}
    \Norm{u}{2 \al_{k+1} p_*, B(\si_{k+1} n)}
    &\;\leq\;
    \bigg(
      c\, 2^{2k}\,
      \frac{1 \vee \Norm{\mu^\om}{p,B(n)}\, \Norm{\nu^\om}{q,B(n)}}
      {(\si -\si')^{2}}\;
    \bigg)^{\!\!1/(2\al_k)}\,
    \Norm{u}{2 \al_k p_*,B(\si_{k} n)}
  \end{align}
  for some $c < \infty$ (note that for all $\al_k \in (0, 1/4)$ the constant in \eqref{eq:EHI:DFsuper} is uniformly bounded in $k$).   Further, since $2p_* \al_0 \leq \al$ by Jensen's inequality we have that $\|u\|_{2\al_0 p_*, B(\si n)} \leq \|u\|_{\al,B(\si n)}$.  By iterating \eqref{eq:EHI:iter_super:0} and \eqref{eq:EHI:iter_super}, respectively, we get
  \begin{align*}
    \Norm{u}{\be, B(\si' n)}
    &\;\leq\;
     \prod_{k=0}^{j-1}
    \Bigg( c\, 2^{2k} \,
      \frac{1 \vee \Norm{\mu^\om}{p,B(n)}\, \Norm{\nu^\om}{q,B(n)}}
      {(\si -\si')^{2}}
    \Bigg)^{\!\!1/(2\al_k)}\;
    \Norm{u}{\al,B(\sigma n)}.
  \end{align*}
  Since $\rho / p_* \geq 1$ and $j \geq 2$ we have that
  \begin{align*}
    \frac{(\rho/p_*)^j}{\beta} \,-\,
    \frac{1}{\be}
    \;\leq\;
    \Big( 1 + \frac{\rho}{p_*} \Big)\,
    \bigg( \frac{1}{\al} - \frac{1}{\be} \bigg)  
  \end{align*}
  and therefore
  \begin{align*}
    \sum_{k=0}^{j-1} \frac{1}{\al_k}
    \;=\;
    \frac{2\, p_* \rho}{\rho - p_*}\,
    \bigg( \frac{(\rho/p_*)^j}{\be} -\frac{1}{\be}  \bigg)
    \;\leq\;
    \frac{ 2\, \rho(\rho+p_*)}{\rho-p_*}\,
    \bigg( \frac{1}{\al} - \frac{1}{\be} \bigg).
  \end{align*}
  Finally, since $\sum_{k=0}^{j-1} (k/\al_k) \leq c \sum_{k=0}^{j-1}1 /\al_k < \infty$ for some $c=c(d,p,q)$ the claim follows.
\end{proof}

\subsection{From mean value inequalities to Harnack principle}
\begin{lemma}\label{lem:EHI:DF_ln}
  Consider a connected, finite subset $B \subset V$ and a function $\eta$ on $V$ with
  \begin{align*}
    \supp \eta \;\subset\; B, \qquad
    0 \;\leq\; \eta \;\leq\; 1 \qquad \text{and} \qquad
    \eta \equiv 0 \quad \text{on} \quad \partial B.
  \end{align*}
  Further, let $u > 0$ be such that $\cL_Y^\om\, u \leq 0$ on $B$. Then,
  \begin{align}\label{eq:EHI:DF_ln}
    \frac{\cE^{\om, \eta}(\ln u)}{|B|}
    \;\leq\;
    8\, \osr\big(\eta^2\big)\;
    \norm{\nabla \eta}{\infty\!}{E}^2\;
    \Norm{\mu^\om}{1,B},
  \end{align}
  where $\osr$ is defined as in Lemma \ref{lem:EHI:DF_ualpha_super}.
\end{lemma}
\begin{proof}
  Since $u > 0$ and $\cL_Y^{\om}\, u \leq 0$ on $B$, a summation by parts yields
  \begin{align}\label{eq:EHI:log:est}
    0
    \;\leq\;
    \scpr{\eta^2 u^{-1}}{-\cL_Y^{\om}\, u}{V,\mu^{\om}}
    \;=\;
    \scpr{\nabla(\eta^2 u^{-1})}{\om\, \nabla u}{E}.
  \end{align}
  In view of the estimates \eqref{eq:A2:beta} (with the choice $\beta=1$) and \eqref{eq:A1:log:ub}, we obtain for any $\{x, y\} \in E$ that
  \begin{align}\label{eq:EHI:log:est1}
    &\nabla\big(\eta^2 u^{-1}\big)(x,y)\, (\nabla u)(x,y)
    \nonumber\\[1ex]
    &\mspace{32mu}\leq\;
    \frac{1}{2}\, \min\big\{\eta^2(x), \eta^2(y)\big\}\,
    \big(\nabla u^{-1}\big)(x,y)\, (\nabla u)(x,y)
    \,+\,
    8\, \osr\big(\eta^2\big)\, (\nabla \eta)^2(x,y)
    \nonumber\\[1ex]
    &\mspace{32mu}\leq\;
    - \frac{1}{2}\, \min\big\{\eta^2(x), \eta^2(y)\big\}\,
    \big(\nabla \ln u \big)^2(x,y)\,
    \,+\,
    8\, \osr\big(\eta^2\big)\, (\nabla \eta)^2(x,y).
  \end{align}
  Thus, by combining the resulting estimate with \eqref{eq:EHI:log:est} the claim \eqref{eq:EHI:DF_ln} follows. 
\end{proof}
\begin{proof} [Proof of Theorem~\ref{thm:EHI}]
  Consider a harmonic function $u > 0$ on $B(n) \equiv B(x_0,n)$.  Pick $\vr$ so that $3/4 \leq \vr < 1$ and set
  \begin{align*}
    \eta(x)
    \;\ldef\;
    \bigg(
      1 \,-\, \frac{\big[ d(x_0, x) - \vr n\, \big]_+ }{(1-\vr)\, n}
    \bigg) \;\vee\; 0.
  \end{align*}
  Obviously, the cut-off function $\eta$ has the properties that $\supp \eta \subset B(n)$, $\eta \equiv 1$ on $B(\vr n)$, $\eta \equiv 0$ on $\partial B(n)$, $\|\nabla \eta\|_{\ell^\infty(E)} \leq 1 / (1-\vr) n$ and $\osr(\eta^2) \leq 4$.  

  Now let $f_1 \ldef  u^{-1}\, \me^{\xi}$ with $\xi = \xi(u) \ldef (\ln u)_{B(n),\eta^2}$.  Our aim is to apply Lemma~\ref{lem:bombieri_giusti} to the function $f_1$ with $U = B(\vr n)$, $U_\sigma=B(\sigma \vr n)$, $\de = 1 / 2 \vr$, $\al^* = \infty$ and $m$ being the counting measure.  Note that Corollary~\ref{cor:EHI:max_est1} implies that the condition (i) of Lemma \ref{lem:bombieri_giusti} is satisfied with $\ga = 4p_* \ka$, where $\ka$ and $p_*$ are chosen as in Proposition~\ref{prop:EHI:mos_it} and
  \begin{align*}
    C_{\mathrm{BG1}}
    \;=\;
    \Big(
      C_4\,
      \Big(1 \vee \Norm{\mu^\om}{p,B(n)}\, \Norm{\nu^\om}{q,B(n)}\Big)
    \Big)^{\!2p_* \ka}\mspace{-26mu}.
  \end{align*}
  To verify the second condition of Lemma~\ref{lem:bombieri_giusti}, we apply the weighted local $\ell^2$-Poincar\'{e} inequality in Proposition~\ref{prop:weight_poinc} to the function $\ln u$.  This yields
  \begin{align*}
    \Norm{\ln u - (\ln u)_{B(n),\eta^2}}{2, B(n), \eta^2}^2
    &\;\leq\;
    C_{\mathrm{PI}}\, M_1\, n^2\, \Norm{\nu^{\om}}{\frac{d}{2}, B(n)}\;
    \frac{\cE_{\eta^2 \wedge \eta^2}^{\om}(\ln u)}{|B(n)|}
    \\[1ex]
    &\;\leq\;
    c\, \Norm{\nu^{\om}}{\frac{d}{2}, B(n)}\, \Norm{\mu^{\om}}{1, B(n)},
  \end{align*}
  where we used \eqref{eq:EHI:DF_ln} in the last step.  Hence,
  \begin{align*}
    \la\, \Norm{\indicator_{\{ \ln f_1 \,>\, \la \}}}{1,B(\vr n)}
    \;\leq\;
    \frac{|B(n)|}{|B(\vr n)|}\;
    \Norm{\ln u - \xi}{1, B(n), \eta^2}
    \;\leq\;
    c\,
    \Big(
      \Norm{\nu^{\om}}{\frac{d}{2}, B(n)}\, \Norm{\mu^{\om}}{1, B(n)}
    \Big)^{\!\frac{1}{2}}\mspace{-8mu}. 
  \end{align*}
  This shows that the second hypothesis of Lemma~\ref{lem:bombieri_giusti} is satisfied by $f_1$.  Thus, we can apply Lemma~\ref{lem:bombieri_giusti} to $f_1$ and, since $U_\delta=B(n/2)$, we conclude that
  \begin{align}\label{eq:appl_bg1}
    \Norm{f_1}{\infty,B(n/2)} \;\leq\; A_1^\om 
    \qquad \Longleftrightarrow \qquad
    \me^{\xi} \;\leq\; A_1^{\om}\, \min_{x \in B(n/2)} u(x).
  \end{align}
  Analogously one can verify that the function $f_2 = u\, \me^{-\xi}$, with $\xi(u) \ldef (\ln u)_{B(n), \eta^2}$, satisfies the conditions of Lemma~\ref{lem:bombieri_giusti} with the choice $U$ and $U_\si$ as before, $\de = 3/4 \vr$ and any fixed $\alpha^*\in (0, p_*/2)$.  Indeed, the first hypothesis is fulfilled by Proposition~\ref{prop:EHI:mos_it_super} with the choice $\be = \al^*$ and the second condition again by the weighted local $\ell^2$-Poincar\'{e} inequality in Proposition~\ref{prop:weight_poinc}. Hence, $\Norm{f_2}{\al^*,B(3n/4)} \; \leq \;A_2^{\om}$ and by Corollary~\ref{cor:EHI:max_est2} we get 
  \begin{align}\label{eq:appl_bg2}
   \max_{x \in B(n/2)} u(x)
   \;\leq\;   
   \Big(
     16\, C_5(\al^*)\,
     \Big( 1 \vee \Norm{\mu^{\om}}{p,B(n)}\, \Norm{\nu^{\om}}{q,B(n)} \Big)
   \Big)^{\ka'}\; 
    A_2^{\om}\, \me^{\xi}.
  \end{align}
  By combining \eqref{eq:appl_bg1} and \eqref{eq:appl_bg2} we obtain the EHI. Finally, the explicit form of the Harnack constant follows from the representation of the constant in Lemma~\ref{lem:bombieri_giusti}.
\end{proof} 

\subsection{H\"{o}lder continuity}
As an immediate consequence from the elliptic Harnack inequality we prove the H\"older-continuity of harmonic functions.
\begin{prop} \label{prop:cont_harmonic}
  Suppose that Assumption \ref{ass:lim_const} holds. Let $x_0 \in V$ and let $s^{\om}(x_0)$ and $C^*_{\mathrm{EH}}$ be as in Assumption~\ref{ass:lim_const}.  Further, let $R \geq s(x_0)$ and  suppose that $u > 0$ is a harmonic function on $B(x_0, R_0)$ for some $R_0 \geq 2 R$. Then, for any $x_1, x_2\in B(x_0,R)$,
  \begin{align*}
    \big| u(x_1) - u(x_2) \big|
    \;\leq\;
    c \cdot \left( \frac{R}{R_0} \right)^{\!\th} \max_{B(x_0, R_0 / 2)} u,
  \end{align*}
  where $\th = \ln\big(2C^*_{\mathrm{EH}}/(2C^*_{\mathrm{EH}}-1)\big)/\ln 2$ and $c$ only depending on $C^*_{\mathrm{EH}}$.
\end{prop}
\begin{proof}
  The proof is based on arguments similar to the ones given in \cite[page~50]{SC97}.  Set $R_k = 2^{-k} R_0$ and let $B_k = B(x_0,R_k)$ and note that $B_{k+1} \subset B_k$.  For $x \in B_k$ we define the function $v_k$ through  
  \begin{align*}
    v_k(x)
    \;=\;
    \frac{u(x) - \min_{B_k}u}{\max_{B_k} u - \min_{B_k}u}
    \;\in\;
    [0,1],
  \end{align*}
  Note that $v_k$ is harmonic on $B_k$ and satisfies $\osc(v_k, B_k) = 1$. (Here the oscillation of a function $u$ on $A$ is defined by $\osc(u,A) = \max_A u - \min_A u$.)  Replacing $v_k$ by $1 - v_k$ if necessary we may assume that $\max_{B_{k+1}} v_k \geq 1/2$. 

  Now for every $k$ such that $R_{k}\geq s^\om(x_0)$, an application of Theorem \ref{thm:EHI} yields
  \begin{align*}
    \frac{1}{2} 
    \;\leq\;
    \max_{B_{k+1}} v_k
    \;\leq\; 
    C_{\mathrm{EH}}^*\, \min_{B_{k+1}} v_k. 
  \end{align*}
  Moreover, for such $k$ it follows that
  \begin{align*}
    \osc(u,B_{k+1})
    &\;=\;
    \frac{\max_{B_{k+1}} u - \min_{B_{k+1}} u}{\osc(u,B_k)}\; \osc(u, B_k)
    \\[1ex]
    &\;=\;
    \bigg(
      1 \,+\,
      \frac{\max_{B_{k+1}} u - \max_{B_k} u}{\osc(u, B_k)}
      \,-\, \min_{B_{k+1}} v_k
    \bigg)\; \osc(u, B_k).
  \end{align*}
  Hence,
  \begin{align*}
    \osc(u,B_{k+1})
    \;\leq\;
    (1-\de)\; \osc(u, B_k)
  \end{align*}
  with $\de = \big(2 C_{\mathrm{EH}}^*\big)^{-1}$. Now we choose $k_0$ such that $R_{k_0} \geq R > R_{k_0 + 1}$.  Then, by iterating the above inequality on the chain of balls $B_1 \supset B_2 \supset \cdots \supset B_{k_0}$ we obtain
  \begin{align*}
    \osc(u, B_{k_0})
    \;\leq\;
    (1 - \de)^{k_0 - 1}\; \osc(u, B_1).
  \end{align*}
  Note that $B(x_0,R)\subseteq B(x_0, R_{k_0})$ and $B_1 = B(x_0, R_0/2)$. Since $(1-\de)^{k_0} \leq c (R/R_0)^\th$, the claim follows.
\end{proof}

\section{Parabolic Harnack inequality} \label{sec:PHI}
In  this section we prove the parabolic Harnack inequality in Theorem~\ref{thm:PHI}.
\subsection{Mean value inequalities}
\begin{lemma} \label{lem:PHI:DF}
  Suppose that $Q = I \times B$, where $I = [s_1, s_2]$ is an interval and $B$ is a finite, connected subset of $V$.  Consider a function $\eta\!: V \to \bbR$ and a smooth function $\ze\!: \bbR \to \bbR$ with
  \begin{align*}
    \supp \eta &\;\subset\; B,
    &0 &\;\leq\; \eta \;\leq\; 1
    &\text{and}&
    &\eta &\;\equiv\; 0 \quad \text{on} \quad \partial B,
    \\
    \supp \ze &\;\subset\; I, 
    &0 &\;\leq\; \ze \;\leq\; 1
    &\text{and}&
    &\ze(s_1) &\;=\; 0.
  \end{align*}
  Further, let $u > 0$ be such that $\partial_t u -\cL_Y^\om\, u \leq 0$ on $Q$.  Then, for all $\al \geq 1$,
  \begin{align}
    \label{eq:PHI:DF}
    \int_{I} \ze(t)\; \frac{\cE_{\eta^2}^{\om}(u_t^{\al})}{|B|} \, \md t
    &\;\leq\;
    64\, \al^2\,
    \bigg( 
      \norm{\nabla \eta}{\infty}{\!E}^2 \,+\,
      \big\| \ze' \big\|_{\raisebox{-0ex}{$\scriptstyle L^{\raisebox{.1ex}{$\scriptscriptstyle \!\infty$}} (I)$}}
    \bigg)\; \int_{I}\; \Norm{u_t^{2 \al}}{1,B, \mu^{\om}}\, \md t
    \nonumber\\
  \end{align}
  \vspace{-5ex}

  \noindent and
  \begin{align}
    \label{eq:PHI:max}
    \max_{t \in I}
    \Big( \ze(t)\; \Norm{(\eta\, u_t^{\al})^2}{1, B, \mu^{\om}} \Big)
    &\;\leq\;
    64\, \al^2\,
    \bigg( 
      \norm{\nabla \eta}{\infty\!}{\!E}^2 \,+\,
      \big\| \ze' \big\|_{\raisebox{-0ex}{$\scriptstyle L^{\raisebox{.1ex}{$\scriptscriptstyle \!\infty$}} (I)$}}
    \bigg)\; \int_{I}\; \Norm{u_t^{2 \al}}{1 , B, \mu^{\om} }\, \md t.
  \end{align}
\end{lemma}
\begin{proof}
  Since $\partial_t u - \cL_Y^\om\, u \leq 0$ on $Q$ we have, for every $t \in I$,
  \begin{align*}
    \frac{1}{2 \al}\, \partial_t \scpr{\eta^2\!}{u_t^{2\al}}{V, \mu^{\om}}
    &\,=\,
    \scpr{\eta^2\, u_t^{2 \al - 1}\!}{\partial_t\, u_t}{V, \mu^{\om}}
    \nonumber\\[.5ex]
    &\,\leq\,
    \scpr{\eta^2\, u_t^{2 \al - 1}\!}{\cL_Y^{\om}\, u_t}{V,\mu^{\om}}
    \,=\,
    -\scpr{\nabla(\eta^2\, u_t^{2\al-1})}{\om\, \nabla u_t}{E}\!.
  \end{align*}
  By choosing $\ve = \frac{1}{8 \al}$ in \eqref{eq:DF:lb} and exploiting the fact that $\al \geq 1$, we obtain
  \begin{align*}
    \scpr{\nabla (\eta^2 u_t^{2\al-1})}{\om\, \nabla u_t}{E}
    \;\geq\;
    \frac{1}{2 \al}\; \cE_{\eta^2}^{\om}(u_t^{\al}) \,-\,
    32\, \al\; |B|\;\norm{\nabla \eta}{\infty\!}{\!E}^2\;
    \Norm{u_t^{2\al}}{1, B, \mu^{\om}}.
  \end{align*}
  Hence,
  \begin{align}\label{eq:PHI:DF:ode}
    \partial_t \Norm{(\eta\, u_t^{\al})^2}{1, B, \mu^{\om}} \,+\,
    \frac{\cE_{\eta^2}^{\om}(u_t^{\al})}{|B|}
    \;\leq\;
    64\, \al^2\;
    \norm{\nabla \eta}{\infty\!}{\!E}^2\;
    \Norm{u_t^{2 \al}}{1, B, \mu^{\om}}.
  \end{align}
  By multiplying both sides of \eqref{eq:PHI:DF:ode} with $\ze(t)$, we get
  \begin{align*}
    &\partial_t
    \Big( \ze(t)\, \Norm{(\eta\, u_t^{\al})^2}{1, B, \mu^{\om}} \Big)
    \,+\,
    \ze(t)\, \frac{\cE_{\eta^2}^{\om}(u_t^{\al})}{|B|}
    \nonumber\\[1ex]
    &\mspace{36mu}\leq\;
    64\, \al^2 \,
    \norm{\nabla \eta}{\infty\!}{\!E}^2\,
    \Norm{u_t^{2 \al}}{1, B, \mu^{\om}}
    \,+\,
    \big\|
      \ze'
    \big\|_{\raisebox{-0ex}{$\scriptstyle L^{\raisebox{.1ex}{$\scriptscriptstyle \!\infty$}} (I)$}}
    \Norm{(\eta\, u_t^{\al})^2}{1, B, \mu^{\om}}.
  \end{align*}
  Integrating this inequality over $[s_1,s]$ for any $s \in I$, we obtain
  \begin{align}\label{eq:PHI:DF:int}
    &\ze(s)\, \Norm{(\eta\, u_s^{\al})^2}{1, B, \mu^{\om}}
    \,+\,
    \int_{s_1}^s\! \ze(t)\; \frac{\cE_{\eta^2}^{\om}(u_t^{\al})}{|B|}\, \md t
    \nonumber\\[1ex]
    &\mspace{36mu}\leq\;
    64\, \al^2\,
    \Big( 
      \norm{\nabla \eta}{\infty\!}{\!E}^2 \,+\,
      \big\|
        \ze'
      \big\|_{\raisebox{-0ex}{$\scriptstyle L^{\raisebox{.1ex}{$\scriptscriptstyle \!\infty$}} (I)$}}
    \Big)\; \int_{I}\; \Norm{u_t^{2 \al}}{1, B, \mu^{\om}}\, \md t.
  \end{align}
  By neglecting the first term on the left-hand side of \eqref{eq:PHI:DF:int} we get \eqref{eq:PHI:DF}, whereas \eqref{eq:PHI:max} follows once we neglect the second term on the left-hand side of \eqref{eq:PHI:DF:int}.
\end{proof}
For any $x_0 \in V$, $t_0 \geq 0$ and $n \geq 1$, we write $Q(n) \equiv [t_0, t_0 + n^2] \times B(n)$ with $B(n) \equiv B(x_0, n)$.  Further, we consider a family of intervals $\{I_{\si} : \si \in [0, 1]\}$, i. e.\
\begin{align*}
  I_{\si}
  \;\ldef\;
  \Big[
    \si\, t_0 + \big( 1 - \si \big) s',
    \big( 1 - \si \big) s'' + \si \big( t_0 + n^2 \big)
  \Big]
\end{align*}
interpolating between the intervals $[t_0, t_0+n^2]$ and $[s', s'']$, where $s' \leq s''$ are chosen such that given $\ve \in (0, 1/4)$ we have $s' - t_0 \geq \ve n^2$ and either $t_0 + n^2 - s'' \geq \ve n^2$ or $s'' = t_0 + n^2$.  Moreover, set $ Q(\si n) \ldef I_{\si} \times B(\si n)$.  In addition, for any sets $I$ and $B$ as in Lemma~\ref{lem:PHI:DF}  let us introduce a $L^p$-norm on functions $u\!:\bbR \times V \to \bbR$ by
\begin{align*}
  \Norm{u}{p,I \times B, \mu^{\om}}
  \;\ldef\;
  \bigg(
    \frac{1}{|I|}\; \int_{I}\; \Norm{u_t}{p, B, \mu^{\om}}^p\; \md t
  \bigg)^{\!\!\frac{1}{p}},
\end{align*}
where $u_t=u(t,.)$, $t\in \bbR$.
\begin{prop}\label{prop:PHI:mos_it}
  Suppose that $s' - t_0 \geq \ve n^2$ and $t_0 + n^2 - s'' \geq \ve n^2$ for some $\ve \in (0, 1/4)$.  Let $u > 0$ be such that $\partial_t u - \cL_Y^{\om}\, u = 0$ on $Q(n)$.  Then, for any $p, q \in (1,\infty]$ with 
  \begin{align}\label{eq:PHI:cond:pq}
    \frac{1}{p} \,+\, \frac{1}{q} \;<\; \frac{2}{d}
  \end{align}
  there exists $\ka \equiv \ka(d,p,q)$ and $C_7 \equiv C_7(d, q)$ such that for all $\be \in [1, \infty]$ and $1/2 \leq \si' < \si \leq 1$, we have
  \begin{align}\label{eq:PHI:mos_it}
    \Norm{u}{2 \be, Q(\si' n),\mu^{\om}}
    \;\leq\;
    C_7\,
    \Bigg(
      \frac{
             \big(1 \vee \Norm{\mu^{\om}}{p,B(n)}\big)\,
             \big(1 \vee \Norm{\nu^{\om}}{q,B(n)}\big)
           }
           {\ve (\si - \si')^2}
    \Bigg)^{\!\!\ka}
    \Norm{u}{2, Q(\si n), \mu^{\om}}.
  \end{align}
\end{prop}
\begin{proof}
We proceed similarly as in the proof of Proposition~\ref{prop:EHI:mos_it}. Let $(B(\si_k n) : k \in \bbN_0)$ be a sequence of balls with radius $\si_k n$ centered at $x_0$, where
  \begin{align*}
    \si_k \;=\; \si' + 2^{-k} (\si - \si')
    \qquad \text{and} \qquad
    \tau_k \;=\; 2^{-k-1} (\si - \si'),
    \quad k \in \bbN_0.
  \end{align*}
  Choose $\al = 1+(\rho-p_*)/\rho$ and set $\al_k = \al^k$ for $k \in \bbN_0$.  Since $1/p + 1/q < 2/d$, $\al > 1$.  Again, we distinguish two cases.

  Consider first the case $\tau_k n \geq 1$.  Let $\eta_k$ be a cut-off function in space such that $\supp \eta_k \subset B(\si_k n)$, $\eta_k \equiv 1$ on $B(\si_{k+1} n)$, $\eta_k \equiv 0$ on $\partial B(\si_k n)$ and $\norm{\nabla \eta_k}{\infty}{E} \leq 1/\tau_{k} n$.  Further, let $\ze_k$ be a cut-off function in time, i.e.\ $\ze_k \in C^{\infty}(\bbR)$, $\supp \ze_k \subset I_{\si_k}$, $\ze_k \equiv 1$ on $I_{\si_{k+1}}$, $\ze_k(\si_k t_0 + (1-\si_k)s') = 0$ and $\big\| \ze_k' \big\|_{\raisebox{-0ex}{$\scriptstyle L^{\raisebox{.1ex}{$\scriptscriptstyle \!\infty$}} ([t_0,t_0+n^2])$}} \leq  1 / \ve \tau_k n^2$.

  Then, by H\"{o}lder's inequality, we have that
  \begin{align*}
    &\Norm{u_t^{2 \al_k}}
    {\al, B(\si_{k+1} n), \mu^{\om}}^{\al}
    \\[1ex]
    &\mspace{36mu}\leq\;
    \frac{|B(\si_k n)|}{|B(\si_{k+1} n)|}\;
    \Norm{(\eta_k\, u_t^{\al_k})^2}{1, B(\si_{k} n),\mu^{\om}}^{1-p_*/\rho}\;
    \Norm{(\eta_k\, u_t^{\al_k})^2}{\rho,B(\si_{k} n)}\;
    \Norm{\mu^{\om}}{p, B(\si_k n)}^{p_*/\rho}.
  \end{align*}
  Since by the volume regularity $|B(\si_k n)| / |B(\si_{k+1} n)| \leq C_{\mathrm{reg}}^2 2^d$, we obtain
  \begin{align*}
    \Norm{u^{2\al_k}}
    {\al, Q(\si_{k+1} n), \mu^{\om}}^{\al}
    &\;\leq\;
    C_{\mathrm{reg}}^2 2^d\, \Norm{\mu^{\om}}{p,B(\si_{k} n)}^{p_*/\rho}
    \bigg(
      \max_{t \in I_{\si_{k+1}}}
      \Norm{(\eta_k\, u_t^{\al_k})^2}
           {1, B(\si_{k} n), \mu^{\om}}^{1 - p_*/\rho}\!
    \bigg)
    \\
    &\mspace{36mu}\times\;
    \bigg(
      \frac{1}{|I_{\si_{k+1}}|}\,
      \int_{I_{\si_{k+1}}}\mspace{-15mu}
      \Norm{(\eta_k\, u_t^{\al_k})^2}{\rho,B(\si_{k} n)}\, \md t
    \bigg).
  \end{align*}
  Now, in view of \eqref{eq:sob:ineq:2}, the integrand can be estimated from above by
  \begin{align*}
    \Norm{(\eta_k\, u_t^{\al_k})^2}{\rho, B(\si_k n)}
    \;\leq\;
    C_{\mathrm{S}}\, n^2\, \Norm{\nu^{\om}}{q,B(\si_k n)}
    \Bigg(
      \frac{\cE_{\eta_k^2}^{\om}\big(u_t^{\al_k}\big)}{|B(\si_k n)|}
      \,+\,
      \frac{1}{(\tau_k n)^2}\,
      \Norm{u_t^{2 \al_k}}{1, B(\si_k n), \mu^{\om}}
    \Bigg).
  \end{align*}
  On the other hand, we obtain from \eqref{eq:PHI:DF} that
  \begin{align*}
    \int_{I_{\si_{k+1}}}\mspace{-8mu}
    \frac{\cE_{\eta_k^2}^{\om}\big(u_t^{\al_k}\big)}{|B(\si_k n)|}\; \md t
    \;\leq\;
    64\, \bigg(\frac{\al_k}{\tau_k n}\bigg)^{\!\!2}\,
    \Big(
      1 + \frac{\tau_k}{\ve}
    \Big)\,
    \int_{I_{\si_k}}\mspace{-6mu}
    \Norm{u_t^{2 \al_k}}{1,B(\si_k n), \mu^{\om}}\, \md t.
  \end{align*}
  Combining the estimates above and using \eqref{eq:PHI:max}, we finally get
  \begin{align}\label{eq:PHI:iter:est2}
    &\Norm{u}{2 \al_{k+1}, Q(\si_{k+1} n), \mu^{\om}}
    \nonumber\\[.5ex]
    &\mspace{36mu}\leq\;
    \bigg(
      c\,
      \frac{2^{2k}\, \al_k^2}{\ve (\si - \si')^2}\;
      \Big(
        1 \vee \Norm{\mu^{\om}}{p,B(n)}^{p_*/\rho}\,
        \Norm{\nu^{\om}}{q,B(n)}
      \Big)
    \bigg)^{\!\!1/(2 \al_k)}
    \Norm{u}{2 \al_k, Q(\si_k n), \mu^{\om}}
  \end{align}
  for some $c < \infty$, where we used that $\tau_k \leq 1$ and $|I_{\si_k}|/|I_{\si_{k+1}}| = \si_k/\si_{k+1}\leq 2$.

  Next we consider the case $\tau_k n < 1$ (i.e.\ $B(\si_{k+1} n) = B(\si_k n)$).  Then, for any $t \in I_{\si_{k+1}}$,
  \begin{align}\label{eq:PHI:est:fixed_t}
    \max_{x \in B(\si_k n)} u_t(x)^{2\al_k}
    &\;\leq\;
    |B(\si_k n)|^{(q+1)/q}\, \Norm{u_t^{2\al_k}}{q/(q+1), B(\si_k n)}
    \nonumber\\[.5ex]
    &\;\leq\;
    |B(\si_k n)|^{(q+1)/q}\, \Norm{u_t^{2\al_k}}{1, B(\si_k n), \mu^{\om}}\; 
    \Norm{\nu^{\om}}{q, B(\si_k n)},
  \end{align}
  where we used in the last step that $\mu^{\om}(x)^{-1} \leq \nu^{\om}(x)$.  The following argument is similar to the one given in \cite{De99}.  Since $u>0$ is caloric on $Q(n)$, we have for every $t \in I_{\si}$ and $x \in B(\si n)$ that $\partial_t u_t(x) = (\cL_Y^{\om} u_t)(x) \geq -u_t(x)$, which implies that $u_{t_2}(x) \geq \me^{-(t_2 - t_1)} u_{t_1}(x)$ for every $t_1, t_2 \in I_{\si}$ with $t_1 < t_2$ and $x \in B(\si n)$.  In particular, $\|u_{t_2}^{2\al_k}\|_{1, B(\si_k n), \mu^{\om}} \geq \me^{-2\al_k (t_2 - t_1)} \|u_{t_1}^{2\al_k}\|_{1, B(\si_k n), \mu^{\om}}$ for all $t_1, t_2 \in I_k$. Now choose $t_1 \in I_{k+1}$ in such a way that
  \begin{align*}
    \Norm{u_{t_1}^{2\al_k}}{1, B(\si_k n), \mu^{\om}}
    \;=\;
    \max_{t \in I_{k+1}}\, \Norm{u_t^{2\al_k}}{1, B(\si_k n), \mu^{\om}}
  \end{align*}
  and set $I_k' \ldef [t_1, (1-\si_k)s'' + \si_k(t_0 + n^2)]$.  Then, 
  \begin{align}  \label{eq:PHI:est:maxtime}
    \max_{t \in I_{k+1}}\, \Norm{u_t^{2\al_k}}{1, B(\si_k n), \mu^{\om}}\;
    &\;\leq\;
    \frac{2\al_k|I_k|}{1 - \me^{-2\al_k|I'_k|}}\, 
    \Norm{u^{2\al_k}}{1,Q(\si_k n), \mu^{\om}} \nonumber
    \\[.5ex]
    &\;\leq\;
    \Big(\frac{1}{\ve \tau_k} + 2\al_kn^2\Big)\, 
    \Norm{u^{2 \al_k}}{1,Q(\si_k n), \mu^{\om}},
  \end{align}
  where we used that $|I_k| \leq n^2$, $|I_k'| \geq \ve \tau_k n^2$ and $x / (1 - \me^{-x}) \leq 1+x$.  Since $n < 1/\tau_k$ and $(1+d(1+1/q)/2)/\al_* \leq 1$, we combine \eqref{eq:PHI:est:fixed_t} and \eqref{eq:PHI:est:maxtime} to obtain
  \begin{align*}
    &\Norm{u^{2\al_k}}{\al, Q(\si_{k+1} n), \mu^{\om}}
    \\[.5ex]
    &\mspace{36mu}\leq\; 
    \bigg(
      \max_{(t,x) \in Q(\si_{k+1} n)} u_t(x)^{2\al_k}
    \bigg)^{\!\!1/\al_*}\,
    \Norm{u^{2\al_k}}{1,Q(\si_k n), \mu^{\om}}^{1/\al}  \\[.5ex]
    &\mspace{36mu}\leq\;
    \bigg(
      |B(\si_k n)|^{(q+1)/q}\, \Norm{\nu}{q, B(\si_k n)}
      \max_{t \in I_{k+1}} \Norm{u_t^{2\al_k}}{1, B(\si_k n), \mu^{\om}}
    \bigg)^{\!\!1/\al_*}\,
    \Norm{u^{2\al_k}}{1,Q(\si_k n), \mu^{\om}}^{1/\al}
    \\[.5ex]
    &\mspace{36mu}\leq\;
    c\, \al_k\, \frac{1 \vee \Norm{\nu^{\om}}{q, B(n)}}{\ve \tau_k^2}\,
    \Norm{u^{2\al_k}}{1,Q(\si_k n), \mu^{\om}}.
  \end{align*}
  Hence, 
  \begin{align}\label{eq:PHI:iter:est1}
    \Norm{u}{2\al_{k+1}, Q(\si_{k+1} n), \mu^{\om}}
    \;\leq\;  
    \bigg(
      c\, \frac{2^{2k}\, \al_k}{\ve (\si-\si')^2}\, 
      \Big( 1 \vee \Norm{\nu^{\om}}{q, B(n)}\Big)
    \bigg)^{\!\!1/(2\al_k)}\;
    \Norm{u}{2\al_k, Q(\si_k n), \mu^{\om}}.
  \end{align}

  By iterating the inequalities \eqref{eq:PHI:iter:est1} and \eqref{eq:PHI:iter:est2}, respectively, and using the fact that $\sum_{k=0}^{\infty} k / \al_k < \infty$ and $\sup_k 2^{(k+1)/2 \alpha_k} < \infty$, there exists $c < \infty$ independent of $n$ and $K$ such that, for any $K \in \bbN$,
  \begin{align}\label{eq:PHI:max:est:0}
    \Norm{u}{2 \al_{K}, Q(\si_K n), \mu^{\om}}
    \;\leq\;
    c
    \prod_{k=0}^{K-1}
    \Bigg(
      \frac{
             \big(1 \vee \Norm{\mu^{\om}}{p,B(n)}\big)\,
             \big(1 \vee \Norm{\nu^{\om}}{q,B(n)}\big)
           }
           {\ve (\si - \si')^2}
    \Bigg)^{\!\!1/2 \al_k}\;
    \Norm{u}{2, Q(\si n), \mu^{\om}},
  \end{align}
  Setting $\ka \ldef \frac{1}{2} \sum_{k=0}^{\infty} 1/\al_k < \infty$ and using that $Q(\si_K n) \downarrow Q(\si' n)$, we get
  \begin{align}\label{eq:PHI:max:est:1}
    \max_{(t,x) \in Q(\si' n)} u(t,x)
    &\;=\;
    \lim_{K \to \infty} \Norm{u}{2 \al_{K}, Q(\si_K n), \mu^{\om}}
    \nonumber\\[.5ex]
    &\;\leq\;
    c
    \Bigg(
      \frac{
             \big(1 \vee \Norm{\mu^{\om}}{p,B(n)}\big)\,
             \big(1 \vee \Norm{\nu^{\om}}{q,B(n)}\big)
           }
           {\ve (\si - \si')^2}
    \Bigg)^{\!\!\ka}\;
    \Norm{u}{2, Q(\si n), \mu^{\om}}
  \end{align}
  Finally, for any $\be \in [1, \infty)$, choose $K < \infty$ such that $\be \in[\al_K, \al_{K+1})$.   Since
  \begin{align*}
    \Norm{u}{2 \be, Q(\si' n), \mu^{\om}} 
    \;\leq\; 
    \Norm{u}{2\al_K, Q(\si' n), \mu^{\om}}^{\al_K / \be}\, 
    \bigg( \max_{(t,x) \in Q(\si' n)} u(t,x) \bigg)^{\!\!1 - \al_K / \be}
  \end{align*}
  and $\|u\|_{2 \al_K, Q(\si'n), \mu^{\om}} \leq c \|u\|_{2 \al_K, Q(\si_K n) ,\mu^{\om}}$, the assertion follows \eqref{eq:PHI:max:est:0} and \eqref{eq:PHI:max:est:1}.
\end{proof}
\begin{corro}\label{cor:PHI:max_est}
  Let $u > 0$ be such that $\partial_t u - \cL^{\om} u = 0$ on $Q(n)$.  Further, suppose that $p, q \in (1, \infty]$ satisfy the condition \eqref{eq:PHI:cond:pq}.
  \begin{enumerate}[(i)]
    \item Suppose that $s' - t_0 \geq \ve n^2$ and either $t_0 + n^2 - s'' \geq \ve n^2$ or $s'' = t_0 + n^2$ for some $\ve \in (0, 1/4)$.  Then, for all $\al \in (0, \infty)$, there exist $C_8 = C_8(d, q)$ and $\ka' = 2 \ka / \al$ such that
      \begin{align}\label{eq:PHI:max_est1}
        \max_{(t,x) \in Q(\si' n)}\! u(t,x)^{-1}
        \;\leq\;
        \Bigg(
          C_8\,
          \frac{
            \big(1 \vee \Norm{\mu^{\om}}{p,B(n)}\big)\,
            \big(1 \vee \Norm{\nu^{\om}}{q,B(n)}\big)
          }
          {\ve (\si - \si')^2}
        \Bigg)^{\!\!\ka'}
        \Norm{u^{-1}}{\al, Q(\si n), \mu^{\om}}.
      \end{align}
    \item Suppose that $s' - t_0 \geq \ve n^2$ and $t_0 + n^2 - s'' \geq \ve n^2$ for some $\ve \in (0, 1/4)$.  Then, for all $\al \in (0, \infty)$, there exist $C_9(\al) = C_9(d, q, \al)$ and $\ka' = (1 \vee \frac{2}{\al})(1 + \ka)$ such that
      \begin{align}\label{eq:PHI:max_est2}
        \max_{(t,x) \in Q(\si' n)} u(t,x)
        \;\leq\;
        \Bigg(
          C_9(\al)\,
          \frac{
            \big(1 \vee \Norm{\mu^{\om}}{p, B(n)}\big)\,
            \big(1 \vee \Norm{\nu^{\om}}{q, B(n)}\big)
          }
          {\ve (\si - \si')^2}
        \Bigg)^{\!\!\ka'}
        \Norm{u}{\al, Q(\si n), \mu^{\om}}.
      \end{align}
  \end{enumerate}
\end{corro}
\begin{proof}
  (i) Since $u > 0$ is caloric on $Q(n)$, the function $f = u^{-\ga}$, $\ga \in (0, \infty)$, is sub-caloric on $Q(n)$, that is $\partial_t f - \cL_{Y}^{\om} f \leq 0$.  Let $\si_k, \tau_k$ and $\al_k$ be defined as in the proof of Proposition~\ref{prop:PHI:mos_it}.  Again, we have to distinguish two different cases.  First consider the case $\tau_k n \geq 1$.  Since the energy estimate (Lemma~\ref{lem:PHI:DF}) is valid for sub-caloric functions, by following literally the corresponding argument in the proof of Proposition~\ref{prop:PHI:mos_it}, we obtain that
  \begin{align}\label{eq:PHI:iter:corro1}
    &\Norm{u^{-\ga}}{2 \al_{k+1}, Q(\si_{k+1} n), \mu^{\om}}
    \nonumber\\[.5ex]
    &\mspace{36mu}\leq\;
    \bigg(
      c\,
      \frac{2^{2k}\, \al_k^2}{\ve (\si - \si')^2}\;
      \Big(
        1 \vee
        \Norm{\mu^{\om}}{p,B(n)}^{p_*/\rho}\,
        \Norm{\nu^{\om}}{q,B(n)}
      \Big)
    \bigg)^{\!\!1/(2 \al_k)}
    \Norm{u^{-\ga}}{2 \al_k, Q(\si_k n), \mu^{\om}}.
  \end{align}
  Next consider the case $\tau_k n < 1$.  Since $u > 0$ is caloric on $Q(n)$, we have for every $t \in I_{\si}$ and $x \in B(\si n)$ that $\partial_t u_t(x) = (\cL_Y^{\om} u_t)(x) \geq -u_t(x)$, which implies that $u_{t_1}^{-1}(x) \geq \me^{-(t_2 - t_1)} u_{t_2}^{-1}(x)$ for every $t_1, t_2 \in I_{\si}$ with $t_1 < t_2$ and $x \in B(\si n)$.  In particular, $\|u_{t_1}^{-2\ga\al_k}\|_{1, B(\si_k n), \mu^{\om}} \geq \me^{-2 \ga \al_k (t_2 - t_1)} \|u_{t_2}^{2\al_k}\|_{1, B(\si_k n), \mu^{\om}}$ for all $t_1, t_2 \in I_k$.  By choosing $t_2 \in I_{k+1}$ in such a way that
  \begin{align*}
    \Norm{u_{t_2}^{-2 \ga \al_k}}{1, B(\si_k n), \mu^{\om}}
    \;=\;
    \max_{t \in I_{k+1}}\, \Norm{u_t^{-2 \ga \al_k}}{1, B(\si_k n), \mu^{\om}}
  \end{align*}
  and setting $I_k'' \ldef [\si_k t_0 + (1 - \si_k)s', t_2]$, we obtain that
  \begin{align*}
    \max_{t \in I_{k+1}}\, \Norm{u_t^{-2 \ga \al_k}}{1, B(\si_k n), \mu^{\om}}\;
    &\;\leq\;
    \frac{2 \ga \al_k|I_k|}{1 - \me^{-2 \ga \al_k|I''_k|}}\, 
    \Norm{u^{-2 \ga \al_k}}{1, Q(\si_k n), \mu^{\om}}.
  \end{align*}
  Since $|I_k''| \geq \ve \tau_k n^2$, following the corresponding argument in the proof of Proposition~\ref{prop:PHI:mos_it} we obtain that
  \begin{align}\label{eq:PHI:iter:corro2}
    &\Norm{u^{-\ga}}{2 \al_{k+1}, Q(\si_{k+1} n), \mu^{\om}}
    \nonumber\\  
    &\mspace{36mu}\leq\;  
    \bigg(
      c\, \frac{2^{2k}\, (1 \vee \ga) \al_k}{\ve (\si-\si')^2}\, 
      \Big( 1 \vee \Norm{\nu^{\om}}{q, B(n)}\Big)
    \bigg)^{\!\!1/(2 \al_k)}\;
    \Norm{u^{-\ga}}{2 \al_k, Q(\si_k n), \mu^{\om}}.
  \end{align}
  Thus, by iterating the estimates \eqref{eq:PHI:iter:corro1} and \eqref{eq:PHI:iter:corro2}, respectively and using the fact that $\sum_{k=0}^{\infty}k/\al_k < \infty$, $\sup_{k} 2^{(k+1)/(2\ga\al_k)} < \infty$ and $Q(\si_k n) \downarrow Q(\si' n)$, we get that there exists $C_8(d,q) < \infty$ which is independent of $\ga$ such that
  \begin{align*}
    \max_{(t, x) \in Q(\si' n)} u^{-1}(t,x)
    &\;=\;
    \lim_{K \to \infty} \Norm{u^{-1}}{2 \ga \al_K, Q(\si_K n), \mu^{\om}}
    \\[.5ex]
    &\;=\;
    \Bigg(
      C_8\,
      \frac{
       \big(1 \vee \Norm{\mu^{\om}}{p,B(n)}\big)\,
        \big(1 \vee \Norm{\nu^{\om}}{q,B(n)}\big)
      }
      {\ve (\si - \si')^2}
    \Bigg)^{\!\!\ka/\ga}\!
    \Norm{u^{-1}}{2 \ga, Q(\si n), \mu^{\om}}, 
  \end{align*}
  where $\ka = \frac{1}{2} \sum_{k=0}^{\infty} 1/\al_k$.  Since $\ga \in (0, \infty)$ was arbitrary, the claim follows.
  
  (ii) For every $\al \geq 2$, the inequality \eqref{eq:PHI:max_est2} follows directly from \eqref{eq:PHI:mos_it} with $\be = \infty$ and the fact that $\|u\|_{2, Q(\si n), \mu^{\om}} \leq (1 \vee \|\mu^{\om}\|_{1, B(n)})\, \|u\|_{\al, Q(\si n), \mu^{\om}}$.  For $\al \in (0, 2)$ the proof goes literally along the lines of the proof of Corollary~\ref{cor:EHI:max_est2}.
\end{proof}
\begin{lemma} \label{lem:PHI:DF:super}
  Suppose that $Q = I \times B$, where $I = [s_1, s_2]$ is an interval and $B$ is a finite, connected subset of $V$.  Consider a function $\eta\!: V \to \bbR$ and a smooth function $\ze\!: \bbR \to \bbR$ with
  \begin{align*}
    \supp \eta &\;\subset\; B,
    &0 &\;\leq\; \eta \;\leq\; 1
    &\text{and}&
    &\eta &\;\equiv\; 0 \quad \text{on} \quad \partial B,
    \\
    \supp \ze &\;\subset\; I, 
    &0 &\;\leq\; \ze \;\leq\; 1
    &\text{and}&
    &\ze(s_2) &\;=\; 0.
  \end{align*}
  Further, let $u > 0$ be such that $\partial_t u -\cL_Y^\om\, u \geq 0$ on $Q$.  Then, for all $\al \in (0,1/2)$,
  \begin{align}
    \label{eq:PHI:DF:super}
    &\int_{I} \ze(t)\; \frac{\cE_{\eta^2}^{\om}(u_t^{\al})}{|B|} \, \md t
    \nonumber\\[.5ex]
    &\mspace{36mu}\leq\;
    \frac{8}{(1 - 2\al)^2}\,
    \bigg(
      \osr\big(\eta^2\big)\,
      \norm{\nabla \eta}{\infty}{\!E}^2 \,+\,
      \big\| \ze' \big\|_{\raisebox{-0ex}{$\scriptstyle L^{\raisebox{.1ex}{$\scriptscriptstyle \!\infty$}} (I)$}}
    \bigg)\; \int_{I}\; \Norm{u_t^{2 \al}}{1,B, \mu^{\om}}\, \md t
    \nonumber\\
  \end{align}
  \vspace{-5ex}

  \noindent and
  \begin{align}
    \label{eq:PHI:max:super}
    &\max_{t \in I}
    \Big( \ze(t)\; \Norm{(\eta\, u_t^{\al})^2}{1, B, \mu^{\om}} \Big)
    \nonumber\\[.5ex]
    &\mspace{36mu}\leq\;
    \frac{16}{1 - 2\al}\,
    \bigg(
      \osr\big(\eta^2\big)\,
      \norm{\nabla \eta}{\infty\!}{\!E}^2 \,+\,
      \big\| \ze' \big\|_{\raisebox{-0ex}{$\scriptstyle L^{\raisebox{.1ex}{$\scriptscriptstyle \!\infty$}} (I)$}}
    \bigg)\; \int_{I}\; \Norm{u_t^{2 \al}}{1 , B, \mu^{\om} }\, \md t.
  \end{align}
\end{lemma}
\begin{proof}
  Since $\partial_t u - \cL_Y^\om\, u \geq 0$ on $Q$ we have, for every $t \in I$ and $\al \in (0, 1/2)$,
  \begin{align*}
    \frac{1}{2 \al}\, \partial_t \scpr{\eta^2\!}{u_t^{2\al}}{V, \mu^{\om}}
    &\,=\,
    \scpr{\eta^2\, u_t^{2 \al - 1}\!}{\partial_t\, u_t}{V, \mu^{\om}}
    \nonumber\\[.5ex]
    &\,\geq\,
    \scpr{\eta^2\, u_t^{2 \al - 1}\!}{\cL_Y^{\om}\, u_t}{V,\mu^{\om}}
    \,=\,
    -\scpr{\nabla(\eta^2\, u_t^{2\al-1})}{\om\, \nabla u_t}{E}\!.
  \end{align*}
  Together with the estimate \eqref{eq:DFsuper:est1}, we obtain that
  \begin{align}\label{eq:PHIsuper:DF:ode}
    -\partial_t \Norm{(\eta\, u_t^{\al})^2}{1, B, \mu^{\om}} \,+\,
    \frac{1-2\al}{\al}\,
    \frac{\cE_{\eta^2}^{\om}(u_t^{\al})}{|B|}
    \;\leq\;
    \frac{32\, \al}{1-2\al}\;
    \osr(\eta^2)\;
    \norm{\nabla \eta}{\infty\!}{\!E}^2\;
    \Norm{u_t^{2 \al}}{1, B, \mu^{\om}}.
  \end{align}
  The difference between \eqref{eq:PHI:DF:ode} and \eqref{eq:PHIsuper:DF:ode} is the minus sign that appears in front of the first term on the left-hand side of \eqref{eq:PHIsuper:DF:ode}.  This is the reason why we choose a cut-off function, $\zeta$, that vanishes at $s_2$.  To finish the proof, it suffices to repeat the remaining arguments in Lemma~\ref{lem:PHI:DF}.
\end{proof}
\begin{prop} \label{prop:PHIsuper:mos_it}
  Let $u > 0$ be such that $\partial_t u - \cL_Y^{\om}\, u \geq 0$ on $Q(n)$.  Further, suppose that $s' - t_0 \geq \ve n^2$ and $t_0 + n^2 - s'' \geq \ve n^2$ for some $\ve \in (0, 1/4)$.  Then, for any $p, q \in (1,\infty]$ satisfying \eqref{eq:PHI:cond:pq} and any $\be \in (0, 1/2)$, there exists $C_{10} \equiv C_{10}(d, q)$ and $\ka \equiv \ka(d,p,q)$ such that for all $\al \in (0, \be / 2)$ and for all $1/2 \leq \si' < \si \leq 1$, we have
  \begin{align}\label{eq:PHIsuper:mos_it}
    \Norm{u}{\be, Q(\si' n),\mu^{\om}}
    \;\leq\;
    \Bigg(
      C_{10}\,
      \frac{
             \big(1 \vee \Norm{\mu^{\om}}{p,B(n)}\big)\,
             \big(1 \vee \Norm{\nu^{\om}}{q,B(n)}\big)
           }
           {\ve (\si - \si')^2}
    \Bigg)^{\!\!\ka\,(\frac{1}{\al} - \frac{1}{\be})}
    \Norm{u}{\al, Q(\si n), \mu^{\om}}.
  \end{align}
\end{prop}
\begin{proof}
  We proceed similar as in \cite[Theorem~5.2.17]{SC97} (c.f.\ Proposition~\ref{prop:EHI:mos_it_super}).  Recall that $\rho > p_*$ since $1/p + 1/q < 2/d$ and therefore $1+(\rho-p_*)/\rho>1$.  W.l.o.g.\ let us assume that $1 + (\rho - p_*) / \rho \leq 2$ (otherwise we replace $\rho$ by $\rho' \ldef 2p_*$ in the following argument).  Then, there exists $j \geq 2$ such that
  \begin{align*}
    \be\, \bigg(1 +  \frac{\rho - p_*}{\rho} \bigg)^{\!\!-j}
    \;\leq\;
    \al
    \;<\;
    \be\, \bigg(1 + \frac{\rho - p_*}{\rho} \bigg)^{\!\!-(j-1)}.
  \end{align*}
  Further, set
  \begin{align*}
    \al_k
    \;\ldef\;
    \frac{\be}{2} \bigg(1 +  \frac{\rho -p_*}{\rho} \bigg)^{\!\!k-j}
    \leq\;
    \frac{\be}{2} 
    \;<\;
    \frac{1}{4},
    \qquad \forall\, k = 0, \ldots, j.
  \end{align*}
  Let $\si_k$, $\tau_k$, $\eta_k$ and $\ze_k$ be defined as in the proof of Proposition~\ref{prop:PHI:mos_it} above.  Again, we have to distinguish two different cases. In case $\tau_k n < 1$, note that $\partial_t u_t(x) \geq \big(\cL_Y^{\om} u_t\big)(x) \geq - u_t(x)$ for all $(t, x) \in Q(n)$.  Hence, by following literally the argument in the proof of Proposition~\ref{prop:PHI:mos_it}, we get
  \begin{align}\label{eq:PHI:iter_super:0}
    \Norm{u}{2\al_{k+1}, Q(\si_{k+1} n), \mu^{\om}}
    \;\leq\;  
    \bigg(
      c\, \frac{2^{2k} \al_k}{\ve (\si-\si')^2}\, 
      \Big( 1 \vee \Norm{\nu^{\om}}{q, B(n)}\Big)
    \bigg)^{\!\!1/(2\al_k)}\;
    \Norm{u}{2\al_k, Q(\si_k n), \mu^{\om}}.
  \end{align}
  Next consider the case $\tau_k n \geq 1$.  Note that $\osr(\eta^2_k) \leq 4$ for all $k$ such that $\tau_k n \geq 1$.  As in the proof of Proposition~\ref{prop:PHI:mos_it} we obtain by using the Sobolev inequality, Lemma~\ref{lem:PHI:DF:super} and the volume regularity that
  \begin{align}\label{eq:PHI:iter_super}
    &\Norm{u}{2 \al_{k+1}, Q(\si_{k+1} n), \mu^{\om}}
    \nonumber\\[.5ex]
    &\mspace{36mu}\leq\;
    \bigg(
      c\,
      \frac{2^{2k} \al_k^2}{\ve (\si -\si')^{2}}\;
      \Big(
        1 \vee \Norm{\mu^\om}{p,B(n)}^{p_*/\rho}\,
        \Norm{\nu^\om}{q,B(n)}
      \Big)
    \bigg)^{\!\!1/(2\al_k)}
    \Norm{u}{2 \al_k, Q(\si_{k} n), \mu^{\om}}.\!
  \end{align}
  Further, since $2 \al_0 \leq \al$, H\"older's inequality implies that
  \begin{align*}
    \Norm{u}{2 \al_0, Q(\si n), \mu^{\om}}
    \;\leq\;
    \Norm{\mu^{\om}}{1, B(\si n)}^{1/(2\al_0) - 1/\al}\,
    \Norm{u}{\al, Q(\si n), \mu^{\om}}.
  \end{align*}
  Thus, by iterating \eqref{eq:PHI:iter_super:0} and \eqref{eq:PHI:iter_super},respectively, and using the fact that $\al_k < 1/4$ for all $k = 0, \ldots, j$ and $\sum_{k=0}^{j-1} (k/\al_k) \leq c \sum_{k=0}^{j-1}1 /\al_k < \infty$ for some $c=c(d,p,q)$, there exists $C_{10} < \infty$ such that
  \begin{align*}
    \Norm{u}{\be, Q(\si' n), \mu^{\om}}
    \;\leq\;
     \prod_{k=0}^{j-1}
    \Bigg( C_{10}\,
      \frac{\big(1 \vee \Norm{\mu^\om}{p,B(n)}\big)\,
        \big(1 \vee \Norm{\nu^\om}{q,B(n)}\big)}{(\si -\si')^{2}}
    \Bigg)^{\!\!1/\al_k}
    \Norm{u}{\al, Q(\si n), \mu^{\om}}.
  \end{align*}
  Since $1 + (\rho - p_*)/\rho \geq 1$ and $j \geq 2$ we have that
  \begin{align*}
    \frac{1}{\beta}\, \Big(1 + \frac{\rho - p_*}{\rho} \Big)^{\!j}
    \,-\, \frac{1}{\be}
    \;\leq\;
    \Big( 2 + \frac{\rho - p_*}{\rho} \Big)\,
    \bigg( \frac{1}{\al} - \frac{1}{\be} \bigg).
  \end{align*}
  Hence,
  \begin{align*}
    \sum_{k=0}^{j-1} \frac{1}{\al_k}
    \;=\;
    2\, \frac{2\rho - p_*}{\rho - p_*}\,
    \bigg(
      \frac{1}{\beta}\, \Big(1 + \frac{\rho - p_*}{\rho} \Big)^{\!j}
      \,-\, \frac{1}{\be}
    \bigg)
    \;\leq\;
    \frac{ 2\, (2\rho-p_*)(3\rho-p_*)}{\rho(\rho-p_*)}\,
    \bigg( \frac{1}{\al} - \frac{1}{\be} \bigg),
  \end{align*}
  and the assertion follows.
\end{proof}

\subsection{From mean value inequalities to Harnack principle}
From now on we denote by $m$ the product measure $m = \md t \times | \cdot |$ on $\bbR_+ \times V$, where $| \cdot |$ is still the counting measure on $V$.
\begin{lemma} \label{lem:PHI:log}
  For $x_0 \in V$, $t_0 \geq 0$ and $n \geq 1$, let $Q(n) = I_1 \times B(n)$ where $I_1 = [t_0, t_0 + n^2]$ and $B(n) \equiv B(x_0, n)$.  Further, let $u > 0$ be such that $\partial_t u -\cL_Y^\om\, u \geq 0$ on $Q(n)$.  Then, there exists $C_{11}=C_{11}(d)$ and $\xi = \xi(u)$ such that for all $\la > 0$ and all $\varrho \in [1/2, 1)$,
  \begin{align*}
    m
    \Big[
      \big\{ (t,z) \in I_1 \times B(\varrho n) :\,
      \ln u(t, z) \,<\, -\la \,-\, \xi \big\}
    \Big]
    \;\leq\;
    C_{11}\,
    M_3\, M_4\,
    m\big[I_1 \times B(\varrho n)\big]\, \la^{-1},
  \end{align*}
  and
  \begin{align*}
    m
    \Big[
      \big\{ (t,z) \in I_1 \times B(\varrho n) :\,
      \ln u(t,z) \,>\, \la \,-\, \xi \big\}
    \Big]
    \;\leq\;
    C_{11}\,
    M_3 \, M_4\,
    m\big[I_1 \times B(\varrho n)\big] \, \la^{-1},
  \end{align*}
  where $M_3 \ldef \Norm{\mu^{\om}}{1, B(n)} /\Norm{\mu^{\om}}{1, B(n/2)}$ and $M_4 \ldef 1 \vee \Norm{\mu^{\om}}{p, B(n)}^2\, \Norm{\nu^{\om}}{q, B(n)} $.
\end{lemma}
\begin{proof}
  Let $\eta$ be defined as
  \begin{align*}
    \eta(x)
    \;=\;
    \bigg[ 1 - \frac{d(x_0,x)}{n} \bigg]_+.
  \end{align*}
  In particular, $\eta \leq 1$, $\max_{x\sim y} |\nabla\eta(x,y)| \leq 1/n$ and $\osr(\eta^2) \leq 4$. Moreover, note that $\eta \geq 1 - \varrho$ on $B(\varrho n)$, and by the volume regularity of $V$ we get
  \begin{align}\label{eq:comp_scpr_eta}
    |B(n)|\,  \Norm{\mu^{\om}}{1, B(n)}
    &\;\geq\; 
    \scpr{\eta^2}{1}{V,\mu^\om}
    \nonumber\\[.5ex]
    &\;\geq\;
    \big(1 - \varrho)^2\, \mu^{\om}\big[B(\varrho n)\big]
    \;\geq\;
    \frac{(1 - \vr)^2}{2^{d}\, C_{\mathrm{reg}}^{2}}\,
    |B(n)|\, \Norm{\mu^{\om}}{1, B(n/2)}.
  \end{align}
  Since $\partial_t u -\cL_Y^\om\, u \geq 0$ on $Q(n)$, we have 
  \begin{align*}
    \partial_t \scpr{\eta^2}{-\ln u_t}{V, \mu^{\om}}
    &\;=\;
    \scpr{\eta^2 u_t^{-1}}{-\partial_t u_t}{V, \mu^{\om}}
    \nonumber\\[.5ex]
    &\;\leq\;
    \scpr{\eta^2 u_t^{-1}}{-\cL_Y^{\om}\, u_t}{V, \mu^{\om}}
    \;=\;
    \scpr{\nabla \big(\eta^2 u_t^{-1}\big)}{\om\, \nabla u_t}{E}.
  \end{align*}
  In view of \eqref{eq:EHI:log:est1} and \eqref{eq:comp_scpr_eta}, we obtain
  \begin{align*}
    \scpr{\nabla \big(\eta^2 u_t^{-1}\big)}{\om\, \nabla u_t}{E}
    \;\leq\;
    -\cE^{\om,\eta}(\ln u_t) \,+\,
    \frac{c_1}{n^2}\, \scpr{\eta^2}{1}{V, \mu^\om} \, M_3.
  \end{align*}
  Hence,
  \begin{align}\label{eq:PHI:log:estimate}
    \partial_t \scpr{\eta^2}{-\ln u_t}{V, \mu^{\om}}
    \;\leq\;
    -\cE^{\om,\eta}(\ln u_t) \,+\,
    \frac{c_1}{n^2}\, \scpr{\eta^2}{1}{V, \mu^\om} \, M_3.
  \end{align}

  From now on the proof goes literally along the lines of \cite[Lemma 5.4.1]{SC02}.  For the reader's convenience we provide the details here.  Let us define
  \begin{align*}
    w(t,x)
    &\;\ldef\;
    -\ln u(t,x),
    \qquad
    &\overline w(t,x)
    &\;\ldef\;
    w(t,x) - \frac{c_1}{n^2}\, M_3 \, (t - t_0),
    \\
    W(t)
    &\;\ldef\;
    \frac{\scpr{\eta^2}{w_t}{V, \mu^\om}}{\scpr{\eta^2}{1}{V,\mu^\om}},
    &\overline W(t)
    &\;\ldef\;
    W(t) - \frac{c_1}{n^2}\, M_3 \, (t - t_0).
  \end{align*}
  With the notation introduced above, the inequality \eqref{eq:PHI:log:estimate} reads  
  \begin{align*}
  \partial_t W(t) \,+\,  \scpr{\eta^2}{1}{V,\mu^\om}^{-1}\,
  \cE^{\om,\eta}(w_t)
  \;\leq\;
  \frac{c_1}{ n^2}\, M_3.
  \end{align*}
  Hence, applying the local $\ell^2$-Poincar\'e inequality with radial cutoff given in Proposition~\ref{prop:weight_poinc}, the left hand side can be estimated from below by
  \begin{align*}
    \partial_t W(t) \,+\,
    \Big(
      C_{\mathrm{PI}}\, M_2\, n^2\,
      \Norm{\mu^{\om}}{p, B(n)}^2\, \Norm{\nu^{\om}}{q, B(n)}
    \Big)^{-1}\,
    \Norm{w_t - W(t)}{2, B(n), \eta^2\mu^\om}^2 ,
  \end{align*}
  where we used that $|B(n)|/\scpr{\eta^2}{1}{V, \mu^\om} \geq  \Norm{\mu^{\om}}{1, B(n)}^{-1} \geq  \Norm{\mu^{\om}}{p, B(n)}^{-1}$ by \eqref{eq:comp_scpr_eta} and Jensen's inequality.  This implies
  \begin{align} \label{eq:est:barW}
    \partial_t \overline W(t) \,+\,
    \Big(c_2\, n^2\, M_3 \, M_4 \Big)^{-1}\,
    \Norm{\overline{w}_t -\overline{W}(t)}{2, B(\varrho n)}^2
    \;\leq\;
    0.
  \end{align}
  In particular, $\partial_t \overline{W}(t) \leq 0$.  Let $\xi = \xi(u) = \overline{W}(t_0)$ and for any $\la > 0$ we define
  \begin{align*}
    D_t^+(\la)
    \;=\;
    \big\{ z \in B(\varrho n) :\, \overline{w}(t,z) > \xi + \la \big\}.
  \end{align*}
  Then, we have for all $t \in (t_0, t_0 + n^2]$ and $z \in D_t^+(\la)$ that
  \begin{align*}
    \overline{w}(t,z) \,-\, \overline{W}(t)
    \;\geq\;
    \la \,+\, \xi \,-\, \overline{W}(t)
    \;>\;
    \la,
  \end{align*}
  because $\xi = \overline{W}(t_0)$ and $\partial_t \overline{W}(t) \leq 0$. Using this in \eqref{eq:est:barW} we obtain
  \begin{align*}
    \partial_t \overline{W}(t) \,+\, 
    \Big(
      c_2\, n^2\, |B(\varrho n)|\, M_3 \, M_4
    \Big)^{-1}\,
    \big| \la + \xi - \overline{W}(t) \big|^2\; \big| D_t^+(\la) \big|
    \;\leq\;
    0
  \end{align*}
  which is equivalent to
  \begin{align*}
    \big| D_t^+(\la) \big|
    \;\leq\;
    - c_2\, n^2\, |B(\varrho n)|\; M_3 \, M_4 \;
    \partial_t \Big( \big| \la + \xi - \overline{W}(t) \big|^{-1} \Big).
  \end{align*}
  Integrating from $t_0$ to $t_0 + n^2$ gives
  \begin{align*}
    m
    \Big[
      \big\{
        (t,z) \in I \times B( \varrho n) :\,
        \overline{w}(t,z) \,>\, \xi + \la
      \big\}
    \Big]
    \;\leq\;
    c_2\, M_3 \, M_4 \, \frac{n^2\, |B(\varrho n)|}{\la}.
  \end{align*}
  On the other hand, an application of the Markov inequality yields
  \begin{align*}
    m
    \Big[
      \big\{
        (t,z) \in I \times B(\varrho n) :\,
        \frac{c_1}{n^{2}}\, M_3 \, (t - s')
        \,>\,
        \tfrac{\la}{2}
      \big\}
    \Big]
    \;\leq\;
    c_1\, M_3 \, \frac{n^2\, |B(\varrho n)|}{\la}.
  \end{align*}
  Since $-\ln u = w = \overline{w} + c_1 \, n^{-2}\, M_3 \, (t - t_0)$, we obtain finally
  \begin{align*}
    &m
    \Big[
      \big\{ (t,z) \in I \times B(\varrho n) :\,
      \ln u(t, z)  \,<\, -\xi - \la \big\}
    \Big]
    \\[.5ex]
    &\mspace{36mu}\leq\;
    m
    \Big[
      \big\{
        (t,z) \in I \times B(\varrho n) :\,
        \overline{w}(t, z) 
        \,>\,
        \xi + \tfrac{\la}{2}
      \big\}
    \Big]
    \\
    &\mspace{72mu}+\;
    m
    \Big[
      \big\{
        (t,z) \in I \times B(\varrho n) :\,
        \frac{c_1}{n^{2}}\, M_3 \, (t - s')
        \,>\,
        \tfrac{\la}{2}
      \big\}
    \Big]
    \\[.5ex]
    &\mspace{36mu}\leq\;
    c\,
    \big(
      M_3\, M_4+ M_3
    \big)\,
    n^2\, |B(\varrho n)|\, \la^{-1}, 
  \end{align*}
  which proves the first inequality. Working instead of $D_t^+(\la)$ with the set $D_t^-(\la) = \big\{ z \in B(\varrho n) :\, \overline{w}(t,z) < \xi - \la \big\}$, the second inequality follows by a similar argument.
\end{proof}
\begin{proof} [Proof of Theorem~\ref{thm:PHI}]
  Consider a caloric function $u > 0$ on space-time cylinder $Q(n) = [t_0, t_0 + n^2] \times B(n)$ with $B(n) = B(x_0, n)$ and fix some $\varrho \in [3/4, 1)$.  Consider the function $f_1 \ldef  u^{-1}\, \me^{-\xi}$, where $\xi = \xi(u)$ be as in Lemma~\ref{lem:PHI:log}.  Our aim is to apply Lemma~\ref{lem:bombieri_giusti} to the function $f_1 \ldef u^{-1}\, \me^{-\xi}$ with $U = [t_0, t_0 + n^2] \times B(\varrho n)$,
  \begin{align*}
    U_{\si}
    \;=\;
    [\si\, t_0 + (1 - \si)\,s'\!,\, t_0 + n^2] \times B(\si \varrho n)
    \qquad \text{with} \quad
    s' = t_0 + \frac{3 \varrho}{4 \varrho - 2}\; n^2,
  \end{align*}
  $\de = 1/(2 \varrho)$, $\al^* = \infty$ and $m = \md t \times | \cdot |$.  Then, the condition (i) of Lemma~\ref{lem:bombieri_giusti} follows from \eqref{eq:PHI:max_est1} and H\"older's inequality with $\ga = 2 \ka'$ and constant $C_{\mathrm{BG1}}$ given by
  \begin{align*}
    C_{\mathrm{BG1}}
    \;=\;
    c\,
    \bigg(
      \Big(1 \vee \Norm{\mu^\om}{p,B(n)}\Big)\;
      \Big(1 \vee \Norm{\nu^\om}{q,B(n)}\Big)
    \bigg)^{\!\!1 + \ka'}.
  \end{align*}
  Since $\ln f_1 > \la$ is equivalent to $\ln u < -\la -\xi$, the condition (ii) of Lemma~\ref{lem:bombieri_giusti} follows immediately from Lemma~\ref{lem:PHI:log}.  Noting that $U_{\de} = Q_+$, an application of Lemma~\ref{lem:bombieri_giusti} gives
  \begin{align} \label{eq:phi_bg1}
    \max_{(t, x) \in Q_+} f_1(t,x) \;\leq\; A_1^\om 
    \qquad \Longleftrightarrow \qquad
    \me^{-\xi} \;\leq\;  A_1^\om \min_{(t,x) \in Q_+} u(t,x).
  \end{align}

  Next, consider the function $f_2 = u\, \me^{\xi}$, where $\xi = \xi(u)$ be again as in Lemma~\ref{lem:PHI:log}.  Now, set $U_{\si} = \big[\si\, t_0 + (1 - \si)\, s', (1 - \si)\, s'' + \si (t_0 + n^2)\big] \times B(\si \varrho n)$ with
  \begin{align*}
    s'
    \;=\;
    t_0 + \frac{\varrho}{4 \varrho - 2}\, n^2
    \qquad \text{and} \qquad
    s''
    \;=\;
    t_0 + \frac{\varrho - 1}{2 \varrho - 1}\, n^2,
  \end{align*}
  where $\de = 1/(2 \vr)$, $\al^* \in (0,1/2)$ arbitrary but fixed and $m$ are as before.  Similarly to the above procedure, the two conditions of Lemma~\ref{lem:bombieri_giusti} are verified due to Proposition~\ref{prop:PHIsuper:mos_it} and Lemma~\ref{lem:PHI:log}. Hence, $\Norm{f_2}{\al^*,U_{\de}}\leq A_2^{\om}$, and  since $U_{1/2 \vr} = Q_-$ we obtain from Corollary~\ref{cor:PHI:max_est} and H\"older's inequality with $\ga = 2(\ka + \ka')$,
  \begin{align} \label{eq:phi_bg2}
    \max_{(t,x) \in Q_{-}} u(t,x)
    \;\leq\;
    c
    \Big(
      C_9(\al^*)
      \Big( 1 \vee \Norm{\mu^{\om}}{p,B(n)}\Big)
      \Big( 1 \vee \Norm{\nu^{\om}}{q,B(n)} \Big)
   \Big)^{\!1+\ka'+\ka}
    A_2^{\om}\; \me^{-\xi}.
  \end{align}
  By combining \eqref{eq:phi_bg1} and \eqref{eq:phi_bg2}, the assertion is immediate. Note that the quantity $M^\om_{x_0,n}$ enters the Harnack constant $C_{\mathrm{PH}}$ through the application of Lemma~\ref{lem:PHI:log}.
\end{proof}

\subsection{Near-diagonal estimates and H\"older-continuity}
It is a well known fact that the parabolic Harnack inequality is a powerful property, which has a number of important consequences. We state here only a few of them, namely it provides immediately near-diagonal estimates on the heat kernel and a quantitative H\"older continuity estimate for caloric functions. 
\begin{prop}[Near-diagonal heat-kernel bounds]\label{prop:hke}
  Let $t > 0$ and $x_1 \in V$.  Then,
  \begin{enumerate}[(i)]
    \item for any $x_2 \in V$,
      \begin{align*}
        q^\om(t, x_1, x_2)
        \;\leq\; 
        \frac{C_{\mathrm{PH}}}{\mu^{\om}\big[B(x_1, \tfrac{1}{2} \sqrt{t})\big]}
        \;\leq\;
        \frac{C_{\mathrm{PH}}\, C_{\mathrm{reg}}}
             {\Norm{\mu^{\om}}{1, B(x_1, \sqrt{t} / 2)}}\;
        2^{d}\,t^{-\frac{d}{2}},
      \end{align*} 
    \item for any $x_2 \in B\big(x_1,\tfrac{1}{2} \sqrt{t}\big)$
      \begin{align*}
        q^\om(t,x_1,x_2)
        \;\geq\;
        \frac{1}{C_{\mathrm{PH}}^{2}\, \big|B(x_1, \frac{1}{2} \sqrt{t}) \big|}
        \;\geq\;
        \frac{C_{\mathrm{reg}}}{C_{\mathrm{PH}}^{2}}\,
        \min\Big\{1, 2^{d}\, t^{-\frac{d}{2}} \Big\}
      \end{align*}
  \end{enumerate}
  where $C_{\mathrm{PH}}\big(M^\om_{x_1,\sqrt{t}}, \|\mu\|_{p, B(x_1, \sqrt{t})}, \|\nu\|_{q, B(x_1, \sqrt{t})}\big)$ is the constant in Theorem \ref{thm:PHI}.
\end{prop}
\begin{proof}
  We will proceed as in \cite[Proposition 3.1]{De99}.  Given $x_1, x_2 \in V$ and $t > 0$ we apply the PHI to the caloric function $q^\om(\cdot, \cdot, x_2)$ on $[t_0, t_0 + t] \times B(x_1, \sqrt t)$ with $t_0 = \frac{3}{4} t$.  This yields for any $z \in B(x_1,\frac{1}{2} \sqrt t)$
  \begin{align*}
    q^{\om}(t, x_1, x_2)
    \;\leq\;
    C_{\mathrm{PH}}\, q^{\om}\big(\tfrac{3}{2} t, z, x_2\big),
  \end{align*}
  with $C_{\mathrm{PH}}= C_{\mathrm{PH}}\big(\|\mu\|_{p,B(x_1, \sqrt{t})}, \|\nu\|_{q, B(x_1, \sqrt{t})} \big)$.  Hence, by multiplying both sides with $\mu^{\om}(z)$, a summation over all $z \in B(x_1, \tfrac{1}{2} \sqrt{t})$ gives
  \begin{align*}
    q^{\om}(t, x_1, x_2)
    \;\leq\;
    \frac{C_{\mathrm{PH}}}{\mu^{\om}\big[B(x_1, \tfrac{1}{2} \sqrt{t})\big]}\,
    \Prob_{x_2}^{\om}\!\Big[Y_{3t/2} \in B(x_1, \tfrac{1}{2} \sqrt{t}) \Big]
    \;\leq\;
    \frac{C_{\mathrm{PH}}}{\mu^{\om}\big[B(x_1, \tfrac{1}{2} \sqrt{t})]},
  \end{align*}
  where we used in the first step both the symmetry and the definition of the transition density.  The assertion (i) follows now from the volume regularity assumption. 
 
  In order to prove (ii), suppose that $x_2 \in B(x_1, \tfrac{1}{2} \sqrt{t})$.  Similar to the above considerations an application of Theorem~\ref{thm:PHI} to the caloric function $q^{\om}(\cdot, \cdot, x_2)$ on $[0, t] \times B(x_1, \sqrt{t})$ yields that $q^{\om}(\tfrac{1}{2} t, z, x_2) \leq C_{\mathrm{PH}}\, q^{\om}(t, x_1, x_2)$ for any $z \in B(x_1, \tfrac{1}{2} \sqrt{t})$ and $C_{\mathrm{PH}}$ as defined above.  In particular,
  \begin{align*}
    \sum_{z \in B(x_1, \sqrt{t} / 2)}\mspace{-15mu}
    q^{\om}\big(\tfrac{1}{2} t, z, x_2\big)
    \;\leq\;
    C_{\mathrm{PH}}\, \big|B(x_1, \tfrac{1}{2} \sqrt{t}) \big|\, q^{\om}(t, x_1, x_2).
  \end{align*}
  Next, consider the function
  \begin{align*}
    u(\tau, x)
    \;\ldef\;
    \left\{
      \begin{array}{ll}
        1, & \text{if } \tau \in [0, t / 2],
        \\[1ex]
        \sum_{z \in B(x_1, \sqrt{t} / 2)}\,
        q^{\om}(\tau - t / 2, x, z),
        & \text{if } \tau \in [t / 2, t].
      \end{array}
    \right.
  \end{align*}
  Notice that the function $u$ is caloric on $[0, t] \times B(x_1, \sqrt{t})$ and by the parabolic Harnack inequality we have that
  \begin{align*}
    1
    \;=\;
    u\big(\tfrac{1}{2} t, x_2\big)
    \;\leq\;
    C_{\mathrm{PH}}\mspace{-15mu} \sum_{z \in B(x_1, \sqrt{t} / 2)} \mspace{-15mu}
    q^{\om}\big(\tfrac{1}{2}t, z, x_2 \big)
    \;\leq\;
    C_{\mathrm{PH}}^2\, \big|B(x_1, \tfrac{1}{2} \sqrt{t}) \big|\,
    q^{\om}(t, x_1, x_2),
  \end{align*}
  where we used the definition of $u$ and the symmetry of the transition density in the second step.  Now, the claim follows from the volume regularity assumption because $\big|B(x_1, \tfrac{1}{2} \sqrt{t}) \big| \leq C_{\mathrm{reg}} \max\big\{1, 2^{-d}\, t^{d/2}\big\}$.
\end{proof}
\begin{prop} \label{prop:cont_caloric}
  Suppose that Assumption \ref{ass:lim_const} holds. Let $x_0 \in V$ and let $s(x_0)$ and $C_{\mathrm{PH}}^*$ be as in Assumption \ref{ass:lim_const}.  Further, let $R\geq s(x_0)$ and $\sqrt T \geq R$. Set $T_0 \ldef T+1$ and $R_0^2 \ldef T_0$ and suppose that $u > 0$ is a caloric function on $[0, T_0] \times B(x_0, R_0)$.  Then, for any $x_1, x_2 \in B(x_0,R)$ we have
  \begin{align*}
    \big|u(T,x_1) - u(T,x_2)\big|
    \;\leq\;
    c\, \bigg( \frac{R}{T^{1/2}} \bigg)^{\!\!\th}\,
    \max_{[3 T_0 / 4, T_0] \times B(x_0, R_0 / 2)} u,
  \end{align*}
  where $\th = \ln (2C^*_{\mathrm{PH}}/(2C^*_{\mathrm{PH}}-1)) / \ln 2$ and $c$ only depending on $C^*_{\mathrm{PH}}$.
\end{prop}
\begin{proof}
  We will follow the arguments in \cite[Lemma 6]{CH08} (cf.\ also \cite[Proposition 3.2]{BH09}), but some extra care is needed due to the special form of the Harnack constant.  Set $R_k = 2^{-k} R_0$. Let
  \begin{align*}
    Q_k
    \;=\;
    [T_0 - R_k^2, T_0] \times B(x_0,R_k), 
  \end{align*}
  and let $Q_k^-$ and $Q_k^+$ be accordingly defined as
  \begin{align*}
    Q_k^-
    \;=\;
    \big[T_0 - \tfrac{3}{4} R_k^2, T_0 - \tfrac{1}{2} R_k^2 \big]
    \times
    B(x_0, \tfrac{1}{2} R_k),
    \mspace{15mu}
    Q_k^+
    \;=\;
    \big[T_0 - \tfrac{1}{4} R_k^2, T_0 \big] \times B(x_0,\tfrac{1}{2} R_k).
  \end{align*}
  In particular, note that $Q_{k+1}\subset Q_k$ and $Q_k^+ = Q_{k+1}$. Setting 
  \begin{align*}
    v_k
    \;=\;
    \frac{u - \min_{Q_k} u}{\max_{Q_k} u - \min_{Q_k}u},
  \end{align*}
  we have that $v_k$ is caloric, $0 \leq v_k \leq 1$ and $\osc(v_k, Q_k) = 1$. (Here $\osc(u,A)=\max_A u - \min_A u$ denotes the oscillation of $u$ on $A$.) Replacing $v_k$ by $1-v_k$ if necessary we may assume that $\max_{Q_k^-} v_k \geq \frac{1}{2}$. 

  Now for every $k$ such that $R_{k} \geq s(x_0)$ we apply the PHI in Theorem \ref{thm:PHI} and using the definition of $C^*_{\mathrm{PH}}$ we get
  \begin{align*}
    \tfrac{1}{2}
    \;\leq\;
    \max_{Q_k^-} v_k
    \;\leq\; 
    C_{\mathrm{PH}}^*\, \min_{Q_k^+} v_k. 
  \end{align*}
  Since $\osc(u,Q_{k+1}) =  \osc(u,Q_k^+)$, it follows that for such $k$,
  \begin{align*}
    \osc(u,Q_{k+1})
    &\;=\;
    \frac{\max_{Q_k^+} u - \min_{Q_k^+} u}{\osc(u, Q_k)} \osc(u,Q_k)
    \\[1ex]
    &\;=\;
    \bigg(
      1 \,+\, \frac{\max_{Q_k^+} u - \max_{Q_k} u }{\osc(u, Q_k)} \,-\,
      \min_{Q_k^+} v_k
    \bigg)\; \osc(u, Q_k).
  \end{align*}
  Hence,
  \begin{align*}
    \osc\big(u, Q_{k+1}\big) \;\leq\; (1 - \de)\, \osc(u, Q_k)
  \end{align*}
  with $\de = (2 C_{\mathrm{PH}}^*)^{-1}$.  Now we choose $k_0$ such that $R_{k_0} \geq R > R_{k_0+1}$. Then, we can iterate the above inequality on the chain of boxes $Q_1 \supset Q_2 \supset \cdots \supset Q_{k_0}$ to obtain
  \begin{align*}
    \osc\big(u, Q_{k_0}\big)
    \;\leq\;
    (1 - \de)^{k_0-1}\, \osc(u,Q_1).
  \end{align*}
  Note that $B(x_0, R) \subseteq B(x_0, R_{k_0})$, $T \in [T_0 - R_{k_0}^2,T_0]$ and $Q_1 = Q_0^+$.  Finally, since $(1 - \de)^{k_0} \leq c (R / T^{1/2})^\th$, the claim follows.
\end{proof}

\section{Local limit theorem} \label{sec:lclt}
In this section, we prove the quenched local limit theorem for the random conductance model as formulated in Theorem~\ref{thm:llt_csrw}.  For this reason, we consider the $d$-dimensional Euclidean lattice $(\bbZ^d, E_d)$ as the underlying graph and assume throughout this section that the random weights $\om$ satisfy the Assumptions~\ref{ass:environment} and \ref{ass:moment}.  Let us recall that $k_t = k_t^{\Si_Y}\!$, defined in \eqref{eq:def:GHK}, denotes the Gaussian heat kernel with diffusion matrix $\Si_Y^2$.

The proof of the local limit theorem is based on the approach in \cite{BH09} and \cite{CH08}.  In the sequel, we will briefly explain the strategy.  For $x \in \bbR^d$ and $\de >0$, let
\begin{align}\label{eq:def:J}
  J(t, n)
  \;\equiv\;
  J(t, x, \de, n)
  \;\ldef\;
  \Prob_0^{\om}\!\Big[ Y^{(n)}_t \in C(x, \de) \Big]
  \,-\, \int_{C(x, \de)} k_t(y)\, \md y.
\end{align}
Here, we denote by $C(x, \de) = x + [-\de, \de]^d$ the cube with center $x$ and side length $2 \de$, i.e.\ $C(x, \de)$ is a ball in $\bbR^d$ with respect to the supremum norm $|\cdot|_{\infty}$.  Further, for any $n \geq 1$, we set $C^n(x, \de) \ldef n C(x, \de) \cap \bbZ^d$.

Now, let us rewrite \eqref{eq:def:J} as $J(t,n) = J_1(t,n) + J_2(t, n) + J_3(t, n) + J_4(t, n)$, where
\begin{align*}
  J_1(t, n)
  &\;\ldef\;
  \sum_{z \in C^n(x,\delta)}
  \Big(
    q^{\om}\big(n^2 t, 0, z\big)
    \,-\, q^{\om}\big(n^2 t, 0, \lfloor nx \rfloor \big)
    \Big)\,
  \mu^{\om}(z),
  \\
  J_2(t, n)
  &\;\ldef\;
  \mu^{\om}\big[C^n(x,\delta) \big]\,
  \Big(
    q^{\om}\big(n^2 t, 0, \lfloor nx \rfloor \big)
    \,-\, n^{-d}\, a\, k_t(x)
  \Big), 
  \\[1.5ex]
  J_3(t, n)
  &\;\ldef\;
  k_t(x)\, \Big( \mu^{\om}\big[C^n(x,\delta) \big]\, n^{-d}\, a - (2\de)^d \Big),
  \\[.75ex]
  J_4(t, n)
  &\;\ldef\;
  \int_{C(x,\delta)} \Big( k_t(x) - k_t(y) \Big)\, \md y,
\end{align*}
and $a = 1 / \mean\big[\mu^{\om}(0)\big]$.  Thus, in order to conclude a quenched local limit theorem, it remains to establish suitable upper bounds on the terms $J_1$, $J_3$ and $J_4$ uniformly for $t \in [T_1, T_2] \subset (0, \infty)$ that are small compared to $\mu^{\om}[C^n(x, \de)] n^{-d}$.  Further, note that by means of the spatial ergodic theorem, cf.\ Lemma~\ref{lemma:ergodic} below, $\mu^{\om}[C^n(x, \de)] n^{-d}$ converges $\prob$-a.s to $\mean[\mu^{\om}(0)] (2 \de)^{d}$ as $n$ tends to infinity.  Hence, we obtain that $\lim_{n \to \infty} |J_3(t, n)| / \mu^{\om}[C^n(x, \de)] n^{-d} = 0$ and $\lim_{n \to \infty} |J(t,n)| / \mu^{\om}[C^n(x, \de)] n^{-d} = 0$ uniformly for all $t \in [T_1, T_2]$.  The latter follows from the quenched invariance principle, see Theorem~\ref{thm:ip}.  Further, note that $q^{\om}(\cdot, 0, \cdot)$ is a caloric function.  Thus, provided that the assumptions of Proposition~\ref{prop:cont_caloric} are satisfied, the H\"{o}lder continuity combined with the near-diagonal 
upper bounds of Proposition~\ref{prop:hke} implies 
\begin{align*}
  \max_{z \in C^{n}(x, \de)}\, n^{d}\,
  \big|
    q^{\om}\big( n^2 t, 0, z \big)
    \,-\,
    q^{\om}\big( n^2 t, 0, \lfloor n x \rfloor \big)
  \big|
  \;\leq\;
  c\, \de^{\th}\, T_1^{-\frac{d}{2} - \frac{\th}{2}}.
\end{align*}
This estimate allows us to deduce $\lim_{n \to \infty } |J_1(t,n)| / \mu^{\om}[C^n(x, \de)] n^{-d}\big| = \cO(\de^{\th})$ uniformly for $t \in [T_1, T_2]$, which vanishes when $\de$ tends to $0$.  The verification of the assumptions of Proposition~\ref{prop:cont_caloric} will be based on the ergodic theorem as well.  Finally, an uniform upper bound of $|J_4(t, n)| / \mu^{\om}[C^n(x, \de)] n^{-d}$ can be deduced from the explicit form of $k_t(x)$.

For the proof Theorem~\ref{thm:llt_csrw} we will apply the general criterion in \cite[Theorem~1]{CH08}, so we will just verify the required assumptions on the underlying sequence of rescaled graphs, the weights and a PHI on sufficiently large balls.

Let us start with some immediate consequences of the ergodic theorem.
\begin{lemma}\label{lemma:ergodic}
  For any $x \in \bbR^d$ and $\de > 0$ we have that, for $\prob$-a.e.\ $\om$,
  \begin{align*}
    \lim_{n \to \infty}\, \Norm{\mu^{\om}}{p, C^n(x, \de)}^p
    \;=\;
    \mean\!\big[\mu^{\om}(0)^p\big]
    \quad\; \text{and} \quad\;
    \lim_{n \to \infty}\, \Norm{\nu^{\om}}{q,C^n(x, \de)}^q
    \;=\;
    \mean\!\big[\nu^{\om}(0)^q\big].
  \end{align*}
  In particular, for every $\ve > 0$ we have $r^{\om}(x, \de, \ve) < \infty$ for $\prob$-a.e.\ $\om$, where
  \begin{align*}
    r^{\om}(x, \de, \ve)
    \;\ldef\;
    \inf
    \Big\{
      n \geq 1 \;:\;\; 
      &\Big|
        \Norm{\mu^{\om}}
          {p, C^{n}(x, \de)}^p \,-\, \mean\!\big[\mu^{\om}(0)^p\big]
      \Big|
      \;<\;
      \ve
      \\[1ex]
      &\mspace{12mu}\text{and}\mspace{12mu}
      \Big|
        \Norm{\nu^{\om}}
          {q, C^{n}(x, \de)}^q \,-\, \mean\!\big[\nu^{\om}(0)^q\big]
      \Big|
      \;<\;
      \ve
    \Big\}.
  \end{align*}
\end{lemma}
\begin{proof}
  The assertions of the lemma are an immediate consequence of the spatial ergodic theorem \cite[Theorem~2.8 in Chapter 6]{Kr85}, i.e.\
  \begin{align*}
    \lim_{n \to \infty} \Norm{\mu^{\om}}{p, C^n(x, \de)}^p
    \;=\;
    \lim_{n \to \infty} \frac{1}{\big| C^n(x, \de) \big|}\,
    \sum_{x \in C^n(x, \de)}\mspace{-12mu}
    \big( \mu^{\tau_x \om}(0) \big)^p
    \;=\;
    \mean\!\big[\mu^{\om}(0)^p\big],
  \end{align*}
  provided we have shown that the sequence $\big\{ C^n(x, \de) : n \in \bbN \big\}$ is regular.  For this purpose, consider $\big\{C^{n}\big(0, |x|_{\infty} + \de \big) : n \in \bbN \big\}$.  Obviously, this sequence of cubes is increasing and have the property that $C^n(x, \de) \subset C^{n}\big(0, |x|_{\infty} + \de \big)$ for all $n \in \bbN$.  Moreover, there exists $K \equiv K(x, \de) < \infty$ such that 
  \begin{align*}
    \big| C^{n}\big(0, |x|_{\infty} + \de \big) \big|
    \;\leq\;
    K\,
    \big| C^n(x, \de) \big|,
    \qquad \forall\, n \in \bbN.
  \end{align*}
  Thus, the sequence $\big\{ C^n(x, \de) : n \in \bbN \big\}$ is regular. 
\end{proof}
\begin{remark}\label{rem:ergodic:balls}
  Note that $\big\{ B(\lfloor x n \rfloor, \de n) : n \in \bbN \big\}$ is also a regular sequence of boxes for every given $x \in \bbR^d$ and $\de > 0$.  Thus, the assertion of Lemma~\ref{lemma:ergodic} also holds if we replace $C^n(x,\de)$ by $B(\lfloor n x \rfloor, \de n)$.
\end{remark}
\begin{proof}[Proof of Theorem~\ref{thm:llt_csrw}]
  We just need to verify Assumption~1 and 4 in \cite{CH08}.  Then, the assertion follows immediately from \cite[Theorem~1]{CH08}.

  Consider the sequence of rescaled Euclidean lattices $\big\{ (\frac{1}{n} \bbZ^d, \frac{1}{n} E_d) : n \in \bbN \big\}$ on the underlying metric space $(\bbR^d, |\cdot - \cdot|_{\infty})$ and set $\al(n) = n$, $\be(n) = n^d / \mean[\mu^{\om}(0)]$ and $\ga(n) = n^2$.  Recall that we denote by $d(x,y)$ the graph distance for $x, y \in \tfrac{1}{n} \bbZ^d$.  Since,
  \begin{align*}
    \al(n)\, |x - y|_{\infty}
    \;\leq\;
    d(x,y)
    \;\leq\;
    d\, \al(n)\, |x - y|_{\infty},
    \qquad \forall\, x, y \in \tfrac{1}{n}\, \bbZ^d,
  \end{align*}
  and for all $\de > 0$
  \begin{align*}
    \lim_{n \to \infty}\, \tfrac{1}{n}\, \bbZ^d \cap C(0, \de) \;=\; C(0, \de)
  \end{align*}
  with respect to the usual Hausdorff topology on non-empty compact subsets of $(\bbR^d, |\cdot - \cdot|_{\infty})$, the Assumption~1(a), (b) and Assumption~4(b) are satisfied.  The Assumption~1(c) requires that for every $x \in \bbR^d$ and $\de > 0$,
  \begin{align*}
    \lim_{n \to \infty}\, \be(n)^{-1}\, \mu^{\om}\big[C^n(x,\de)\big]
    \;=\;
    \big(2 \de\big)^d
    \;=\;
    \int_{C(x,\de)} \md y.
    \qquad \prob\text{-a.s.}
  \end{align*}
  But, in view of Lemma~\ref{lemma:ergodic} and the definition of $\be(n)$, this assumption is satisfied.  Further, the Assumption~1(d) requires that for any compact interval $I \subset (0, \infty)$, $x \in \bbR^d$ and $\de > 0$
  \begin{align*}
    \lim_{n \to \infty}
    \Prob_0^{\om}\!\Big[\tfrac{1}{n} Y_{\ga(n) t} \in C(x, \de) \Big]
    \;=\;
    \int_{C(x, \de)} k_t(y)\, \md y,
    \qquad \prob\text{-a.s.}
  \end{align*}
  uniformly for all $t \in I$.  As we have discussed above, this is a direct consequence of Theorem~\ref{thm:ip}.

  It remains to verify the Assumption~4(c): For every $x_0 \in \frac{1}{n} \bbZ^d$, $n \geq 1$, there exists $s_n(x) < \infty$ such that the parabolic Harnack inequality with constant $C_{\mathrm{PH}}^*$ holds on $Q(r) \equiv [0, r^2] \times B(x_0, r)$ for all $r \geq s_n(x_0)$ and for every $x \in \bbR^d$ we have $\al(n)^{-1} s_n(\lfloor n x \rfloor) \to 0$ as $n \to \infty$.
  
  By the ergodic theorem we know that $\|\mu^{\om}\|_{p, B(x_0, n)}$ and $\|\nu^{\om}\|_{q, B(x_0, n)}$ converge as $n$ tends to infinity for $\prob$-a.e.\ $\om$.  In particular, by Lemma~\ref{lemma:ergodic} and Remark~\ref{rem:ergodic:balls}, we have  for all $x \in \bbR^d$ and $\ve > 0$ that $\prob$-a.s.\ there exists $r^{\om}(x, 1, \ve) < \infty$ such that
  \begin{align*}
    \Norm{\mu^{\om}}{p, B(\lfloor n x \rfloor, n)}
    \;\leq\;
    \big( \ve \,+\, \mean\!\big[\mu^{\om}(0)^p\big] \big)^{\frac{1}{p}},
    \qquad 
    \Norm{\nu^{\om}}{q, B(\lfloor n x \rfloor, n)}
    \;\leq\;
    \big( \ve \,+\, \mean\!\big[\nu^{\om}(0)^q\big] \big)^{\frac{1}{q}}
  \end{align*}
  for all $n \geq r^{\om}(x, 1, \ve)$.  Since the Harnack constant appearing in Theorem~\ref{thm:PHI} depends monotonically on $\|\mu^{\om}\|_{p, B(\lfloor n x \rfloor, n)}$ and $\|\nu^{\om}\|_{q, B(\lfloor n x \rfloor, n)}$, this implies that for any $r \geq r^{\om}(x, 1, \ve)$ the parabolic Harnack inequality holds on $Q(r)$ with constant $C_{\mathrm{PH}}^* \ldef C_{\mathrm{PH}}\big((\ve + \mean[\mu^{\om}(0)^p])^{1/p}, (\ve + \mean[\nu^{\om}(0)^q])^{1/q} \big)$.  Thus, the Assumption~4(c) of \cite{CH08} is satisfied as well.
\end{proof}
We now give an example, communicated to us by Noam Berger, which shows that the moment condition in Assumption~\ref{ass:moment} is optimal for a local limit theorem. For preparation we first note in the following lemma that a local limit theorem implies certain near-diagonal bounds on the heat kernel.
\begin{lemma} \label{lem:near-diag}
  Suppose that the local limit theorem holds in the sense that
  \begin{align}  \label{eq:lclt_restat}
    \lim_{n \,\to\, \infty}  
    \sup_{|x| \leq 1}\, \sup_{t \in [1,4]}\,
    \Big|
      n^{d}\, q^{\om}\big(n^2 t, 0, \lfloor nx \rfloor\big) - a\, k_t(x)
    \Big|
    \;=\;
    0,
    \qquad  \prob \text{-a.s.} 
  \end{align}
  with $a \ldef 1 / \mean[\mu^{\om}(0)]$.  Then, there exists $N_0 = N_0(\om)$ and constants $C_{12}$ and $C_{13}$ such that $\prob$-a.s.\ for all $t \geq N_0$ and $x \in B(0,\sqrt{t})$,
  \begin{align}\label{eq:near-diag-bound}
    C_{12} t^{-d/2}
    \;\leq\;
    q^\om(t,0,x)
    \;\leq\;
    C_{13} t^{-d/2}.
  \end{align}
\end{lemma}
\begin{proof}
  Setting $n=\lfloor \sqrt{t} \rfloor$ and $h(t)=t/\lfloor \sqrt t \rfloor^2$ we write
  \begin{align*}
    q^\om(t,0,x)
    \;=\;
    n^{-d}\,
    \left(
      n^d q^\om\big(n^2 h(t), 0, \lfloor n (x/n) \rfloor \big) - a\, k_{h(t)}(x/n)
    \right)
    + n^{-d}\, a k_{h(t)}(x/n).
  \end{align*}
  Since $|x|/n \leq 1$ and $h(t) \in [1,4]$ for $t \geq 1$ we have $c^{-1} \leq k_{h(t)}(x/n)\leq c$ for some $c > 0$ independent of $t$, $x$ and $n$.  Further, by \eqref{eq:lclt_restat} for $\prob$-a.e.\ $\om$,
  \begin{align*}
    &\Big|
      n^d\, q^\om(n^2 h(t),0, \lfloor n (x/n) \rfloor ) \,-\, a\, k_{h(t)}(x/n)
    \Big|
    \\[.5ex]
    &\mspace{36mu}\leq\;
    \sup_{|y| \leq 1} \sup_{s\in [1,4]}
    \Big| n^d\, q^\om(n^2s, 0, \lfloor n y\rfloor) \,-\, a\, k_s(y) \Big|
    \;\rightarrow\;
    0   
  \end{align*}
  as $n \to \infty$, and we obtain \eqref{eq:near-diag-bound}.
\end{proof}
\begin{theorem} \label{thm:ex_lclt}
  Let $d \geq 2$ and  $p,q \in (1,\infty]$ be such that $1/p + 1/q > 2/d$.  There exists an ergodic law $\prob$ on $(\Om, \cF)$ with $\mean[\om(e)^p] < \infty$ and $\mean[\om(e)^{-q}] < \infty$  such that the conclusion of Lemma~\ref{lem:near-diag} does not hold.
\end{theorem}
The strategy of the proof is to construct an ergodic environment such that in every large box one can find a sufficiently deep trap for the walk. 
\begin{proof}
  \emph{Step 1: Definition of the environment}.  We first construct an ergodic law on $(\Om, \cF)$.  Let $\{\xi(e): e\in E_d\}$ be a collection of i.i.d.\ random variables uniformly distributed over $[0,1]$. Further let $k_0\in \bbN$ be such that $2^{-k_0}<p_c$, where $p_c$ denotes the critical probability for bond percolation on $\bbZ^d$. Then, for $k\geq k_0$ we define $\{\th^k(e) : e \in E_d\}$ by
  \begin{align*}
    \th^k(e) 
    \;=\;
    \begin{cases} 
      1  & \text{if $\xi(e)\leq 2^{-k}$,} \\
      0 & \text{otherwise.}
    \end{cases}
  \end{align*}
  Further, set $\cO^k \ldef \{ e \in E_d: \th^k(e)=1 \}$ and let $\{\cH^k_j \}_{j \geq 1}$ be the collection of connected components of the graph $(\bbZ^d, \cO^k)$ and $\cH^k \ldef \bigcup_{j\geq 1} \cH^k_j$, $k\geq k_0$. Note that $\prob$-a.s.\  every component  $\cH^k_j$ is finite by our choice of $k_0$. Notice that by construction the components are decreasing in $k$, i.e.\ for each $\cH^{k+1}_j$ there exists a unique $i$ such that $\cH^{k+1}_j \subseteq \cH^k_i$.  Since we have assumed $1/p + 1/q > 2/d$ there exist $\al < 1/p$, $\be < 1/q$ and $\varepsilon >0$ such that $(\al + \be-\varepsilon)d /2 > 1$.
  \begin{figure}[t]
    \begin{center}
      \includegraphics{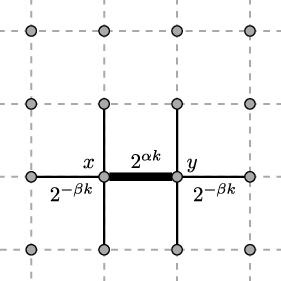}
    \end{center}
    \caption{Illustration a $k$-trap in $d = 2$.}
    \label{fig:traps}
  \end{figure}
  We now define the conductances by an iteration procedure over $k$ starting from $k_0$.  We initialize $\om(e) = 1$ for all $e \in E_d$. In the iteration step $k$ for all $j \geq 1$ we  choose in each component $\cH^k_j$ one edge $e_{k,j}$ uniformly at random and set $\om(e_{k,j}) = 2^{\al k}$ and $\om(e') = 2^{-\be k}$ for all $e' \sim e_{k,j}$. Since the components $\{\cH^k_j\}$ are separated by at least one edge, this procedure is well-defined and does not allow high conductances on neighbouring edges.  That is, on each component $\cH^k_j$ we put one trap (see Figure~\ref{fig:traps}).  Note that during the iteration it might happen that some traps are replaced by even deeper traps and that some of the weaker bonds are left when one of the associated bonds with high conductance has been overwritten. 
  
  Henceforth, we will call an edge $e \in E_d$ a \emph{$k$-trap} if $\om(e) \geq 2^{\al k}$ and  $\om(e') \leq 2^{-\be k}$ for all $e' \sim e$.  Note that if the CSRW starts at a $k$-trap $\{x,y\} \in E_d$ it will spend at least $O(2^{(\al+\be)k})$units of time  in $\{x,y\}$.  In particular, there exists $C_{14} > 0$ independent of $k$ such that for $s\leq 2^{(\al+\be)k}$,
  \begin{align}\label{eq:hke:estimate}
    \Prob^{\om}_x\big[Y_s = x\big] \;\geq\; C_{14}
    \qquad \text{for all } k \geq k_0.
  \end{align}
   
  \emph{Step 2: Ergodicity}.  The environment constructed in Step~1 is ergodic.  First, the stationarity is obvious from the construction.  To see that it is ergodic, consider two local functions $f, g \in L^2(\Om)$, meaning that $f$ and $g$ only depend on a finite set of coordinates $\{\om(e): e \in K \}$, for any $K \subset E_d$ finite. Set
  \begin{align*}
    D_K
    \;\ldef\;
    \max_{j: \cH^{k_0}_j \cap K \ne \emptyset} \diam \cH_j^{k_0} 
    \;<\;
    \infty
    \qquad \prob\text{-a.s.}
  \end{align*}
  Then, note that by construction the conductances on edges contained in $K$ are independent of the conductances on edges having distance at least $D_K + \diam K$. This implies
  \begin{align*}
    \mean\big[ f \cdot g\circ \tau_x \big]
    \;\rightarrow\;
    \mean[f] \, \mean[g]
    \qquad \text{as } |x| \to \infty.
  \end{align*}
  By a standard density argument  the same holds for all $f,g \in L^2(\Omega)$.  Thus, for any shift-invariant set $A\in \mathcal{F}$,
  \begin{align*}
    \prob\big[ A \big]
    \;=\;
    \prob\big[ A \cap \tau_x(A) \big]
    \;\rightarrow\;
    \prob[A]^2 
    \qquad \text{as } |x| \to \infty,
  \end{align*}
  and therefore $\prob[A] \in \{0,1\}$.
 
  \emph{Step 3: Moment condition}.  The environment satisfies the required moment condition.  Indeed, by construction we have for any edge $e$ the obvious bounds
  \begin{align*}
    \min_{e' \sim e} \xi(e')^{\be}
    \;\leq\;
    \min_{e' \sim e} \min_{k: \xi(e') \leq 2^{-k}} 2^{-\be k}
    \;\leq\;
    \om(e)\leq \max_{k: \xi(e) \leq 2^{-k}} 2^{\alpha k}
    \;\leq\;
    \xi(e)^{-\al},
  \end{align*}
  and therefore 
  \begin{align*}
    \mean \big[ \om(e)^p \big]
    &\;\leq\;
    \mean \big[\xi(e)^{-\al p} \big]
    \;=\;
    \int_0^1 u^{-\al p} \, \md u
    \;<\;
    \infty,
    \\
    \mean \big[\om(e)^{-q}\big]
    &\;\leq\;
    \mean \big[ \max_{e'\sim e} \xi(e')^{-\be q} \big]
    \;\leq\;
    2d\, \int_0^1 u^{- \be q}\, \md u
    \;<\;
    \infty,
  \end{align*} 
  where we used that $\al < 1/p$ and $\be < 1/q$.  
  
  \emph{Step 4:}  Let $A_n(e) \ldef \{\om(e)\geq 2^{n \al} \}$ and $R_n = (n 2^n)^{1/d}$.  In this step we show that almost surely for infinitely many $n$ we find an $n$-trap in $B(R_n) \equiv B(0,R_n)$. First,
  \begin{align*}
    &\prob \Big[\bigcap\nolimits_{e \in E(R_n)} A_n(e)^c\Big]
    \\
    &\mspace{36mu}\leq\;
    \prob\big[\cH^n \cap B(R_n) = \emptyset\big]
    \,+\,
    \prob\Big[
      \cH^n \cap B(R_n) \ne \emptyset,\, \bigcap\nolimits_{e\in E(R_n)} A_n(e)^c
    \Big],
  \end{align*}
  where $E(R_n)$ is the set of all edges with at least one endpoint contained in $B(R_n)$. The first term on the right hand side can be easily computed as
  \begin{align} \label{eq:nocluster}
    \prob\big[\cH^n \cap B(R_n)  = \emptyset \big]
    &\;=\;
    \prob\big[\th^n(e)=0, \, \forall\, e \in E(R_n)\big]
    \nonumber\\
    &\;=\;
    \prob\big[\xi(e)>2^{-n}\big]^{c R_n^d}
    \;=\;
    (1-2^{-n})^{cR_n^d}.
  \end{align}
  For the second term we obtain
  \begin{align*}
    &\prob\Big[
      \cH^n \cap B(R_n) \ne \emptyset,\; \bigcap\nolimits_{e \in E(R_n)} A_n(e)^c
    \Big]
    \\[.5ex]
    &\mspace{36mu}\leq\;
    \prob\big[\cH^n \cap B\big(\tfrac{1}{2} R_n\big) = \emptyset \big]
    \\
    &\mspace{72mu}+\, 
    \prob\Big[
      \cH^n \cap B\big(\tfrac{1}{2} R_n\big) \ne \emptyset,\,
      \cH^n \cap B(R_n) \not=\emptyset,\, \bigcap\nolimits_{e\in E(R_n)} A_n(e)^c
    \Big]
    \\[.5ex]
    &\mspace{36mu}\leq\;
    \prob\big[\cH^n \cap B\big(\tfrac{1}{2} R_n\big) = \emptyset\big]
    \,+\,
    \prob\big[
      \exists\, j \geq 1 : \mathrm{diam} \cH_j^n\geq \tfrac{1}{2} R_n
    \big]
    \\[.5ex]
    &\mspace{36mu}\leq\;
    \big(1 - 2^{-n}\big)^{cR_n^d} \,+\,  |\partial B(R_n)|\, c\, \me^{-c R_n},
  \end{align*}
  where in the last step we used a similar argument as in \eqref{eq:nocluster} and a tail estimate for the size of the holes in supercritical percolation clusters (see Theorems 6.10 and 6.14 in \cite{Gr99}).  Thus,
  \begin{align*}
    \prob\Big[\bigcap\nolimits_{e\in E(R_n)} A_n(e)^c\Big]
    \;\leq\;
    2 (1 - 2^{-n})^{cR_n^d} \,+\,  |\partial B(R_n)|\, c\, \me^{-c R_n},
  \end{align*}
  which by our choice of $R_n$ is summable over $n$.  In particular, by Borel-Cantelli there exists $n_0 = n_0(\om)$ such that $\prob$-a.s.\ for all $n \geq n_0$ there exists $\bar e_n \in E(R_n)$ such that the event $A_n(\bar e_n)$ occurs, that is for each $n \geq n_0$ we find an $n$-trap $\bar e_n$ in $B(R_n)$.
 
  \emph{Step 5:} Now suppose that the heat kernel bounds in Lemma~\ref{lem:near-diag} hold. For any   positive $\de < C_{14}$ there exists $n_{\de}$ such that, setting $t_1 = t_1(n) \ldef R_n^2$ and $t_2 = t_2(n) \ldef t_1^{(\al+\be-\varepsilon)d/2}$, we have for any $x \in B(R_n)$,
  \begin{align} \label{eq:hkratio}
    \frac{\Prob_0^\om[Y_{t_1+t_2} = x]}{\Prob_0^\om[Y_{t_1} = x]}
    \;=\;
    \frac{q^\om(t_1+t_2, 0, x)}{ q^\om(t_1, 0, x)}
    \;\leq\;
    c\, \Big(1 + \frac{t_2}{t_1}\Big)^{-\frac d 2}
    \;\leq\;
    \de,
    \qquad \forall\, n \geq n_{\de},
  \end{align}
  since $t_2/t_1 = R_n^{(\al+\be-\varepsilon)d-2}\gg 1$ as $(\al + \be-\varepsilon)d - 2 > 0$ by our choice of $\al$, $\be$ and $\varepsilon$.  On the other hand, note that $t_2\leq 2^{n(\al+\be)}$ for $n$ large enough. Thus, using the Markov property for $n \geq n_0$ and \eqref{eq:hke:estimate},
  \begin{align*}
    \frac{\Prob_0^\om[Y_{t_1+t_2} = \bar e_n^+]}{\Prob_0^\om[Y_{t_1} = \bar e_n^+]}
    \;\geq\;
    \Prob_0^\om\big[Y_{t_1+t_2} = \bar e_n^+ \, | \, Y_{t_1}=\bar e_n^+ \big]
    \;=\;
    \Prob^\om_{\bar e_n^+}[Y_{t_2} = \bar e_n^+]
    \;\geq\;
    C_{14},
  \end{align*}
  which for $n \geq \max\{ n_0, n_\delta\}$ contradicts \eqref{eq:hkratio} as $\de < C_{14}$.
\end{proof}

\section{Improving the lower moment condition} \label{sec:imprLB}
It is naive to expect that for general ergodic distributions, one could formulate a necessary condition for a local limit theorem based on the marginal distribution of the conductances only.  In fact, for any edge $\{x, y\} \in E_d$, the lower bound condition
\begin{align}\label{eq:discussion_LB}
  \mean\!\big[ \om(x,y)^{-q} \big]
  \;<\;
  \infty
\end{align}
some $q > 0$ can be improved in several situations to a more general condition based on the heat kernel of the VSRW 
\begin{align}\label{eq:improve1_LB}
  \mean\!\Big[ p^{\om}(t, x, y)^{-q'} \Big]
  \;<\;
  \infty
\end{align}
or on the heat kernel of the CSRW
\begin{align}\label{eq:improve2_LB}
  \mean\!\Big[\big( \mu^{\om}(x)\, q^{\om}(t, x, y)\, \mu^{\om}(y) \big)^{-q'} \Big]
  \;<\;
  \infty,
\end{align}
respectively, for some $q' > 0$.
\begin{prop}
  Suppose that $d \geq 2$ and that Assumptions \ref{ass:environment} and \ref{ass:moment} hold with the moment condition $\mean\!\big[\om(x,y)^{-q}\big] < \infty$ replaced by
  \begin{align*}
    \mean\!\big[ p^{\om}(t, x, y)^{-q}\big]
    \;<\;
    \infty
    \qquad \text{or} \qquad
    \mean\!\Big[
      \big( \mu^{\om}(x)\, q^{\om}(t, x, y)\, \mu^{\om}(x) \big)^{-q}
    \Big]
    \;<\;
    \infty
  \end{align*}
  for any $t>0$ and all $x, y \in \bbZ^d$ such that $\{x, y\} \in E_d$. Then, the statements of Theorem~\ref{thm:ip}  and Theorem~\ref{thm:llt_csrw} hold.
\end{prop}
\begin{proof}
  We only consider the moment condition on the heat kernel $q(t,x,y)$ of the CSRW.  For $x \in \bbZ^d$ and $t > 0$, set $\tilde{\nu}^\om(x) \ldef \sum_{y \sim x} \big(\mu^{\om}(x)\, q^{\om}(t,x,y)\, \mu^{\om}(y)\big)^{-1}\mspace{0mu}$.  Our main objective is to establish a Sobolev inequality similar to the one in \eqref{eq:sob:ineq:2}.  Indeed, one can show along the lines of the proof of this Sobolev inequality in \cite[Proposition 3.5]{ADS13} that
  \begin{align*}
    \Norm{(\eta u)^2}{\rho,B}
    \;\leq\;
    c\, |B|^{\frac{2}{d}}\, \Norm{\tilde{\nu}^{\om}}{q, B}\,
    \frac{1}{|B|}\, \sum_{\{x,y\} \in E}\mspace{-9mu} \mu^{\om}(x)\, q^{\om}(t,x,y)\,
    \mu^{\om}(y) \big( (\eta u)(x) - (\eta u)(y) \big)^2
  \end{align*}
  (cf.\ equation (3.13) in \cite{ADS13}).  Since the sum on the right hand side is smaller than
  \begin{align*}
    \cE^{q_t}(\eta u)
    \;\ldef\;
    \frac{1}{2}\,
    \sum_{x,y}\, \mu^{\om}(x)\, q^{\om}(t,x,y)\,
    \mu^{\om}(y) \big((\eta u)(x) - (\eta u)(y) \big)^2
  \end{align*}
  and $t^{-1} \cE^{q_t}_{\eta^2}(u) \nearrow \cE^{\om}_{\eta^2}(u)$  as $t \searrow 0$, cf.\ \cite[p. 250]{CKS87}, following the arguments in \cite[Proposition 3.5]{ADS13} we conclude that \eqref{eq:sob:ineq:2} holds with $\nu^{\om}$ replaced by $\tilde{\nu}^{\om}$, but with a constant depending now also on $t$.  Finally, since $\tilde{\nu}^{\om}(x) = \tilde{\nu}^{\tau_x \om}(0)$ for every $x \in \bbZ^d$, an application of the ergodic theorem gives that
  \begin{align*}
    \lim_{n \to \infty}\, \Norm{\tilde{\nu}^\om}{q,B(n)}^q
    \;=\;
    \mean\!\big[\tilde{\nu}^\om(0)^q\big]
    \;<\;
    \infty.
  \end{align*}
  The statements of Theorem~\ref{thm:ip}  and Theorem~\ref{thm:llt_csrw} now follow as before.
\end{proof}
Of course, it is not easy to verify \eqref{eq:improve2_LB} (or \eqref{eq:improve1_LB}), because this would require apriori lower bounds on the heat kernel.  However, note that
\begin{align} \label{eq:expansion_q}
  \mu^{\om}(y)\, q^{\om}(t,x,y)
  \;=\;
  \Prob_x^{\om}\big[ Y_t = y \big]
  \;=\;
  \me^{-t}\, \sum_{n=0}^{\infty}\, \frac{t^n}{n!}\, (p^{\om})^{n}(x,y),
\end{align}
where still $p^{\om}(x,y) = \om(x,y) / \mu^{\om}(x)$ and $(p^{\om})^{n+1}(x,y) = \sum_z (p^{\om})^{n}(x,z)\, p^{\om}(z,y)$.  In particular, since $(p^{\om})^{n}(0, x) = 0$ for every $x \in \bbZ^d$ with $|x|_{\infty} = 1$ and even $n$, the summation is only over odd $n$. 

\begin{figure}[t]
  \begin{center}
    \includegraphics{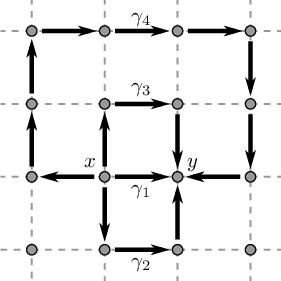}
  \end{center}
  \caption{Illustration of the $4$ different disjoint nearest neighbour paths $\{\ga_i : i = 1, \ldots 4\}$ from $x$ to $y$ for $d = 2$.}
  \label{fig:paths}
\end{figure}

Our first result deals with an example, where \eqref{eq:improve2_LB} and \eqref{eq:discussion_LB} are equivalent, i.e.\ no improvement can be achieved.
\begin{prop} \label{prop:no_impr}
  Let either $\om(x,y) = \th^{\om}(x) \wedge \th^{\om}(y)$ or $\om(x,y) = \th^{\om}(x)\, \th^{\om}(y)$, where $\{\th^{\om}(x) : x \in \bbZ^d\}$ is a family of stationary random variables, i.e.\ $\th^{\om}(x) = \th^{\tau_x \om}(0)$ for all $x \in \bbZ^d$.  Then, for every $\{x, y\} \in E_d$ and all $q \geq 1$,
  \begin{align*} 
    \mean\!\big[ \om(x, y)^{-q} \big]
    \;<\;
    \infty
    \quad \Longleftrightarrow \quad
    \mean\!\Big[\big( \mu^{\om}(x)\, q^{\om}(t, x, y)\, \mu^{\om}(y) \big)^{-q} \Big]
    \;<\;
    \infty.
  \end{align*}
\end{prop}
\begin{proof}
It suffices to show that for every $n$ odd and $\{x, y\} \in E_d$ 
  \begin{align}\label{eq:est_om}
    \mu^{\om}(x) (p^{\om})^{n}(x,y)
    \;\leq\;
    (2d)^n\, \om(x,y).
  \end{align}

  Once the claim \eqref{eq:est_om} is proven, the assertion of the Proposition is immediate.  In order to see this, mind that \eqref{eq:expansion_q} implies the following trivial estimate
  \begin{align*}
    \mu^{\om}(x)\, q^\om(t, x, y)\, \mu^\om(y)
    \;\geq\;
    t\, e^{-t}\, \mu^{\om}(x)\, p^{\om}(x, y)
    \;=\;
    t\, e^{-t}\, \om(x, y).
  \end{align*}
  On the other hand, using again \eqref{eq:expansion_q} and \eqref{eq:est_om} we obtain
  \begin{align*}
    \mu^{\om}(x)\, q^{\om}(t, x, y)\, \mu^{\om}(y)
    \;\leq\;
    \me^{-t}\, \om(x,y)\, \sum_{k \text{ odd}} \frac{t^k}{k!}\, (2d)^k
    \;\leq\;
    \me^{(2d-1)t}\, \om(x, y).
  \end{align*}

  In what follows, we will prove the claim \eqref{eq:est_om}.  For this purpose, we distinguish two cases depending on the definition of the weights $\om$.  To start with, suppose that $\om(x, y) = \th^{\om}(x) \wedge \th^{\om}(x)$ for $\{x, y\} \in E_d$.  Since $\mu^{\om}(x) \leq 2 d\, \th^{\om}(x)$ we obtain that
  \begin{align*}
    \mu^{\om}(x)\, (p^{\om})^{n}(x, y)
    \;\leq\;
    2d\, \th^{\om}(x).
  \end{align*}
  On the other hand, by exploiting the detailed balance condition we get
  \begin{align*}
    \mu^{\om}(x)\, (p^{\om})^{n}(x, y)
    \;=\;
    \mu^{\om}(y)\, (p^{\om})^{n}(y, x)
    \;\leq\;
    2d\, \th^{\om}(y).
  \end{align*}
  Thus, \eqref{eq:est_om} follows by combing these two estimates.

  Let us now consider the case $\om(x, y) = \th^{\om}(x)\, \th^{\om}(y)$.  By introducing the abbreviation $\al^{\om}(x) = \sum_{z \sim x} \th(z)$, we can rewrite
  \begin{align}
    \mu^{\om}(x)
    \;=\;
    \th^{\om}(x)\, \al^{\om}(x)
    \qquad \text{and} \qquad
    p^{\om}(x, y)
    \;=\;
    \frac{\th^{\om}(y)}{\al^{\om}(x)}.
  \end{align}
  Thus, for $n = 3$,  we have that
  \begin{align*}
    \mu^{\om}(x)\, (p^{\om})^{3}(x, y)
    \;=\;
    \th^{\om}(x)\, \th^{\om}(y)\,
    \sum_{z_1 \sim x}\, \sum_{z_2 \sim y}\,
    \frac{\th^{\om}(z_1)}{\al^{\om}(z_2)}\, \frac{\th^{\om}(z_2)}{\al^{\om}(z_1)}
    \;\leq\;
    (2 d)^2\, \om(x,y),
  \end{align*}
  where we used that $\th^{\om}(z_1) / \al^{\om}(z_2) \leq 1$ and $\th^{\om}(z_2) / \al^{\om}(z_1) \leq 1$ since $z_1 \sim z_2$.  This implies \eqref{eq:est_om} for $n = 3$.  For general odd $n$ the claim \eqref{eq:est_om} follows inductively.
\end{proof}
Note that in Proposition~\ref{prop:no_impr} no assumptions have been made on the distribution of the random variables $\big\{\th^{\om}(x) : x \in \bbZ^d\big\}$.  In general the computations can be rather intricate, in particular due to the quotient in the definition of $p^\om(x,y)$. For simplicity we will now focus on the case with conductances bounded from above, that is $p = \infty$.
\begin{prop}\label{prop:impr1}
  Assume that the conductances $\{\om(e) : e \in E_d\}$ are independent and uniformly bounded from above, i.e.\ $\prob$-a.s.\ $\om(e) \leq 1$.   Then, if \eqref{eq:discussion_LB} holds for any $q > 0$, we have that \eqref{eq:improve2_LB} holds for all $q' < 2 d q$.  In particular, if $\mean\big[ \om(x,y)^{-q} \big] < \infty$ for some $q > 1/4$ and all $\{x, y\} \in E_d$ then there exists $q' > \frac{d}{2}$ such that
  \begin{align*}
    \mean\!%
    \Big[\big( \mu^{\om}(x)\, q^{\om}(t, x, y)\, \mu^{\om}(y) \big)^{-q'} \Big] 
    \;<\;
    \infty.
  \end{align*}
\end{prop}
For the proof we will need the following simple lemma.
\begin{lemma} \label{lem:momZ}
  Let $Z_1, Z_2, \ldots, Z_N$ be independent random variables such that
  \begin{align*}
    \prob\!\big[0 < Z_i \leq 1\big] \;=\; 1
    \qquad \text{and} \qquad
    c_i(\al)
    \;=\;
    \mean\!\big[ Z_i^{-\al}\big]
    \;<\;
    \infty, 
  \end{align*}
  for some $\al > 0$.  Then
  \begin{align*}
    \mean\!\Big[\big( Z_1 + Z_2 + \ldots + Z_N \big)^{-\be}\Big]
    \;<\;
    \infty
  \end{align*}
  for all $\be < N \al$.
\end{lemma}
\begin{proof}
  Note that $Z_1 + Z_2 + \ldots + Z_N \geq \max_{i = 1, \ldots, N} Z_i$, and \v{C}eby\v{s}ev's inequality yields
  \begin{align*}
    \prob\!\Big[\Big( \max_{i = 1, \ldots, N} Z_i \Big)^{\!-\be} >\, t \Big]
    &\;=\;
    \prob\!\Big[\max_{i=1,...,N} Z_i \,<\, t^{-\frac{1}{\be}}\Big]
    \;\leq\;
    t^{-\frac{n \al}{\be}}\, \prod_{i=1}^n c_i(\al).
  \end{align*}
  Since $\prob[0 < Z_i \leq 1] = 1$, the assertion follows by integrating the estimate above in $t$ over the interval $[1, \infty)$.
\end{proof}
\begin{proof}[Proof of Proposition~\ref{prop:impr1}]
  Suppose that \eqref{eq:discussion_LB} holds for any $q > 0$.  Since the conductances are assumed to be bounded from above $\om(x,y) \leq 1$ we have $\mu^{\om}(x) \leq 2d$ and therefore $(p^\om)^n(x,y) \geq (2d)^{-n} \om^{(n)}(x,y)$, where we define $\om^{1}(x, y) \ldef \om(x, y)$ and $\om^{n+1}(x,y) = \sum_z \om^{n}(x,z)\, \om(z,y)$.  Thus, in view of \eqref{eq:expansion_q} we get
  \begin{align*}
    \mu^{\om}(x)\, q^{\om}(t, x, y)\, \mu^{\om}(y)
    \;\geq\;
    c\, \big( \om(x, y) \,+\, \om^{3}(x, y) \,+\, \om^{9}(x, y) \big),
  \end{align*}
  for some positive constant $c \equiv c(d, t)$. Hence, it suffices to show that
  \begin{align} \label{eq:mom_om3}
    \mean\!%
    \Big[
      \big( \om(x, y) \,+\, \om^{3}(x, y) \,+\, \om^{9}(x, y) \big)^{-q'}
    \Big]
    \;<\;
    \infty.
  \end{align}
  Note that for $\{x, y\} \in E_d$ there are at least $2d$ different disjoint nearest neighbour oriented paths $\{\ga_i : i = 1, \ldots, 2d\}$ from $x$ to say $y$, meaning that the sets of edges along the paths $(\gamma_i)$ are disjoint: one direct path of length $1$, $2(d-1)$ paths of length $3$ and one path of length $9$ (see Figure~\ref{fig:paths}). That is,
  \begin{align*}
    \om(x, y) + \om^{3}(x, y) + \om^{9}(x, y)
    \;\geq\;
    Z_1 + Z_2 + \ldots + Z_{2d},
  \end{align*}
  where $Z_i = \prod_{e \in \ga_i} \om(e)$ for $i = 1, \ldots, 2d$.  Using the independence we see that $\mean\big[Z_i^{-q}\big] < \infty$ for $i = 1, \ldots, 2d$. This implies \eqref{eq:mom_om3} by Lemma~\ref{lem:momZ}.
\end{proof}
\begin{prop}\label{prop:impr2}
  Let either $\om(x,y) = \th^{\om}(x) \vee \th^{\om}(y)$ or $\om(x,y) = \th^{\om}(x) + \th^{\om}(y)$, where $\{\th^{\om}(x) : x \in \bbZ^d\}$ is a family of uniformly bounded i.d.d.\ random variables, that is $p = \infty$. Then, if \eqref{eq:discussion_LB} holds for any $q > 0$, we have that \eqref{eq:improve2_LB} holds for all $q' < \frac{2}{3} (d+1) q$.  In particular, if $\mean\big[ \om(x,y)^{-q} \big] < \infty$ for some $q > \frac 3 4\, \frac{d}{d+1}$ and all $\{x, y\} \in E_d$ then there exists $q' > \frac{d}{2}$ such that
  \begin{align*}
    \mean\!%
    \Big[\big( \mu^{\om}(x)\, q^{\om}(t, x, y)\, \mu^{\om}(y) \big)^{-q'} \Big] 
    \;<\;
    \infty.
  \end{align*}
\end{prop}
\begin{lemma} \label{lem:momX}
  Let $X_1$ and $X_2$ be i.i.d.\ random variables such that
  \begin{align*}
    \prob\!\big[0 < X_i \leq 1\big] \;=\; 1,
    \qquad \text{and} \qquad
    c(\al)
    \;=\;
    \mean\!\big[ X_i^{-\al}\big]
    \;<\;
    \infty, \qquad i=1,2,
  \end{align*}
  for some $\al > 0$.  Then
  \begin{align*}
    \mean\!\Big[ \big(X_1^2 X_2 \vee X_1 X_2^2 \big)^{-\be} \Big]
    \;<\;
    \infty
  \end{align*}
  for all $\be < \frac{2}{3} \al$.
\end{lemma}
\begin{proof}
  First of all, note that for any $s>0$ we have
  \begin{align*}
    \prob\!\big[X_1^2 X_2 \vee X_1 X_2^2 < s\big]
    &\;=\;
    \prob\!\Big[X_1< s^{\frac{1}{3}},\, X_2 < s^{\frac{1}{3}} \Big]
    \\
    &\mspace{36mu}+\;
    \prob\!\Big[
      X_1\geq  s^{\frac{1}{3}},\, X_2 < s^{\frac{1}{3}},\, X_2< s X_1^{-2}
    \Big]
    \\
    &\mspace{36mu}+\;
    \prob\!\Big[
      X_1<  s^{\frac{1}{3}},\, X_2 \geq s^{\frac{1}{3}},\, X_1< s X_2^{-2}
    \Big]
    \\
    &\;=\;
    \prob\!\Big[X_1< s^{\frac{1}{3}} \Big]\, \prob\!\Big[X_2 < s^{\frac{1}{3}}\Big]
    \,+\, 2\, \prob\!\Big[X_1\geq  s^{ \frac{1}{3}}, X_2< s X_1^{-2} \Big],
  \end{align*}
  where the second equality follows by a symmetry argument.  Let us denote by $\widehat{\mean}_{\al}$ the expectation w.r.t.\ $\md\widehat{\prob}_{\al} \ldef c_1(\al)^{-1} X_1^{-\al}\, \md\prob$.  Then, we obtain
  \begin{align*}
    \mean\!\Big[ \indicator_{\{X_1 \geq  s^{1/3} \}}\, X_1^{-2 \al} \Big]
    \;=\;
    c(\al)\; \widehat{\mean}_{\al}\!%
    \Big[ \indicator_{\{ X_1 \geq  s^{1/3} \}}\, X_1^{-\al} \Big]
    \;\leq\;
    c(\al)\; s^{\frac{1}{3} \al}.
  \end{align*}
  Thus, \v{C}eby\v{s}ev's inequality yields the following upper bound
  \begin{align*}
    \prob\!\Big[X_1\geq  s^{\frac 1 3},\, X_2< s X_1^{-2} \Big] 
    \;=\;
    \mean\!%
    \Big[
      \indicator_{\{X_1 \geq  s^{1/3}\}}\,
      \prob\!\big[ X_2^{-\al} > s^{-\al} X_1^{2 \al}  \,\big|\, X_1 \big]
    \Big]
    &\;\leq\;
    c(\al)^2\, s^{\frac{2}{3} \al}.
  \end{align*}
  Hence, for any $t > 0$, again by \v{C}eby\v{s}ev's inequality
  \begin{align*} 
    \prob\!\Big[ \big(X_1^2 X_2 \vee X_1 X_2^2 \big)^{-\beta} > t \Big]
    &\;=\;
    \prob\!\Big[ X_1 < t^{-\frac{1}{3 \be}} \Big]\,
    \prob\!\Big[ X_2 < t^{-\frac{1}{3 \be}} \Big]
    \\
    &\mspace{36mu}+\;
    2\,
    \prob\!\Big[
      X_1 \geq  t^{-\frac{1}{3\be}},\, X_2 < t^{-\frac{1}{\be}}\, X_1^{-2}
    \Big]
    \;\leq\;
    3\, c(\al)^2\, t^{-\frac{2 \al}{3 \be}}\!.
  \end{align*}
\end{proof}
\begin{proof}[Proof of Proposition~\ref{prop:impr2}]
  Since
  \begin{align*}
    \frac{1}{2}\, \Big( \th^{\om}(x) + \th^{\om}(y) \Big)
    \;\leq\;
    \th^{\om}(x) \vee \th^{\om}(y)
    \;\leq\;
    \th^{\om}(x) + \th^{\om}(y)
  \end{align*}
  it is enough to consider the case $\om(x,y) = \th^{\om}(x) \vee \th^{\om}(y)$.  Again it suffices to show \eqref{eq:mom_om3}.  First note that we have $\prob\!\big[ \th^{\om}(x) < s \big]^2 = \prob\!\big[\om(x,y) < s \big]$ for any edge $\{x, y\} \in E_d$ and therefore
  \begin{align*}
    \prob\!\big[\th^{\om}(x)^{-r} > t \big]
    \;=\;
    \prob\!\Big[ \th^{\om}(x) < t^{-\frac{1}{r}} \Big]
    \;=\;
    \prob\!\Big[ \om(x,y)^{-q} > t^{\frac{q}{r}} \Big]^{\frac{1}{2}}
    \;\leq\;
    t^{-\frac{q}{2r}}\, \mean\!\big[ \om(x,y)^{-q} \big]^{\frac{1}{2}}.
  \end{align*}
  In particular, $\mean\!\big[ \th^{\om}(x)^{-r} \big] < \infty$ for all $r < q/2$.  Let $x = x_0, x_1, \ldots, x_{2k}, x_{2k+1} = y$ be a disjoint nearest neighbour path from $x$ to $y$. Then,
  \begin{align} \label{eq:low_est_path}
    \prod_{j=0}^{2k} \om(x_j,x_{j+1})
    &\;=\;
    \prod_{j=0}^{2k} \th^{\om}(x_j) \vee \th^{\om}(x_{j+1})
    \nonumber\\
    &\;\geq\;
    \prod_{j=1}^{k} \th^{\om}(x_j) \cdot 
    \big( \th^{\om}(x_k) \vee \th^{\om}(x_{k+1}) \big) \cdot\!\!
    \prod_{j = k+1}^{2k} \th^{\om}(x_j).
  \end{align}
  In particular, note that the right hand side does not depend on $\th^{\om}(x)$ and $\th^{\om}(y)$.  Moreover, by Lemma~\ref{lem:momX}
  \begin{align*}
    \mean\!%
    \bigg[
      \Big( \prod\nolimits_{j=0}^{2k} \om(x_j,x_{j+1}) \Big)^{\!-\be}
    \bigg]
    \;<\;
    \infty,
    \qquad \forall\, \be < \frac{q}{3}.
  \end{align*}
  In order to prove \eqref{eq:mom_om3} we consider again $2d$ disjoint nearest neighbour paths from $x$ to $y$ and define $Z_i$, $i = 1, \ldots, 2d$ as above in the proof of Proposition~\ref{prop:impr1}.  Then, using \eqref{eq:low_est_path} we have $Z_i \geq \tilde{Z}_i$ for $i=2,\ldots ,2d$, where $\{\tilde{Z}_i : i = 2, \ldots, d \}$ is a collection of independent random variables with $\mean\!\big[\tilde{Z}_i^{-\be} \big] < \infty$ for all $\be < \frac{q}{3}$.  Furthermore, $\{\tilde{Z}_i \}$ is independent of $Z_1 = \om(x, y)$.  Thus,
  \begin{align*}
    \om(x, y) + \om^{3}(x, y) + \om^{9}(x, y)
    \;\geq\;
    Z_1 + \tilde{Z}_2 + \ldots + \tilde{Z}_{2d}
  \end{align*}
  and by Lemma~\ref{lem:momZ} we conclude that \eqref{eq:mom_om3} holds for all $q'< q + (2d-1) \be$, which implies the claim.
\end{proof}

\appendix
\section{Technical estimates}
In this section we collect some technical estimates needed in the proofs.  In a sense some of them may be seen as a replacement for a discrete chain rule.
\begin{lemma}
  \begin{enumerate}[ (i) ]
    \item For all $a, b \geq 0$ and any $\al, \be \ne 0$
      \begin{align}\label{eq:A1:chain:ub1}
        \big| a^\al -  b^\al \big|
        \;\leq\;
        \Big|\frac{\al}{\be}\Big|\;
        \big| a^\be - b^\be \big|\,
        \big(\,|a|^{\al-\be} + |b|^{\al-\be} \big).
      \end{align}
    \item For all $a,b \geq 0$ and any $\al > 1/2 $
      \begin{align}\label{eq:A1:pol:ub}
        \big(a^{\al} - b^{\al}\big)^2
        \;\leq\;
        \bigg|\frac{\al^2}{2\al-1}\bigg|\, \big(a - b\big)\,
        \big(a^{2\al-1} - b^{2\al-1}\big).
      \end{align}
      Moreover, if $a,b > 0$ then
      \begin{align}\label{eq:A1:log:ub}
        \big(\ln a - \ln b\big)^2
        \;\leq\;
        -\big(a^{-1} - b^{-1}\big)\, \big(a - b\big).
      \end{align} 
    \item For all $a,b \geq 0$ and any $\al,\be \geq 0$
      \begin{align}\label{eq:A1:chain:lo}
        \big(a^{\al} + b^{\al}\big)\,
        \big|a^{\be} -\, b^{\be}\big|
        \;\leq\;
        2\, \big|a^{\al+\be} -\, b^{\al+\be}\big|.
      \end{align}
    \item For all $a,b \geq 0$ and any $\al \geq 1/2$
      \begin{align}\label{eq:A1:chain:ub2}
        \big(a^{2\al-1} + b^{2\al-1}\big)\, \big|a - b\big|
        \;\leq\; 
        4\,
        \big| a^{\al} -  b^{\al}\big|\,
        \big(a^{\al} + b^{\al}\big).
      \end{align}
    \item For all $a,b\geq 0$ and $\al \in (0,1/2)$,
      \begin{align} \label{eq:chain_super}
        \big(a^{2\al-1}-b^{2\al -1}\big) (a-b)
        \;\leq\;
        \frac{2\al -1}{\al^2} \Big( a^{\al}-b^{\al}\Big)^2.
      \end{align}
  \end{enumerate}
\end{lemma}
\begin{proof}
  \textit{(i)} W.l.o.g.\ assume that $a \geq b$.  Then, for all $\al, \be \ne 0$,
  \begin{align*}
    \big|a^{\al} - b^{\al}\big|
    \;=\;
    \bigg|
      \al\, \int_{b}^{a} t^{\be-1}\, t^{\al-\be}\, \md t
    \bigg|
    \;\leq\;
    \Big|\frac{\al}{\be}\Big|\, \bigg(\max_{t \in [b, a]} t^{\al-\be}\bigg)
    \big| a^{\be} - b^{\be} \big|.
  \end{align*}
  Since $t \longmapsto t^{\al-\be}$ for $t \geq 0$ is monotone decreasing if $\al -\be < 0$ and monotone increasing if $\al-\be > 0$, the maximum is attained at one of the boundary points.  In particular, $\max_{t \in [b, a]} t^{\al-\be} \leq |a|^{\al-\be} + |b|^{\al-\be}$.

  \textit{(ii)} Notice that for all $a,b \geq 0$ and for any $\al > 0$,
  \begin{align*}
    \big(a^\al - b^\al\big)^2
    \;=\;
    \bigg(\al\, \int_b^a t^{\al-1} \, \md t\bigg)^{\!\!2}.
  \end{align*}
  Hence, for any $\al > 1/2$, an application of the Cauchy-Schwarz inequality yields
  \begin{align*}
    \bigg( \al \int_{b}^{a} t^{\al-1} \md t\bigg)^{\!\!2}
    \;\leq\;
    \al^2\, \int_{b}^{a} \md t\, \int_{b}^{a} t^{2\al-2} \md t
    \;=\;
    \frac{\al^2}{2\al-1}\, \big(a - b\big)\,
    \big(a^{2\al-1} - b^{2\al-1}\big).
  \end{align*}
  Similarly, for $a, b > 0$ we have
  \begin{align*}
    (\ln a -\ln b)^2
    \;=\;
    \left( \int_b^a t^{-1}\, \md t \right)^{\!\!2}
    \;\leq\;
    \int_b^a \md t\, \int_b^a t^{-2}\, \md t
    \;=\;
    -(a-b)\, \big( a^{-1} - b^{-1} \big).
\end{align*}

  \textit{(iii)}  Since the assertion is trivial, if $a = 0$, we assume that $a \ne 0$ and set $z = b/a \geq 0$.  Then, the left-hand side of \eqref{eq:A1:chain:lo} reads $\big(1 + z^{\al}\big)\, \big|1 - z^{\be}\big|$.  Provided that $\al \geq 0$, by distinguishing two cases, $z \in [0,1)$ or $z \geq 1$, we obtain
  \begin{align*}
    \big(1 + z^{\al}\big)\, \big|1 - z^{\be}\big|
    \;=\;
    2 \big|1 - z^{\al+\be}\big| \,-\, \big|1 - z^{\al}\big|\,
    \big(1 + z^{\be}\big)
    \;\leq\;
    2 \big|1 - z^{\al+\be}\big|.
  \end{align*}
  This completes the proof.

  \textit{(iv)} The assertion follows immediately from \eqref{eq:A1:chain:ub1} and \eqref{eq:A1:chain:lo}.

  \textit{(v)} W.l.o.g.\ assume that $a \geq b$.  Then, by the Cauchy-Schwarz inequality
  \begin{align*}
    \big( a^{\al} - b^{\al} \big)^2
    &\;=\;
    \al^2\,  \bigg( \int_b^a t^{\al -1}\, \md t \bigg)^{\!\!2}
    \;\leq\;
    \al^2\, (a - b)\, \int_b^a t^{2\al-2} \, \md t
    \\[.5ex]
    &\;=\;
    \frac{\al^2}{2\al-1}\, (a-b) \big( a^{2\al-1} - b^{2\al-1} \big),
  \end{align*}
  which gives the claim since $2\al - 1 < 0$.
\end{proof}
\begin{lemma}
  For all $x, y > 0$, $b \geq  a \geq 0$ and $\be \in (0, 1]$
  \begin{align}\label{eq:A2:beta}
    \bigg( \dfrac{b^2}{y^{\be}} - \dfrac{a^2}{x^{\be}} \bigg)\, \big( y - x \big)
    \;\leq\;
    \begin{cases}
      \dfrac{a^2}{2}\,
      \bigg(\dfrac{1}{y^{\be}} - \dfrac{1}{x^{\be}}\bigg)\,
      \big(y - x \big)
      \,+\,
      \dfrac{8}{\be}\, \dfrac{b^2}{a^2}\, \big(b - a\big)^2\, y^{1-\be},
      & a > 0,\\[2ex]
      b^2\, y^{1-\be},
      & a = 0.
    \end{cases}
  \end{align}
\end{lemma}
\begin{proof}
  First of all, the case $a = 0$ is trivial.  Further, notice that \eqref{eq:A2:beta} holds true if $x = y$ or $a = b$ where we use in the latter that $\big(y^{-\be} - x^{-\be}\big)\, (y - x) \leq 0$.  Therefore, we assume in the sequel that $x \neq y$ and $0 < a < b$.  Since
  \begin{align*}
    \bigg(\frac{b^2}{y^\be} - \frac{a^2}{x^\be}\bigg)\, \big(y - x\big)
    \;=\;
    a^2\, \bigg(\frac{1}{y^\be} - \frac{1}{x^\be}\bigg)\, \big(y - x\big)
    \,+\,
    \frac{1}{y^\be}\, \big(b^2 - a^2\big)\, \big(y - x\big),
  \end{align*}
  the assertion \eqref{eq:A2:beta} is immediate without the need of the second positive term, if $y^{-\be}\, \big(b^2 - a^2\big)\, (y-x) \leq - a^2 \big(y^{-\be} - x^{-\be}\big) \big(y - x\big) / 2$.  Thus, we are left with the case
  \begin{align}\label{eq:A2:case2}
    \frac{1}{y^{\be}}\,\big(b^2 - a^2\big)\, \big(y - x\big)
    \;\geq\;
    -\frac{a^2}{2}\, \bigg( \frac{1}{y^{\be}} - \frac{1}{x^{\be}}\bigg)\,
    \big(y - x\big).
  \end{align}
  Since the right-hand side of \eqref{eq:A2:case2} is positive and $a < b$, we deduce from \eqref{eq:A2:case2} that $y > x$.  Thus, an elementary computation yields that \eqref{eq:A2:case2} is equivalent to
  \begin{align*}
    \frac{a^2}{2}\, \frac{1}{x^{\be}}\, \Big( y^{\be} - x^{\be} \Big)
    \;\leq\;
    b^2 - a^2.
  \end{align*}
  Hence, by exploiting the fact that $x < y$, we obtain
  \begin{align*}
     \frac{1}{y^{\be}}\, \big(b^2 - a^2\big)\, \big(y - x\big)
     \;\leq\;
     \frac{2}{a^2}\, \big( b^2 - a^2 \big)^2\, y^{1-\be}\,
     \frac{1 - x/y}{(x/y)^{-\be} - 1}
     \;\leq\; \frac{8}{\be}\, \frac{b^2}{a^2}\, (b-a)^2\, y^{1-\be},
  \end{align*}
  where we used that $\max_{z \in [0,1]} (1-z)/(z^{-\be} - 1) = 1/\be$.  Thus, together with the fact that $(y^{-\be} - x^{-\be})(y-x) \leq 0$, the assertion \eqref{eq:A2:beta} follows.
\end{proof}

\subsubsection*{Acknowledgment}
We thank N.~Berger for the example in Theorem \ref{thm:ex_lclt} and M.~Barlow, A.~Chiarini, T.~Kumagai and P.~Mathieu for useful discussions and valuable comments.  M.S. gratefully acknowledge financial support from the DFG Forschergruppe 718 ''Analysis and Stochastics in Complex Physical Systems''.
\bibliographystyle{plain}
\bibliography{literature}

\end{document}